\documentclass[12pt]{article}

\newcommand{\blind}{1}

\usepackage{tikz}
\usepackage{adjustbox}

\usepackage[round]{natbib}

\usepackage{amsmath, amsfonts, amssymb, amsthm, dsfont}

\usepackage{mathrsfs}

\usepackage{soul}

\usepackage{bm}

\usepackage{booktabs} 

\definecolor{mypurple}{RGB}{0,0,0}
\definecolor{myblue}{RGB}{0,87,120}
\definecolor{aqua}{RGB}{0,0,0}
\definecolor{myorange}{RGB}{0,0,0}
\definecolor{mygrey}{RGB}{255,255,255}

\usepackage{ulem}

\usepackage{graphicx}
\usepackage{color}
\usepackage{nicefrac}

\usepackage{scrextend}

\usepackage{titlesec}
\titleformat{\section}[hang]{\large\center\scshape}{\thesection.}{1em}{}
\titleformat{\subsection}[hang]{\large}{\thesubsection.}{1em}{}
\titleformat{\subsubsection}[hang]{}{\thesubsubsection.}{1em}{}

\usepackage{xr-hyper}
\usepackage[colorlinks=true, linkcolor=red, urlcolor=blue, citecolor=blue]{hyperref}

\makeatletter
\newcommand*{\addFileDependency}[1]{
  \typeout{(#1)}
  \@addtofilelist{#1}
  \IfFileExists{#1}{}{\typeout{No file #1.}}
}
\makeatother
 
\newcommand*{\myexternaldocument}[1]{%
    \externaldocument{#1}%
    \addFileDependency{#1.tex}%
    \addFileDependency{#1.aux}%
}

\myexternaldocument{main_jasa_app}

\newtheoremstyle{mytheoremstyle} 
    {0.3cm}                      
    {0cm}                        
    {\itshape}                   
    {}                           
    {\scshape}                   
    {: }                          
    {0em}                       
    {}  

\theoremstyle{mytheoremstyle}
\newtheorem{Theorem}{Theorem}
\newtheorem{Definition}{Definition}
\newtheorem{Lemma}{Lemma}
\newtheorem{Corollary}{Corollary}
\newtheorem{Proposition}{Proposition}

\renewenvironment{proof}{{\noindent \sc Proof:}}{\qed}

\newtheoremstyle{myExampleRemarkstyle} 
    {0.3cm}                    
    {0cm}                           
    {\itshape}                   
    {}                           
    {\scshape}                   
    {: }                          
    {0em}                       
    {}  

\theoremstyle{myExampleRemarkstyle}
 
\newtheorem{Remark}{Remark}
\newtheorem{Assumption}{Assumption}

\renewcommand{\theAssumption}{\Alph{Assumption}}

\newtheorem{innercustomgeneric}{\customgenericname}
\providecommand{\customgenericname}{}
\newcommand{\newcustomtheorem}[2]{%
  \newenvironment{#1}[1]
  {%
   \renewcommand\customgenericname{#2}%
   \renewcommand\theinnercustomgeneric{##1}%
   \innercustomgeneric
  }
  {\endinnercustomgeneric}
}

\newcustomtheorem{customthm}{Theorem}
\newcustomtheorem{customlemma}{Lemma}

\newtheoremstyle{simuStyle}
{0.3cm} 
{0cm} 
{} 
{} 
{\bfseries} 
{.} 
{0em} 
{} 

\theoremstyle{simuStyle}

\newtheoremstyle{stratStyle}
{0.3cm} 
{0cm} 
{} 
{} 
{\scshape} 
{: } 
{0em} 
{} 

\theoremstyle{stratStyle}

\DeclareSymbolFont{lettersA}{U}{txmia}{m}{it}
\DeclareMathSymbol{\real}{\mathord}{lettersA}{"92}
\DeclareMathSymbol{\field}{\mathord}{lettersA}{"83}

\def\real{{\rm I\!R}}

\DeclareMathOperator*{\argmin}{argmin}

\def\0{{\bf 0}}

\def\btheta{{\bm{\theta}}}
\def\hbtheta{\hat{\bm{\theta}}}
\def\bbeta{{\bm{\beta}}}

\def\hbpi{\hat{{\bm{\pi}}}}

\def\bOmega{{\bm{\Omega}}}

\def\X{{\bf X}}

\def\Y{\mathbf{Y}}

\DeclareMathOperator*{\var}{var}

\DeclareMathOperator*{\argzero}{argzero}

\def\ph{\hat{\bm{\pi}}}
\def\bt{\bm{\theta}}
\def\bto{\bt_0}
\def\bT{\bm{\Theta}}
\def\bD{\bm{\Delta}}
\def\hbt{\hat{\bt}_{(n,H)}}
\def\hbtseq{\left\{\hbt^{(k)}\right\}_{k\in\mathbb{N}}}
\def\cn{\mathbf{c}(n)}
\def\bw{\bm{\omega}}
\def\bwo{\bw_0}

\def\Ln{\mathbf{L}(n)}

\def\N{\mathbb{N}}
\def\Ns{\N^\ast}

\newcommand{\hp}[2]{\ph\left({#1},n,{#2}\right)}
\newcommand{\hpH}[1]{\frac{1}{H}\sum_{h=1}^H\ph\left({#1},n,\bw_{h}\right)}
\renewcommand{\a}[1]{\mathbf{a}\left({#1}\right)}
\renewcommand{\r}[1]{\mathbf{r}\left({#1},n\right)}
\renewcommand{\r}[1]{\mathbf{r}\left({#1},n\right)}
\renewcommand{\v}[2]{\mathbf{v}\left({#1},n,{#2}\right)}
\newcommand{\vH}[1]{\frac{1}{H}\sum_{h=1}^H\mathbf{v}\left({#1},n,\bw_{h}\right)}

\begin{document}

\def\spacingset#1{\renewcommand{\baselinestretch}%
 {#1}\small\normalsize} \spacingset{1}

\let\refBKP\ref
\renewcommand{\ref}[1]{{\upshape\refBKP{#1}}}

\if1\blind
{
    \title{\bf Phase~Transition~Unbiased Estimation in High~Dimensional Settings}
    \author{St\'ephane Guerrier,
        Mucyo Karemera,\\
        Samuel Orso  \& 
        Maria-Pia Victoria-Feser \\
        \vspace{0.01cm}\\
        Research Center for Statistics, University of Geneva
    }\date{}
    \maketitle
} \fi

\if0\blind
{
    \bigskip
    \bigskip
    \bigskip
    \begin{center}
        {\LARGE\bf Phase~Transition~Unbiased Estimation \\ \vspace{.2cm}
        in High~Dimensional Settings}
    \end{center}
    \medskip
} \fi


\begin{abstract}
     An important challenge in statistical analysis concerns the control of the finite sample bias of estimators. This problem is magnified in high dimensional settings where the number of variables $p$ diverge with the sample size $n$. However, it is difficult to establish whether an estimator $\hat{\bm{\theta}}$ of $\bm{\theta}_0$ is unbiased and the asymptotic order of $\mathbb{E}[\hat{\bm{\theta}}]-\bm{\theta}_0$ is commonly used instead. We introduce a new property to assess the bias, called phase transition unbiasedness, which is weaker than unbiasedness but stronger than asymptotic results. An estimator satisfying this property is such that $\big\lVert\mathbb{E}[\hbtheta]-\bto\big\rVert_2=0$, for all $n$ greater than a finite sample size $n^\ast$. We propose a phase transition unbiased estimator by matching an initial estimator computed on the sample and on simulated data. It is computed using an algorithm which is shown to converge exponentially fast. The initial estimator is not required to be consistent and thus may be conveniently chosen for computational efficiency or for other properties. We demonstrate the consistency and the limiting distribution of the estimator in high dimension. Finally, we develop new estimators for logistic regression models, with and without random effects, that enjoy additional properties such as robustness to data contamination and to the problem of separability.
\end{abstract}

\noindent%
\vfill


\newpage
\spacingset{1.5} 

\section{Introduction}

An important challenge in statistical analysis concerns the control of the (finite sample) bias of estimators. For example, the Maximum Likelihood Estimator (MLE) has a bias that can result in a significant inferential loss. Indeed, under usual regularity conditions, the MLE is consistent and has a bias of asymptotic order $\mathcal{O}(n^{-1})$, where $n$ denotes the sample size. More generally, we qualify an estimator as \textit{asymptotically unbiased of order $\alpha$} if it has a bias of asymptotic order $\mathcal{O}(n^{-\alpha})$, elementwise, where $\alpha > 0$. Thus, the MLE is typically asymptotically unbiased of order 1 and its bias vanishes as $n$ diverges. However, when $n$ is finite, the bias of the MLE may be large, and this bias is transferred to quantities such as test statistics where the MLE is used as a plug-in value (see e.g.~\citealp{kosmidis2014bias} for a review). This problem is typically magnified in situations where the number of variables $p$ is large and possibly allowed to increase with $n$. For example, \cite{sur2019modern} show that the MLE for the logistic regression model can be severely biased in the case where $n$ and $p$ become increasingly large (given a fixed ratio).

The bias of an estimator can often be reduced by using bias correction methods that have recently received substantive attention. Bias correction is often achieved by simulation methods such as the jackknife \citep{Efro:82} or the bootstrap \citep{Efro:79}, or can be approximated using asymptotic expansions (see e.g. \citealp{CoVa:97}, \citealp{CoTU:08} and the references therein
). For example, the bootstrap bias correction \citep{efron1994introduction} allows to obtain an asymptotically unbiased estimator of order 2, provided that the initial estimator is consistent and asymptotically unbiased of order 1 \citep[see][]{hall1988bootstrap}. An alternative approach is to correct the bias of an estimator by modifying its associated estimating equations. For example, \cite{Firt:93} provides an adjustment of the MLE score function and this approach has been successively adapted and extended to Generalized Linear Models (GLM), among others, by \cite{MeMa:95,BuMaGr:02,KoFi:09,KoFi:11,Kosm:14}
. Several other bias correction methods have been proposed
and a review on existing methods can, for example, be found in \cite{kosmidis2014bias}.

It is in general difficult to establish whether an estimator (possibly resulting from a bias correction technique) is \textit{unbiased} for all $n$ (assuming that $n$ is large enough so that the estimator can be computed), and therefore the \textit{asymptotic} order of the bias is used to ``quantify'' its magnitude. However, the asymptotic order of the bias has mainly been studied in low-dimensional settings and its extension in asymptotic regimes where both $n$ and $p$ can tend to infinity is often unclear. To address this, we introduce a new property to quantify the bias of an estimator in possibly high dimensional settings called \textit{Phase Transition (PT) unbiasedness}. This property, presented in Definition \ref{pt-unbiased} (below), can serve as a middle-ground between unbiasedness and the asymptotic order of the bias.

\begin{Definition}[Phase Transition unbiasedness]
\label{pt-unbiased}
An estimator $\hat{\bm{\theta}}$ of $\bm{\theta}_0$ is said to be PT-unbiased if there exists a $n^\ast \in \mathbb{N}^* \equiv \mathbb{N} \backslash\left\{0\right\}$ such that for all $n \in \mathbb{N}^*$ with $n \geq n^*$, we have $\big\lVert	\mathbb{E} [\hbtheta] - \bto  \big\rVert_2 = 0$.
\end{Definition}

In short, this property implies that if an estimator is PT-unbiased, we have $\big\lVert	\mathbb{E} [\hbtheta] - \bto  \big\rVert_2 = 0$ for all $n$ greater than a finite sample size $n^\ast$. Among others, this interpretation also explains the name of this property since, starting from a certain $n^\ast$, the estimator transitions from a biased phase to an unbiased phase. Therefore, PT-unbiasedness is weaker than the classical notion of unbiasedness but stronger than the asymptotic unbiasedness of order $\alpha$. Indeed, suppose that $\tilde{\bm{\theta}}$ is an asymptotically unbiased estimator of order $\alpha$, then we can write $ \mathbb{E} [\tilde{\bm{\theta}}] = \bto + \mathcal{O}\left(n^{-\alpha}\right)$, elementwise, for a \textit{specific} $\alpha > 0$. If $\hat{\bm{\theta}}$ is a PT-unbiased estimator, we have $\mathbb{E} [\hat{\bm{\theta}}] = \bto + \mathcal{O}\left(n^{-\beta}\right)$, elementwise, for \textit{all} $\beta > 0$ and thus for all $\beta \geq \alpha$. 

In this article, we propose a class of PT-unbiased estimators, called \textit{IB-estimators}, building on the ideas of the Iterative Bootstrap (IB)
\citep{kuk1995asymptotically} and of indirect inference \citep{smith1993estimating,gourieroux1993indirect}. In short, these simulation-based methods allow to obtain consistent estimators from inconsistent \textit{initial} ones. Moreover, \cite{guerrier2018simulation} study the order of the asymptotic bias of these estimators when the initial estimator is consistent. In particular, in low-dimensional settings, the latter work shows that these techniques present some strong similarities and can lead to asymptotically unbiased estimators of order $\delta$ where $\delta \in (2, 3]$. In Section \ref{sec:main}, we further study these simulation-based techniques in high dimensional settings and establish a connection between IB-estimators and PT-unbiasedness. In particular, we demonstrate, under suitable (and reasonable) conditions given in Section \ref{sec:setting} that when the initial estimator is consistent or ``slightly'' asymptotically biased, the IB-estimator is a PT-unbiased estimator. In addition, we propose a two-step approach to attain PT-unbiasedness from inconsistent initial estimators. These new results have noticeable advantages in practical applications. As an illustration, in Section \ref{Sec_applica} we present a simulation study in a logistic regression setting (with and without random effects). We consider several slightly inconsistent initial estimators with different properties (e.g. ``robustness'' to separation and/or data contamination) and illustrate the effectiveness, in terms of finite sample bias and mean squared error, of the resulting IB-estimators. These simulation results are in line with the PT-unbiasedness property of IB-estimators. Moreover, they illustrate that IB-estimators can have additional properties (such as robustness) that are transferred from their initial estimator. 

\section{Mathematical Setup}
\label{sec:setting}

Let $\mathbf{X}(\bm{\theta}_0, n)\in\real^n$ denote a random sample generated under model $F_{\bm{\theta}_0}$ (possibly conditional on a set of fixed covariates), where $\bm{\theta}_0\in\bm{\Theta}\subset\real^p$ is the parameter vector we wish to estimate. Moreover, we consider \textit{simulated} samples that we denote as $\mathbf{X}_h^*(\bm{\theta}, n)\in\real^n$ for $\bm{\theta}\in\bm{\Theta}$ and where the subscript $h = 1,\ldots,H$ identifies the distinct samples. Next, we define $\hat{\bm{\pi}}\left(\bm{\theta}_0, n\right)$ and $\hat{\bm{\pi}}_h^*\left(\bm{\theta}, n\right)$ as estimators of $\bm{\theta}_0$ and $\bm{\theta}$ based on the samples $\mathbf{X}(\bm{\theta}_0, n)$ and $\mathbf{X}^*_h(\bm{\theta}, n)$, respectively. In general, we assume that $\hat{\bm{\pi}}\left(\bm{\theta}_0, n\right)$ is a biased estimator of $\bm{\theta}_0$ where this bias may be asymptotic and/or finite sample in nature.

In order to correct the bias of $\hat{\bm{\pi}}(\bm{\theta}_0,n)$, similarly to \cite{guerrier2018simulation}, we propose to use the IB to produce the sequence $\Big\{\hat{\bm{\theta}}^{(k)}\Big\}_{k \in \mathbb{N}}$ defined as
\begin{equation}
		\hat{\bm{\theta}}^{(k)} \equiv \hat{\bm{\theta}}^{(k-1)} + \left[ \hat{\bm{\pi}}(\bm{\theta}_0, n) - \frac{1}{H} \sum_{h = 1}^H  \hat{\bm{\pi}}_h^*\left(\hat{\bm{\theta}}^{(k-1)}, n\right) \right],
		\label{eq:iterboot}
\end{equation}
where $\hbtheta^{(0)}\in\bT$ and $H$ is an integer chosen \textit{a priori}. For example, when $\hbpi(\bto,n)\in\bT$, then this estimator can be used as the initial value of the above sequence. When this iterative procedure converges, we define the IB-estimator $\hat{\bm{\theta}} \in \real^p$ as the limit in $k$ of $\hat{\bm{\theta}}^{(k)}$. It is possible that the sequence is stationary in that there exists a $k^* \in \mathbb{N}$ such that for all $k \geq k^*$ we have $\hat{\bm{\theta}} = \hat{\bm{\theta}}^{(k)}$. Moreover, when $\hat{\bm{\pi}}(\bm{\theta}_0, n)$ in (\ref{eq:iterboot}) is a consistent estimator of $\btheta_0$, the first step of the IB sequence (i.e. $k=1$) is equivalent to the standard bootstrap bias correction proposed by \cite{efron1994introduction} which can be used to reduce the bias of a consistent estimator (under some appropriate conditions). 
\begin{Remark}
\label{Remark:notation}
    Before proceeding to the presentation of the main results of this article, the authors would like to make the reader aware of the fact that the notation used in the article is slightly different from the one used in its appendices where most of the proofs are presented. Indeed, a simplified notation is used in the article to enhance readability while a more detailed and precise notation is used in the appendices to prove the results and avoid misleading conclusions. The notation used in the appendices is presented and explained in Appendix \ref{app:notations}.
\end{Remark}
\vspace{0.25cm}
In order to present the approach, it is useful to define the set $\widehat{\bm{\Theta}}$ as
\begin{equation}
\widehat{\bm{\Theta}} \equiv \argzero_{\bm{\theta} \in
      \bm{\Theta}} \;  \hat{\bm{\pi}}(\bm{\theta}_0, n) - \frac{1}{H} \sum_{h = 1}^H  \hat{\bm{\pi}}^*_h(\bm{\theta}, n),
    \label{eq:indirectInf:hTimesN}
\end{equation}
which will allow us to study the properties of the sequence $\hat{\bm{\theta}}^{(k)}$. As one may expect, under some suitable conditions, we have that $\hat{\bm{\theta}} \in \widehat{\bm{\Theta}}$ thereby explaining why some of our assumptions concern $\bm{\Theta}$.
\begin{Remark}
	\label{rem:ind:inf}
	The elements of $\widehat{\bm{\Theta}}$ are in fact special cases of an indirect inference estimator. Indeed, using the previous notation, an indirect inference estimator  can be defined as 
	\begin{equation}
	    \hat{\bm{\theta}}^\ast \in \widehat{\bm{\Theta}}^\ast \equiv \argmin_{\bm{\theta} \in
	      \bm{\Theta}} \; \Big\| \hat{\bm{\pi}}(\bm{\theta}_0, n) - \frac{1}{H} \sum_{h = 1}^H  \hat{\bm{\pi}}^\ast_h(\bm{\theta}, n) \Big\|^2_{\bm{\Phi}} \; ,
	      \label{eq:indirectInf:rem}
	 \end{equation}
where $\bm{\Phi}$ is a positive-definite matrix. In the case where $\widehat{\bT}\neq \emptyset$, we have that $\widehat{\bm{\Theta}} = \widehat{\bm{\Theta}}^\ast$, which clearly shows the equivalence between \eqref{eq:indirectInf:hTimesN} and \eqref{eq:indirectInf:rem}, for any positive-definite matrix $\bm{\Phi}$. Note that $\hat{\bm{\pi}}(\bm{\theta}, n)$ is referred to as the auxiliary estimator in the indirect inference framework.
\end{Remark}
\vspace{0.25cm}
Before presenting the properties of the IB sequence, we first describe the assumptions we consider.

\setcounter{Assumption}{0}
\renewcommand{\theHAssumption}{otherAssumption\theAssumption}
\renewcommand\theAssumption{\Alph{Assumption}}
\begin{Assumption}
\label{assum:A}
Let $\bm{\Theta}$ be a convex and compact subset of $\real^p$ such that $\btheta_0 \in \operatorname{Int}(\bm{\Theta})$ and
		\begin{equation*}
			 \argzero_{\bm{\theta} \in
		      \real^p \setminus \bm{\Theta}} \; 
		      \hat{\bm{\pi}}(\bm{\theta}_0, n) - \frac{1}{H} \sum_{h = 1}^H  \hat{\bm{\pi}}^*_h(\bm{\theta}, n)
		      = \emptyset\,.
		\end{equation*}
\end{Assumption}
\vspace{0.25cm}
Assumption \ref{assum:A} is quite mild although stronger than necessary. Indeed, the convexity condition and the fact that $\bto$ is required to be in the interior of $\bm{\Theta}$ are convenient to ensure that expansions can be made between $\bto$ and an arbitrary point in $\bm{\Theta}$. Similarly, the compactness is particularly convenient as it allows to bound certain quantities. Finally, the last part of the assumption essentially ensures that no solution outside of $\bm{\Theta}$ exists for the estimator defined in (\ref{eq:indirectInf:hTimesN}). Therefore, the solution set $\widehat{\bm{\Theta}}$ may equivalently be written as
\begin{equation*}
	    \widehat{\bm{\Theta}} = \argzero_{\bm{\theta} \in
		      \real^p} \;  \hat{\bm{\pi}}(\bm{\theta}_0, n) - \frac{1}{H} \sum_{h = 1}^H  \hat{\bm{\pi}}^\ast_h(\bm{\theta}, n).
\end{equation*}
Next, we impose some conditions on the initial estimator $\hat{\bm{\pi}}(\bm{\theta},n)$. Denoting $\mathbf{a}_j$ as being the $j$-th entry of a generic vector $\mathbf{a} \in \real^p$, the assumption is as follows.

\setcounter{Assumption}{1}
    \renewcommand{\theHAssumption}{otherAssumption\theAssumption}
    \renewcommand\theAssumption{\Alph{Assumption}}
\begin{Assumption}
	\label{assum:B}
	For all $\left(\bm{\theta}, \, n\right) \in \bm{\Theta} \times \mathbb{N}^\ast$, the expectation
			$\bm{\pi} \left( \bm{\theta}, n\right) \equiv \mathbb{E}\left[\hat{\bm{\pi}}(\bm{\theta},n)\right]$
		exists and is finite, i.e. $\lvert\bm{\pi}_j \left(\bm{\theta}, n \right)\rvert < \infty$ for all $j = 1, \ldots, p$.
		Moreover, for all $j = 1, \ldots, p$ $\bm{\pi}_j \left( \bm{\theta}\right) \equiv \displaystyle{\lim_{n \to \infty}} \bm{\pi}_j \left( \bm{\theta}, n\right)$ exists.
\end{Assumption}
\vspace{0.25cm}
Assumption \ref{assum:B} is likely  to be satisfied in the majority of practical situations 
and is particularly useful as it allows to decompose the estimator $\hat{\bm{\pi}}(\bm{\theta},n)$ into a non-stochastic component $\bm{\pi} \left( \bm{\theta}, n\right)$ and a random term $\mathbf{v} \left(\bm{\theta},n\right)$. Indeed, using Assumption \ref{assum:B}, we can write:
	\begin{equation}
	   \hat{\bm{\pi}}(\bm{\theta},n)  = \bm{\pi} \left( \bm{\theta}, n\right)+ \mathbf{v} \left(
	    \bm{\theta},n\right),
	  \label{bias:estim}
	\end{equation}
where $\mathbf{v} \left(\bm{\theta},n\right) \equiv \hat{\bm{\pi}}(\bm{\theta},n) - \bm{\pi} \left( \bm{\theta}, n\right)$ is a zero-mean random vector. It must be noticed that this assumption does not imply that the functions $\bm{\pi} \left( \bm{\theta}, n\right)$ and $\hat{\bm{\pi}}(\bm{\theta},n)$ are continuous in $\btheta$. Consequently, when using the IB (see Theorem \ref{thm:conv:consit:iter:boot} below), it is possible to make use of initial estimators  that are not continuous in $\btheta$. In \cite{guerrier2018simulation}, the continuity of these functions was (implicitly) assumed while the weaker assumptions presented in this article extend the applicability of this framework also to models for discrete data such as logistic regression, which is discussed in Section~\ref{Sec_applica}. We now move onto Assumption \ref{assum:C} (below) that imposes some restrictions on the random vector $\mathbf{v} \left(\bm{\theta},n \right)$. 

\setcounter{Assumption}{2}
    \renewcommand{\theHAssumption}{otherAssumption\theAssumption}
    \renewcommand\theAssumption{\Alph{Assumption}}
	\begin{Assumption}
		\label{assum:C}
		The second moment of $\mathbf{v}(\bm{\theta},n)$ exists and there is a real $\alpha > 0$ such that for every $\bm{\theta} \in \bm{\Theta}$ and for all $j = 1,\,\ldots,\,p$, we have
		\begin{equation*}
			\mathbf{v}_j(\bm{\theta},n) = \mathcal{O}_{\rm p}(n^{-\alpha})\;\;\;\;
			 \text{and}\;\;\;\;
			 \lim_{n \to \infty} \; \frac{p^{\nicefrac{1}{2}}}{n^\alpha} = 0.
		\end{equation*}
	\end{Assumption}
\vspace{0.25cm}
Assumption \ref{assum:C} is frequently employed and typically very mild as it simply requires that the variance of $\mathbf{v}_j \left(\bm{\theta},n\right)$ exists and goes to zero as $n$ increases. For example, if $\hat{\bm{\pi}}(\bm{\theta},n)$ is \mbox{$\sqrt{n}$-consistent} (towards $\bm{\pi} \left( \bm{\theta}\right)$) then we would have $\alpha = \nicefrac{1}{2}$ and Assumption \ref{assum:C} would simply require that $p = o(n)$. Next, we consider the bias of $\hat{\bm{\pi}}(\bm{\theta},n)$ for $\bm{\theta}$ and we let $\mathbf{d}\left(\bm{\theta}, n\right) \equiv \bm{\pi}(\bm{\theta},n) - \bm{\theta}$. Using this definition we can rewrite (\ref{bias:estim}) as follows 
	\begin{equation*}
	   \hat{\bm{\pi}}(\bm{\theta},n)  = \bm{\theta} + \mathbf{d} \left( \bm{\theta}, n\right)+ \mathbf{v} \left(
	    \bm{\theta},n\right).
	\end{equation*}
	Moreover, the bias function $\mathbf{d}\left(\bm{\theta}, n\right)$ can always be expressed as follows
	\begin{equation}
	  \mathbf{d}\left(\bm{\theta}, n\right) = \mathbf{a}(\bm{\theta})+\mathbf{c}(n) + \mathbf{b} \left(\bm{\theta} , n\right),
	  \label{def:fct:bias:d}
	\end{equation}
	where $\mathbf{a}(\btheta)$ is defined as the asymptotic bias in the sense that $\mathbf{a}_j(\bm{\theta}) \equiv \displaystyle{\lim_{n \to \infty}} \; \mathbf{d}_j\left(\bm{\theta}, n\right)$ for $j = 1, \ldots, p$, while
	$\mathbf{c}(n)$ and $\mathbf{b} \left(\bm{\theta} , n\right)$ are used to represent the finite sample bias.
	More precisely, $\mathbf{a}(\btheta)$ contains all the terms that are strictly functions of $\btheta$, $\mathbf{c}(n)$ the terms that are strictly functions of $n$ and $\mathbf{b} \left(\bm{\theta}, n\right)$ the rest. This definition implies that if $\mathbf{d}\left(\bm{\theta}, n\right)$ contains a constant term, it is included in $\mathbf{a}(\btheta)$.  Moreover, the function $\mathbf{b} \left(\bm{\theta} , n\right)$ can always be decomposed into a linear and a non-linear term in $\btheta$, i.e.
	\begin{equation}
		\mathbf{b}\left(\bm{\theta}, n\right) = \mathbf{L}(n) \bm{\theta} + 
		\mathbf{r} \left(\bm{\theta} , n\right),
		\label{def:fct:bias:b}
	\end{equation}
	where $\mathbf{L}(n) \in \real^{p \times p}$ and $\mathbf{r}\left(\bm{\theta}, n\right)$ does not contain any linear term in $\bm{\theta}$. Denoting $\mathbf{A}_{j,l}$ as the entry in the $j$-th row and $l$-th column of a generic matrix $\mathbf{A} \in \real^{p \times p}$, we consider Assumption~\ref{assum:D} which imposes some restrictions on the different terms of the bias function. 

\setcounter{Assumption}{3}
\renewcommand\theAssumption{\Alph{Assumption}}
\begin{Assumption}
\label{assum:D}
	The bias function $\mathbf{d}\left(\bm{\theta}, n\right)$ is as follows:
			\begin{enumerate}
		    \item The function $\mathbf{a}(\bm{\theta})$ is a contraction map in that for any $\bm{\theta}_1,\bm{\theta}_2\in\bm{\Theta}$ such that $\bm{\theta}_1 \neq \bm{\theta}_2$ we have 
		    \begin{equation*}
		        \big\lVert \mathbf{a}(\bm{\theta}_2)-\mathbf{a}(\bm{\theta}_1)  \big\rVert_2 <  \big\rVert \bm{\theta}_2-\bm{\theta}_1 \big\lVert_2\,.
		    \end{equation*}
		    
		    \item There exist real $\beta, \gamma > 0$ such that for all $\bm{\theta} \in \bm{\Theta}$ and any $j,l=1,\dots,p$, we have
		    \begin{equation*}
		    \begin{aligned}
			  \mathbf{L}_{j,l}(n) = \mathcal{O}(n^{-\beta}), \;\;\;  \mathbf{r}_j\left(\bm{\theta}, n\right) = \mathcal{O}(n^{-\gamma}),\;\;\; \lim_{n \to \infty} \; \frac{p^{\nicefrac{3}{2}}}{n^\beta} = 0 \;\;\;\;
		      \text{and} \;\;\; \lim_{n \to \infty} \; \frac{p^{\nicefrac{1}{2}}}{n^\gamma} = 0.
		    \end{aligned}
		    \end{equation*}
		    \item Defining $c_n \equiv \displaystyle{\max_{j=1,\dots,p}} \mathbf{c}_j(n)$ for all $n \in \mathbb{N}^*$, we require that the sequence $\left\{c_n\right\}_{n\in\mathbb{N}^*}$ is such that 
            \begin{equation*}
	        \lim_{n \to \infty} \; p^{\nicefrac{1}{2}}c_n = 0.
	        \end{equation*}
		\end{enumerate}
\end{Assumption}
\vspace{0.25cm}
The first part of Assumption \ref{assum:D} is reasonable provided that the asymptotic bias is relatively ``small''  compared to $\btheta$ (up to a constant term). For example, if $\mathbf{a}(\btheta)$ is sublinear in $\btheta$, i.e. $\mathbf{a}(\btheta) = \mathbf{M} \btheta + \mathbf{s}$, then the first part of Assumption \ref{assum:D} would be satisfied if the Frobenius norm is such that $||\mathbf{M}||_F < 1$ since
	\begin{equation}\label{contract:map}
	    \Big\lVert \mathbf{a}(\bm{\theta}_2)-\mathbf{a}(\bm{\theta}_1)  \Big\rVert_2 = \Big\lVert \mathbf{M}(\bm{\theta}_2-\bm{\theta}_1)  \Big\rVert_2 \leq  \big\lVert\mathbf{M}\big\rVert_F \big\lVert \bm{\theta}_2-\bm{\theta}_1 \big\rVert_2 < \big\lVert \bm{\theta}_2-\bm{\theta}_1 \big\rVert_2 .
	\end{equation}
	Moreover, the function ${\bm{\pi}}(\bm{\theta})$ (as defined in Assumption \ref{assum:B}) is often called the asymptotic binding function in the indirect inference literature (see e.g. \citealp{gourieroux1993indirect}). To ensure the consistency of such estimators, it is typically required for ${\bm{\pi}}(\bm{\theta})$ to be continuous and injective. In our setting, this function is given by ${\bm{\pi}}(\bm{\theta}) = \btheta + \mathbf{a}(\btheta)$ and its continuity and injectivity are directly implied by the first part of Assumption \ref{assum:D}. Indeed, since $\mathbf{a}(\bm{\theta})$ is a contraction map, it is continuous. Moreover, taking $\bm{\theta}_1,\bm{\theta}_2\in\bm{\Theta}$ with
	${\bm{\pi}}(\bm{\theta}_1) = {\bm{\pi}}(\bm{\theta}_2)$, then $
	    \big\lVert {\mathbf{a}}(\bm{\theta}_1) - {\mathbf{a}}(\bm{\theta}_2) \big\rVert = \big\lVert \bm{\theta}_1 - \bm{\theta}_2 \big\rVert$,
	which is only possible if $\btheta_1 = \btheta_2$. Thus, $\bm{\pi}(\bm{\theta})$ is injective.
\begin{Remark}
	\label{Remark:kickIB}
	   In situations where $\mathbf{a}(\bm{\theta})$ is not a  contraction map, a possible solution is to modify the sequence considered in the IB as follows:
	   \begin{equation*}
	       	\hat{\bm{\theta}}^{(k)} \equiv \hat{\bm{\theta}}^{(k-1)} + \varepsilon_k\left[ \hat{\bm{\pi}}(\bm{\theta}_0, n) - \frac{1}{H} \sum_{h = 1}^H   \hat{\bm{\pi}}^*_h(\hat{\bm{\theta}}^{(k-1)}, n)\right],
	   \end{equation*}
	   with $\varepsilon_k \in (0, 1]$ for all $k \in \mathbb{N}$. If $\varepsilon_k=\varepsilon$ (i.e. a constant), $\mathbf{a}(\bm{\theta})$ does not need to be a contraction map. Indeed, if $\bm{\Theta}$ is bounded and $\mathbf{a}(\bm{\theta})$ is differentiable, it is always possible to find an $\varepsilon$ such that $\varepsilon \, \mathbf{a}(\bm{\theta})$ is a contraction map. A formal study on the influence of $\varepsilon_k$ on the IB algorithm is, however, left for further research.  
	\end{Remark}
\vspace{0.25cm}
While more general, the second part of Assumption \ref{assum:D} would be satisfied, for example, if $\mathbf{b}(\bm{\theta},n)$ is a sufficiently smooth function in $\bm{\theta}$ and/or $n$, thereby allowing a Taylor expansion, as considered, for example, in \cite{guerrier2018simulation}. Moreover, this assumption is typically less restrictive than the approximations that are commonly used to describe the bias of estimators. Indeed, a common assumption is that the bias of a consistent estimator (including the MLE), can be expanded in a power series in $n^{-1}$ (see e.g. \citealp{kosmidis2014bias} and \citealp{hall1988bootstrap} in the context of the iterated bootstrap), i.e. 
\begin{equation}
    \mathbf{d}(\bm{\theta},n)  = \sum_{j = 1}^{m} \frac{\mathbf{h}^{(j)}(\bm{\theta})}{n^{j}} + \mathbf{g}(\btheta, n),
    \label{Eq_bais-standard}
\end{equation}
where $\mathbf{h}^{(j)}(\bm{\theta})$ is $\mathcal{O}(1)$ elementwise, for $j = 1, \ldots, m$, and $\mathbf{g}(\btheta, n)$ is $\mathcal{O}\left(n^{-(m+1)}\right)$ elementwise, for some $m \geq 1$. The bias function $\mathbf{d}(\bm{\theta},n)$ given in (\ref{Eq_bais-standard}) clearly satisfies the requirements of Assumption~\ref{assum:D}.
Moreover, under the form of the bias postulated in (\ref{Eq_bais-standard}), we have that $\beta,\gamma\geq 1$. If the initial estimator is $\sqrt{n}$-consistent, we have $\alpha = \nicefrac{1}{2}$, and therefore the requirements of Assumptions \ref{assum:C} and \ref{assum:D}, i.e.
	\begin{equation*}
	    \lim_{n \to \infty} \; \max \left( \frac{p^{\nicefrac{1}{2}}}{n^{\min(\alpha,\gamma)}},  \frac{p^{\nicefrac{3}{2}}}{n^\beta} \right) = 0,
	\end{equation*}
	are satisfied if
	\begin{equation}\label{equ:p/n}
	    \lim_{n \to \infty} \; \frac{p^{\nicefrac{3}{2}}}{n} = 0.
	\end{equation}
The last part of Assumption \ref{assum:D} is particularly mild. Indeed, it simply requires that $\mathbf{c}(n)$, i.e. the part of bias that only depends on the sample size $n$, goes to $\mathbf{0}$ faster than $\sqrt{p}$. For the vast majority of estimators the function $\mathbf{c}(n)$ is either very small or equal to $\mathbf{0}$ as the bias generally depends on $\bm{\theta}$.

The assumptions presented in this section are mild and likely satisfied in most practical situations. Moreover, compared to conditions considered in many instances (see e.g. \citealp{guerrier2018simulation} and the references therein), our assumption framework allows to relax various requirements: \textit{(i)} the initial estimator may have an asymptotic bias that depends on $\bm{\theta}$; \textit{(ii)} the initial estimator may also be discontinuous in $\bm{\theta}$ (as is commonly the case when considering discrete data models); \textit{(iii)} the finite sample bias of the initial estimator is allowed to have a more general expression with respect to previously defined bias functions; \textit{(iv)} some of the technical requirements (for example on the topology of $\bm{\Theta}$) are relaxed \textit{(v)} $p$ is allowed to increase with the sample size $n$. Nonetheless, our assumption framework is not necessarily the weakest possible in theory and may be further relaxed (as discussed in Remark \ref{rem:weaker:assum:for:thms} of Appendix \ref{app:main:res}). However, we do not attempt to pursue the weakest possible conditions to avoid overly technical treatments in establishing the theoretical results of the following section.

\section{Main results}
\label{sec:main}

In this section, we study the convergence of the IB sequence along with the properties of the IB-estimator. In particular, in Section~\ref{sec:consistency} we discuss the convergence of the IB sequence together with the consistency of the IB-estimator. In Section~\ref{sec:bias}, under the conditions set in Section \ref{sec:setting}, we show that the IB-estimator is PT-unbiased while the asymptotic distribution of the estimator is discussed in Section \ref{sec:asym:norm}. The asymptotic results presented in this section are somewhat unusual as we always consider arbitrarily large but finite $n$ and $p$ and our results may not be valid when taking the limit in $n$ (and therefore in $p$). Indeed, the (usual) notion of limit we are considering in this article comes from the usual topology (and its induced metric) of $\real^p$ for finite $p$. In order to ensure that these results hold when $p$ is $\infty$, a more detailed topological discussion is needed. However, the difference between the considered framework and others where limits are studied is rather subtle. A more detailed discussion on our asymptotic framework is provided in Appendix~\ref{app:gen:frame}.

\subsection{Consistency}
\label{sec:consistency}

Theorem \ref{thm:conv:consit:iter:boot} below shows that the IB sequence converges (exponentially fast) to the IB-estimator and that the latter is identifiable and consistent.

\begin{Theorem}
		\label{thm:conv:consit:iter:boot}
		Under Assumptions \ref{assum:A}, \ref{assum:B}, \ref{assum:C} and \ref{assum:D}, for all $H \in  \mathbb{N}^\ast$,
		
		\begin{enumerate}
		    
		    \item There exists a $n^* \in \mathbb{N}^\ast$ such that for all $n \in \Ns$ with $n \geq n^*$, $\left\{\hbtheta\right\} = \widehat{\bm{\Theta}}$, i.e. the set $\widehat{\bm{\Theta}}$ is a singleton.
		    
		    \item There exists a $n^* \in \mathbb{N}^\ast$ such that for all $n \in \Ns$ with $n \geq n^*$, the sequence $\left\{\hbtheta^{(k)}\right\}_{k \in \mathbb{N}}$ has the following limit
		    \begin{equation*}
            \lim_{k \to \infty} \; \hbtheta^{(k)} = \hbtheta.
         	\end{equation*}
         	
         	 Moreover, there exists a real $ \epsilon \in (0, \, 1)$ such that for any $k \in \mathbb{N}^\ast$ 
         	       \begin{equation*}
         	           \left\lVert\hbtheta^{(k)} - \hbtheta\right\rVert_2  =\mathcal{O}_{\rm p}({p}^{\nicefrac{1}{2}}\,\epsilon^k).
         	       \end{equation*}
         	
         	\item $\hbtheta$ is a consistent estimator of $\bto$, i.e.,
         	           $\big\lVert\hbtheta - \bto\big\rVert_2  =o_{\rm p}(1)$.
         	
		\end{enumerate} 
	\end{Theorem}
	\vspace{0.25cm}
The proof of Theorem \ref{thm:conv:consit:iter:boot} is given in Appendix \ref{app:main:res}. The three results of this theorem are derived separately. Indeed, the first and the second result of the theorem, which correspond to Lemma \ref{lemma:unique:iter:boot} and Proposition \ref{thm:iter:boot} in Appendix \ref{app:conv:IB}, are proved under a certain set of assumptions that is different from the set of assumptions needed to prove consistency (third result of the theorem). In addition, consistency is proved in two different ways, using two different sets of assumptions (see Proposition \ref{THM:consistency} and Corollary \ref{coro:consist} in Appendices \ref{app:consist:IB} and \ref{app:unbias:IB}, respectively). As a result, Assumptions \ref{assum:A} to \ref{assum:D} condense all the requirements from the previously mentioned assumptions that are found in the appendices, removing possible redundancies between them. A precise account of the dependence structure between the assumptions is given in Figure \ref{Fig_assumptions-mess} of Appendix \ref{app:main:res}.

An important consequence of Theorem \ref{thm:conv:consit:iter:boot} is that the IB provides a computationally efficient algorithm to solve the optimization problem in (\ref{eq:indirectInf:hTimesN}). Moreover, in practical settings, the IB can often be applied to the estimation of complex models where standard optimization procedures used to solve (\ref{eq:indirectInf:hTimesN}) may fail to converge numerically (see e.g. \citealp{guerrier2018simulation}). In practice, the IB procedure is computationally efficient, which is in line with Theorem \ref{thm:conv:consit:iter:boot} showing that $\hbtheta
^{(k)}$ converges to $\hat{\bm{\theta}}$ (in norm) at an exponential rate. Even though the convergence of the algorithm may be slower when $p$ is large, in practical situations, the number of iterations necessary to reach a suitable neighbourhood of the solution appears to be relatively small. 
%
%
For example, if we define $k_n\equiv \mathcal{O}(\log(p n))$, then we have
\begin{equation*}
    \left\lVert\hbtheta^{(k_n^*)} - \bm{\theta}_0\right\rVert_2 \leq \left\lVert\hbtheta^{(k_n^*)} - \hbtheta\right\rVert_2 + \left\lVert\hbtheta - \bm{\theta}_0\right\rVert_2 = o_{\rm p}(1),
\end{equation*}
and therefore $\hbtheta^{(k_n^*)}$ is also a consistent estimator for $\bto$, where $k_n^* \equiv \lfloor k_n \rfloor$.

\subsection{Phase Transition Unbiasedness}
\label{sec:bias}

In the previous section, we studied the consistency of the IB-estimator $\hat{\bm{\theta}}$ and provided a computationally efficient strategy to obtain it. We now investigate the bias of $\hat{\bm{\theta}}$ and we show that this estimator achieves PT-unbiasedness. This result, combined with the guarantee of consistency, is of particular interest in many practical settings. To obtain this property, additional requirements are needed on the function $\mathbf{d}(\bm{\theta}, n)$. For this reason, we introduce Assumption \ref{assum:D'} which combines Assumption \ref{assum:D} and these additional requirements. 

\setcounter{Assumption}{3}
\renewcommand\theAssumption{\Alph{Assumption}$^\ast$}
\begin{Assumption}
\label{assum:D'}
	The bias function $\mathbf{d}\left(\bm{\theta}, n\right)$ is such that:
		\begin{enumerate}
		    \item The asymptotic bias function $\mathbf{a}(\bt)$ can be written as \begin{equation*}
		        \mathbf{a}(\btheta) = \mathbf{M}\btheta + \mathbf{s}, \end{equation*} where $\mathbf{M} \in \real^{p \times p}$ with $|| \mathbf{M} ||_F < 1$ and $\mathbf{s} \in \real^{p}$.
		        
		    \item There exists a $n^* \in \mathbb{N}^*$ such that for all $n \in \mathbb{N}^*$ satisfying $n \geq n^*$, the matrix $(\mathbf{M} + \mathbf{L}(n)+\mathbf{I})^{-1}$ exists.
		    \item There exist real $\beta, \gamma > 0$ such that for all $\bm{\theta} \in \bm{\Theta}$ and any $j,l=1,\dots,p$, we have 
		    \begin{equation*}
			  \mathbf{L}_{j,l}(n) = \mathcal{O}(n^{-\beta}), \;\;\;  \mathbf{r}_j\left(\bm{\theta}, n\right) = \mathcal{O}(n^{-\gamma}),\;\;\; \lim_{n \to \infty} \; \frac{p^{\nicefrac{3}{2}}}{n^\beta} = 0 \;\;\;\;
		      \text{and} \;\;\; \lim_{n \to \infty} \; \frac{p^2}{n^{\gamma}} = 0.
		    \end{equation*}
		    \item The Jacobian matrix $\mathbf{R} (\bm{\theta}, n) \equiv \frac{\partial}{\partial \, \btheta^T}  \mathbf{r}(\bm{\theta}, n) \in \real^{p \times p}$ exists and is continuous in $\btheta\in\bm{\Theta}$ 
		    for all $n \in \Ns$ satisfying $n \geq n^\ast$.
		    \item Defining $c_n \equiv \displaystyle{\max_{j = 1,\dots,p}} \mathbf{c}_j(n)$ for all $n \in \mathbb{N}^*$, the sequence $\left\{c_n\right\}_{n\in\mathbb{N}^*}$ is such that 
            \begin{equation*}
	        \lim_{n \to \infty} \; p^{\nicefrac{1}{2}}c_n = 0.
	        \end{equation*}
		\end{enumerate}
\end{Assumption}
\vspace{0.25cm}
To avoid redundancy, we only discuss the additional requirement of Assumption \ref{assum:D'}. The first one requires $\mathbf{a}(\btheta)$ to be a sublinear function of $\btheta$ with its linear term $\mathbf{M} \in \real^{p \times p}$ satisfying $|| \mathbf{M} ||_F < 1$. As discussed after Assumption \ref{assum:D}, this condition implies that $\mathbf{a}(\btheta)$ is a contraction map. Moreover, the PT-unbiasedness of $\hbtheta$ is still guaranteed even if the condition $|| \mathbf{M} ||_F < 1$ is not satisfied but the IB sequence may not converge. More details are given in Appendices \ref{app:consist:IB} and \ref{app:unbias:IB}.  

 The second additional assumption requires that the matrix $(\mathbf{M} + \mathbf{L}(n)+\mathbf{I})^{-1}$ exists when $n$ is sufficiently large. This requirement is quite general and is, for example, satisfied if $\mathbf{a}(\btheta)$ is relatively ``small'' compared to $\bm{\theta}$ or if $\hat{\bm{\pi}}(\bm{\theta}_0, n)$ is a consistent estimator of $\bm{\theta}_0$. Interestingly, this part of the assumption can be interpreted as requiring that the matrix $(\mathbf{M} + \mathbf{I})^{-1}$ exists, which directly implies that the binding function $\bm{\pi}(\btheta)$ is injective. 
 
 The third one requires that $\displaystyle\lim_{n \to \infty}p^{2}n^{-\gamma} = 0$. In many practical settings it is reasonable\footnote{For example using the bias function proposed in (\ref{Eq_bais-standard}) and assuming that the first term of the expansion $\mathbf{h}^{(1)}(\btheta)$ is linear.} to assume that $\gamma =2$ which implies that this condition is satisfied if {$\displaystyle\lim_{n \to \infty} p/n = 0$}. This latter condition is less strict than \eqref{equ:p/n}. It is also easy to see that $\displaystyle\lim_{n \to \infty}p^{2}n^{-\gamma} = 0$ implies the similar requirement used in Assumption \ref{assum:D}, namely $\displaystyle\lim_{n \to \infty} \; p^{\nicefrac{1}{2}}n^{-\gamma} = 0$.

Finally, the last additional requirement implies that $\mathbf{R} (\bm{\theta}, n)$ is continuous in $\btheta\in\bm{\Theta}$ when $n$ is large enough. A useful consequence of this requirement combined with the compactness of $\bm{\Theta}$ (Assumption \ref{assum:A}) is that $\mathbf{r}_j(\hbtheta, n)$ and $\mathbf{R}_{j,l} (\hbtheta, n)$ are bounded random variables for all $j,l = 1,\dots,p$.

Before stating Theorem \ref{THM:bias} which ensures that $\hat{\bm{\theta}}$ is a PT-unbiased estimator, we introduce a particular asymptotic notation and state a lemma that provides a strategy for proving PT-unbiasedness. Let $\mathcal{A} \subset \N$ and let $f(n)$ and $g(n)$ be real-valued functions with $g(n)$ being strictly positive for all $n \in \N^\ast$. We write $f(n) ~=~ \mathcal{O}_{\delta \in \mathcal{A}}\left(g(n)^\delta\right)$ if and only if :
\begin{equation}
\label{def:newO}
\exists n^\ast > 0, \;\; \forall \delta \in \mathcal{A}, \;\; \exists M_\delta > 0 \;\; \text{such that} \;\; \forall n \geq n^\ast, \;\; |f(n)| \leq M_\delta g(n)^\delta.
\end{equation}
When $\mathcal{A}$ is a finite set we have 
\begin{equation*}
   \mathcal{O}_{\delta \in \mathcal{A}}\left(g(n)^\delta\right) \;\;\;\Longleftrightarrow \;\;\;\mathcal{O}\left(g(n)^\delta\right) \;\text{for all}\; \delta \in \mathcal{A},
\end{equation*}
but this equivalence is not preserved when $\mathcal{A}$ is an infinite set. Indeed, the notation
\begin{equation*}
    f(n)=\mathcal{O}\left(g(n)^\delta\right) \;\text{for all}\; \delta \in \mathcal{A},
\end{equation*}
means that
\begin{equation*}
   \forall \delta \in \mathcal{A}, \;\; \exists n_\delta > 0, \;\; \exists M_\delta > 0  \;\; \text{such that} \;\; \forall n \geq n_\delta \;\; |f(n)| \leq M_\delta g(n)^\delta.
\end{equation*}
Therefore, this definition does not guarantee the existence of a $n^\ast \in \N^*$ that bounds the sequence $\left\{n_\delta\right\}_{\delta \in \mathcal{A}}$ as it may diverge. However, if $\mathcal{A} \subset \N^\ast$ is an infinite set, the following equivalence holds true %
\begin{equation}
\label{equiv:newO}
   \mathcal{O}_{\delta \in \mathcal{A}}\left(g(n)^\delta\right) \;\;\;\Longleftrightarrow \;\;\; \mathcal{O}_{\delta \in \N}\left(g(n)^\delta\right).
\end{equation}
Lemma \ref{lemma:newO} below shows how this particular notation is relevant to show PT-unbiasedness. The proof of this lemma is given below the statement as it is short and provides insight on this property.

\begin{Lemma}
\label{lemma:newO}
Let $g(n)$ be a strictly positive real-valued function such that 
$\displaystyle\lim_{n\to\infty} g(n) = 0$.
If\\ $f(n)~=~\mathcal{O}_{\delta \in \N}\left(g(n)^\delta\right)$
then  there exists a $n^\ast \in \N^\ast$ such that $f(n) = 0$ for all $n \geq n^\ast$.
\end{Lemma}
\vspace{0.25cm}

\begin{proof}
Let $\delta\in\N$. Since $f(n) = \mathcal{O}_{\delta \in \N}\left(g(n)^\delta\right)$ there exists $n^\ast \geq 0$ and $M_\delta > 0$ such that $|f(n)| \leq M_\delta \,g(n)^\delta$ for all $n \geq n^\ast$. We can therefore consider the following quantity
\begin{equation*}
  \bar{M}_\delta \equiv \inf\left\{M > 0 \;\; \Bigg| \;\; \begin{array}{l}
      \exists n^*>0, \, \text{such that} \;\; |f(n)| \leq M \, g(n)^\delta,\;\; \forall n\geq n^\ast  \\
     \text{and} \;\; \forall \delta'\hspace{-0.1cm}\in\hspace{-0.1cm}\N\backslash\hspace{-0.1cm}\left\{\delta\right\}, \, \exists M_{\delta'} > 0, \, |f(n)| \leq M_{\delta'} \, g(n)^{\delta'},\;\; \forall n\geq n^\ast
  \end{array}\right\}.
\end{equation*}
Defining $\bar{M}_{\delta+1}$ in a similar fashion, we have that for all $n\geq n^\ast$ (where $n^\ast$ is the max of the ones in the definitions of $\bar{M}_{\delta}$ and $\bar{M}_{\delta+1}$),
\begin{equation*}
    |f(n)| \leq \bar{M}_{\delta+1}\, g(n)^{\delta+1} = \bar{M}_{\delta+1} \,g(n)\, g(n)^\delta. 
\end{equation*}
By definition of $\bar{M}_\delta$, we have 
\begin{equation}
    0\leq \bar{M}_\delta \leq \bar{M}_{\delta+1}\, g(n),
    \label{eq:ineq:lem}
\end{equation}
for all $n \geq n^\ast$, which implies that $\bar{M}_\delta = 0$. Indeed, supposing the contrary of \eqref{eq:ineq:lem}, there exists $n_0 \geq n^\ast$ such that $\bar{M}_{\delta} > \bar{M}_{\delta+1}\, g(n_0)$. Without loss of generality we can suppose that the function $g(n)$ is decreasing\footnote{Indeed, we can define a decreasing step function $\tilde{g}(n)$ such that $\displaystyle\lim_{n\to 0}\tilde{g}(n) = 0$ and $g(n) \leq \tilde{g}(n)$, $\forall n \in \N^\ast$.}. Therefore, setting $M'_\delta \equiv \bar{M}_{\delta+1}\, g(n_0) > 0$  and redefining $n^\ast$ to be $n_0$, we have that for all $n \geq n^\ast$
\begin{equation*}
     |f(n)| \leq \bar{M}_{\delta+1}\, g(n)^{\delta+1} \leq M'_\delta \, g(n)^\delta, 
\end{equation*}
which contradicts the minimality of $\bar{M}_\delta$. Since $\delta$ was chosen arbitrarily in $\N$, we have that $\bar{M}_\delta = 0$ for all $\delta \in \N$, which concludes the proof.  
\end{proof}
\vspace{0.25cm}

The conditions of Lemma \ref{lemma:newO} may be considered as being too strong \textit{per se}, but as the following two examples show, these conditions are actually necessary.
\begin{enumerate}
    
    \item Suppose $\displaystyle\lim_{n\to\infty} g(n) = 0$ and $f(n)~=~\mathcal{O}\left(g(n)^\delta\right)$ for all $\delta \in \N$. 
    
    In this case, we can consider $g(n) \equiv n^{-1}$ and $f(n) \equiv \exp(-n)$ for all $n \in \N^\ast$. 
     Clearly, $f(n)~=~\mathcal{O}\left(g(n)^\delta\right)$ for all $\delta \in \N$ but $f(n) \neq 0$ for all $n \in \N^\ast$.
    
    \item Suppose $\displaystyle\lim_{n\to\infty} g(n) \neq 0$ and $f(n)~=~\mathcal{O}_{\delta \in \N}\left(g(n)^\delta\right)$.
    
    In this case, we consider a decreasing function $g(n)$ such that $0 < \displaystyle\lim_{n\to\infty} g(n) = c < 1$ and we define $f(n) \equiv d > 0$ for all $n \in \N^\ast$. Since $g(n)$ is decreasing and $\displaystyle\lim_{n\to\infty} g(n) < 1$, there exists $n^\ast \in \N$ such that $0 < g(n) < 1$ for all $n \geq n^\ast$. Setting $M_\delta \equiv d  c^{-\delta}$ for all $\delta \in \N$, we have for all $n \geq n^\ast$ 
    \begin{equation*}
        \vert f(n) \vert = d = d g(n)^{-\delta} g(n)^\delta \leq d c^{-\delta} g(n)^\delta = M_\delta g(n)^\delta.
    \end{equation*}
    Therefore, $f(n)~=~\mathcal{O}_{\delta \in \N}\left(g(n)^\delta\right)$ but $f(n) \neq 0$ for all $n \in \N^\ast$.
\end{enumerate}

Theorem \ref{THM:bias} (below) delivers the conditions that guarantee $\hbtheta$ is a consistent and a PT-unbiased estimator of $\btheta_0$. 

\begin{Theorem}
\label{THM:bias}
	Under Assumptions \ref{assum:A}, \ref{assum:B}, \ref{assum:C} and \ref{assum:D'}, and for all $H \in \mathbb{N}^*$, the IB-estimator $\hbtheta$ is consistent and PT-unbiased. Therefore:
	\begin{enumerate}
	    \item $\hbtheta$ is a consistent estimator of $\bto$, i.e., $\big\lVert\hbtheta - \bto\big\rVert_2  =o_{\rm p}(1)$.
	    \item There exists $n^* \in \mathbb{N}^*$ such that for all $n\in\mathbb{N}^*$ with $n\geq n^*$, we have 
		$\big\lVert	\mathbb{E} [\hbtheta] - \bto  \big\rVert_2 = 0$. 
\end{enumerate}

\end{Theorem}
\vspace{0.25cm}

The proof of this result is given in Appendix \ref{app:main:res} where we show
\begin{equation}
    \Big\lVert	\mathbb{E} [\hbtheta] - \bto  \Big\rVert_2 ~=~\mathcal{O}_{\delta \in \mathbb{N}}\left((p^2n^{-\gamma})^\delta\right).
    \label{eq:bias}
\end{equation}
By Assumption \ref{assum:D'} we have that $p^{2}n^{-\gamma} \to 0$ as $n \to \infty$ and therefore Lemma \ref{lemma:newO} can be applied to \eqref{eq:bias} to conclude that $ \big\lVert	\mathbb{E} [\hbtheta] - \bto  \big\rVert_2 = 0$ for all $n \geq n^\ast$. In practice, as illustrated within the simulation settings in Section \ref{Sec_applica}, the value $n^*$ appears to be quite small as  the IB-estimator does indeed appear to be unbiased.

In Appendix \ref{app:unbias:IB} we demonstrate the PT-unbiasedness result presented in Theorem \ref{THM:bias} under weaker assumptions. With this set of weaker assumptions the consistency of the estimator is not guaranteed. For example, the rate of convergence of the initial estimator (discussed in Assumption \ref{assum:C}) has to be considered to establish consistency results (in particular when $p$ diverges) but is not needed to study the bias of $\hat{\bm{\theta}}$.

Finally, Theorem \ref{THM:bias} relies on the additional assumption that the asymptotic bias $\mathbf{a}(\btheta)$ is sublinear, which can be quite restrictive.  However, if $\hbpi(\bm{\theta}_0,n)$ is a consistent estimator of $\bm{\theta}_0$, then this restriction is automatically satisfied. Nevertheless, considering the third part of Theorem \ref{thm:conv:consit:iter:boot}, a possible approach that would guarantee a PT-unbiased estimator is to obtain the IB-estimator $\hbtheta$ from an inconsistent initial estimator $\hat{\bm{\pi}}(\bm{\theta},n)$ and then to compute a new IB-estimator with the initial estimator $\hat{\bm{\pi}}^{\text{NEW}}(\bm{\theta}_0, n) \equiv \hbtheta$. In practice, however, this computationally intensive approach is probably unnecessary since the (one step) IB-estimator appears to eliminate the bias almost completely when considering inconsistent initial estimators $\hat{\bm{\pi}}(\bm{\theta},n)$ (see e.g. Section \ref{sec:logistic}). In fact, we conjecture that for a larger class of asymptotic bias functions, Theorem \ref{THM:bias} remains true, although the verification of this conjecture is left for further research.

\subsection{Asymptotic Normality}
\label{sec:asym:norm}

In this section we study the asymptotic distribution of the IB-estimator. In the case where $p \to c < \infty$ as $n \to \infty$, the results of \cite{gourieroux1993indirect} can be directly applied to obtain
\begin{equation}
    \sqrt{n} \left(\hat{\bm{\theta}}-\bto\right) \xrightarrow{d} \mathcal{N}\left(\mathbf{0}, \bm{\Xi}_H\right),
\label{eq:asym:norm:p:fixed}
\end{equation}
where
\begin{equation*}
    \bm{\Xi}_H \equiv \left(1+\frac{1}{H}\right) \left(\mathbf{B}(\bto)^{-1}\right)^T \bm{\Sigma}(\bto) \mathbf{B}(\bto)^{-1},
    \label{Eq_ass-var-JINI}
\end{equation*}
and where
\begin{equation*}
	\mathbf{B}(\bto) \equiv \mathbf{I} + \left.\frac{\partial}{\partial \bm{\theta}^T} \mathbf{a}(\bm{\theta})\right|_{\bm{\theta} = \bm{\theta}_0}, \;\;\;\;\; \bm{\Sigma}(\bto) \equiv \lim_{n \to \infty} \; \sqrt{n} \var\left(\mathbf{v}(\bm{\theta}_0, n)\right).
\end{equation*}
Apart from technical requirements (such as the continuity of $\mathbf{B}(\bm{\theta})$ in $\bto$), the main condition ensuring the validity of (\ref{eq:asym:norm:p:fixed}) is the asymptotic normality of $\mathbf{v}(\bm{\theta}_0, n)$ (or $\hat{\bm{\pi}}(\bm{\theta}_0, n)$), i.e.
\begin{equation*}
    \sqrt{n} \, \mathbf{v}(\bm{\theta}_0, n) \xrightarrow{d} \mathcal{N}\left(\mathbf{0}, \bm{\Sigma}(\bto)\right).
\end{equation*}
However, in the case where $p \to \infty$ as $n \to \infty$, such a condition requires to be considerably modified by assuming instead that $\mathbf{v}(\bm{\theta}_0, n)$ satisfies a specific Gaussian approximation. Moreover, several other additional technical requirements are needed to adapt (\ref{eq:asym:norm:p:fixed}) to the case where $p$ diverges. To avoid an overly technical discussion in the main text, these conditions are stated and discussed in Appendix \ref{app:asymp:norm:IB} together with Theorem \ref{app:thm:main:3} in Appendix \ref{app:main:res} which proposes a possible high dimensional counterpart of (\ref{eq:asym:norm:p:fixed}). Unfortunately, many of the conditions needed to derive these results can be strong and difficult to verify for a specific model and initial estimator. In the case where $p \to c < \infty$ as $n \to \infty$, Theorem \ref{app:thm:main:3} becomes equivalent to (\ref{eq:asym:norm:p:fixed}) and relies on very similar requirements as those needed for the asymptotic results in \cite{gourieroux1993indirect}. Informally, Theorem \ref{app:thm:main:3} states that for all $\mathbf{u} \in \real^p$ such that $||\mathbf{u}||_2 = 1$, we have
\begin{equation}
     \sqrt{n} \left(1 + \frac{1}{H}\right)^{-1/2} \mathbf{u}^T \bm{\Sigma}(\bto, n)^{-1/2} \mathbf{B}(\bto, n) \left(\hat{\bm{\theta}}-\bto\right) \xrightarrow{d} \mathcal{N}\left({0}, 1\right),
\label{eq:asym:norm:p:div}
\end{equation}
where
\begin{equation*}
	\mathbf{B}(\bto, n) \equiv \left.\frac{\partial}{\partial \bm{\theta}^T} \bm{\pi}(\bm{\theta},n)\right|_{\bm{\theta} = \bm{\theta}_0}, \;\;\;\;\; \bm{\Sigma}(\bto, n) \equiv \; \sqrt{n} \var\left(\mathbf{v}(\bm{\theta}_0, n)\right).
\end{equation*}
One interesting difference in the conditions needed to derive (\ref{eq:asym:norm:p:fixed}) and (\ref{eq:asym:norm:p:div}), is that the former is valid for all $H \in \mathbb{N}^*$ while the latter requires $H=\mathcal{O}(p^4)$. Therefore, it suggests that the ``quality'' (and validity) of the approximation depends on $H$ when $p$ diverges. However, the conditions of Theorem \ref{app:thm:main:3} are sufficient but may not be necessary thereby implying that  $H=\mathcal{O}(p^4)$ may not always be needed as a condition for this approximation.

In practice, the estimation of the variance of the IB-estimator can be obtained through different methods. A simple approach takes advantage of the results of the last iteration of the IB sequence to construct a parametric bootstrap estimator of the variance of the initial estimator $\hat{\bm{\pi}}(\bm{\theta}_0,n)$. Indeed, in the last iteration of the IB sequence, $H$ samples are simulated under $\hat{\bm{\theta}}$ (which is a consistent estimator), which allows to compute the following quantity (assuming $H$ to be sufficiently large)
\begin{equation*}
     \widehat{\var}(\hat{\bm{\pi}}(\bm{\theta}_0,n)) \equiv \frac{1}{H - 1} \sum_{h = 1}^H   \left[ \hat{\bm{\pi}}^*_h(\hat{\bm{\theta}}, n) - \bar{\bm{\pi}}(\hat{\bm{\theta}}, n)\right]  \left[ \hat{\bm{\pi}}^*_h(\hat{\bm{\theta}}, n) - \bar{\bm{\pi}}(\hat{\bm{\theta}}, n)\right]^T,  \end{equation*}
where $\bar{\bm{\pi}}(\hat{\bm{\theta}}, n) \equiv \frac{1}{H} \sum_{h = 1}^H  \hat{\bm{\pi}}^*_h(\hat{\bm{\theta}}, n)$. Then, an estimator of the covariance matrix of $\hat{\bm{\theta}}$ can be obtained as follows
\begin{equation*}
    \widehat{\var}(\hat{\bm{\theta}}) \equiv  \left(1+\frac{1}{H}\right) \left(\widehat{\mathbf{B}}^{-1}\right)^T \widehat{\var}(\hat{\bm{\pi}}(\bm{\theta}_0,n)) \widehat{\mathbf{B}}^{-1},
\end{equation*}
where $\widehat{\mathbf{B}}$ can be obtained by numerical derivation of $\hat{\bm{\pi}}(\bm{\theta}, n)$ evaluated at $\hat{\bm{\theta}}$.

\section{Application: Logistic Regression Model}
\label{Sec_applica}

In this section, we apply the methodology developed in Sections \ref{sec:setting} and \ref{sec:main} to investigate the performance of IB-estimators in three different settings. As previously mentioned, our conditions allow the initial estimator to be discontinuous in $\bm{\theta}$ and it is therefore of interest to consider the logistic regression model, which may be the most commonly used model for binary (response) data. To illustrate the flexibility of the IB estimation approach, as initial estimators we select slightly modified versions of the MLE and of a robust estimator. As explained further on, these modifications were introduced to allow the estimators to be ``robust'' to the problem of data separation. To compare these IB-estimators, we consider the MLE, the robust estimator as well as the bias corrected estimator proposed by \cite{KoFi:09}. Their respective performance is studied when the data is generated at the model but also under slight model misspecification where ``outliers'' are randomly created. Then, we extend the logistic regression model to include a random intercept, a special case of Generalized Linear Mixed Models (GLMM) for which there is no closed form expression for the likelihood function. In this case, the initial estimator is selected to be a numerically simple approximation to the MLE, and its finite sample behaviour is compared to the ones of the MLE computed using several (more precise) approximation methods.

\subsection{Bias Corrected Estimators for the Logistic Regression Model}
\label{sec:logistic}

\subsubsection{Introduction}

The Logistic Regression Model (LRM) \citep{NeWe:72,McCuNe:89} is one of the most frequently used models for binary response variables conditioned on a set of covariates. However, it is well known that in some quite frequent practical situations, the MLE is biased or its computation can become very unstable, especially when performing some type of resampling scheme for inference. The underlying reasons are diverse, but the main ones are the possibly large $p/n$ ratio, separability (leading to regression slope estimates of infinite value) and data contamination (robustness). 

The first two sources are often confounded and practical solutions are continuously sought to overcome the difficulty in performing ``reasonable'' inference. For example, in medical studies, the bias of the MLE together with the problem of separability has led to a rule of thumb called the number of Events Per Variable (EPV), that is the number of occurrences of the least frequent event over the number of predictors, which is used in practice to choose the maximal number of predictors one is ``allowed'' to use in a LRM (see e.g. \citealp{AuSt:17} and the references therein
).  

The problem of separation or near separation in logistic regression is linked to the existence of the MLE which is not always guaranteed \citep[see][]{Silv:81,AlAn:84,SaDu:86,Zeng:17}. Alternative conditions for the existence of the MLE have been recently developed in \cite{CaSu:19} \citep[see also][]{sur2019modern}. In order to detect separation, several approaches have been proposed (see for instance, \citealp{LeAl:89,Kolo:97,CHRISTMANN2001}). The adjustment to the score function proposed by \cite{Firt:93} has been implemented in \cite{KoFi:09} \citep[see also][]{kosmidis2010} to propose a bias corrected MLE for GLM which has the additional natural property that it is not subject to the problem of separability
. Moreover, the MLE is known to be sensitive to slight model deviations that take the form of outliers in the data leading to the proposal of several robust estimators for the LRM and more generally for GLM (see e.g. \citealp{CaRo:01b,Cize:08,HeCaCoVF:09}  and the references therein
.

Despite all the proposals for finite sample bias correction, separation and data contamination problems, no estimator has so far been able to handle the three potential sources of bias jointly. In this section, we make use of the IB-estimator which is built through a simple adaptation of available estimators. Although the latter estimators might possibly not be the best ones for this problem, they at least guarantee, at a reasonable computational cost, a reduced finite sample bias which is comparable, for example, to that of the bias reduced MLE of \cite{KoFi:09} in uncontaminated data settings as well as reducing the bias in contaminated data settings. Moreover, in both cases, we adapt the initial estimator so that it is not affected by the problem of separability.

\subsubsection{Bias Corrected Estimators}
\label{Sec_IB-logistic}

Consider the LRM with response $\mathbf{Y}(\bm{\beta}_0, n)$  and linear predictor $\mathbf{X}\boldsymbol{\beta}$, where $\mathbf{X}$ is an $n\times p$ matrix of fixed covariates with rows $\mathbf{x}_i,i=1,\ldots,n$, and with logit link $\mu_i(\bm{\beta}) \equiv \mathbb{E}[\Y_i(\bm{\beta}, n)]  =\exp(\mathbf{x}_i\boldsymbol{\beta})/(1+\exp(\mathbf{x}_i\boldsymbol{\beta}))$. The MLE for $\bbeta$ is given by 
\begin{equation}
        \hat{\bm{\pi}}(\bm{\beta}_0, n) \equiv \argzero_{\bbeta \in \real^p}  \frac{1}{n}\sum_{i=1}^n\mathbf{x}_i\left[\Y_i(\bm{\beta}_0, n)-\mu_i(\bm{\beta})\right], 
    \label{Eq_MLE-logistic}
\end{equation}
and can be used as an initial estimator in order to obtain the IB-estimator in (\ref{eq:iterboot}). 

To  avoid the (potential) problem of separation, we follow the suggestion of \cite{Rousseeuw2003} to transform the observed responses $\Y_i(\bm{\beta}_0, n)$ to get pseudo-values
\begin{equation}
    \widetilde{\Y}_i(\bm{\beta}_0, n) = (1-\delta)\Y_i(\bm{\beta}_0, n) +\delta \left(1-\Y_i(\bm{\beta}_0, n)\right),
    \label{Eq_pseudo-val}
\end{equation}
for all $i=1,\dots,n$, where $\delta \in \left[0,0.5\right)$ is a (fixed) ``small'' (i.e. close to 0) scalar. For a discussion on the choice of $\delta$ and also possible asymmetric transformations, see  \cite{Rousseeuw2003}. Jointly using these two approaches (i.e. MLE in \eqref{Eq_MLE-logistic} computed on the pseudo-values) as an initial estimator, we denote the estimator resulting from the IB procedure  as the IB-MLE.

As a robust initial estimator, we consider the robust $M$-estimator proposed by \cite{CaRo:01b}, with general estimating function (for GLMs) given by
\begin{equation}
    \bm{\psi}_{\hbpi}\left(\bbeta, \Y_i(\bm{\beta}_0, n)\right) \equiv \psi_c\left(r\left(\bbeta, \Y_i(\bm{\beta}_0, n)\right)\right)w\left(\mathbf{x}_i\right)V^{-1/2}\left(\mu_i(\bbeta)\right)(\partial/\partial\bbeta)\mu_i(\bbeta)-\mathbf{a}\left(\bbeta\right),
    \label{Eq_rob-glm}
\end{equation}
with $r\left(\bbeta, \Y_i(\bm{\beta}_0, n)\right) \equiv \left(\Y_i(\bm{\beta}_0, n)-\mu_i(\bbeta)\right)V^{-1/2}\left(\mu_i(\bbeta)\right)$ being the Pearson residuals and with consistency correction factor
\begin{equation}
    \mathbf{a}\left(\bbeta\right)\equiv\frac{1}{n}\sum_{i=1}^n\mathbb{E}\left[    \psi_c\left(r\left(\bbeta, \Y_i(\bm{\beta}_0, n)\right)\right)w\left(\mathbf{x}_i\right)V^{-1/2}\left(\mu_i(\bbeta)\right)(\partial/\partial\bbeta)\mu_i(\bbeta)\right],
    \label{Eq_rob-consist}
\end{equation}
where the expectation is taken over the (conditional) distribution of the responses $\Y_i(\bm{\beta}_0, n)$ (given $\mathbf{x}_i$). For the LRM, we have  $V\left(\mu_i(\bbeta)\right) \equiv \mu_i(\bbeta)(1-\mu_i(\bbeta))$. We 
compute the robust initial estimator $\hbpi(\bm{\beta}_0, n)$ on the pseudo-values (\ref{Eq_pseudo-val}), using the implementation in the \texttt{glmrob} function of the \texttt{robustbase} package in R \citep{robustbase2018}, with $\psi_c$ in (\ref{Eq_rob-glm}) being the Huber loss function (with default parameter $c$) (see \citealp{huber1964robust}) and $w\left(\mathbf{x}_i\right)=\sqrt{1-h_{ii}}$, $h_{ii}$ being the $i$th diagonal element of $\mathbf{X}\left(\mathbf{X}^T\mathbf{X}\right)^{-1}\mathbf{X}^T$. The resulting robust estimator is taken as the initial estimator in (\ref{eq:iterboot}) to get a robust IB-estimator that we call IB-ROB. The resulting IB-estimator is in fact robust in the sense that it has a bounded influence function \citep[see][]{HaRoRoSt:86,Hamp:74}. 

Both initial estimators are not consistent estimators, even if we expect the asymptotic bias to be very small (for small values of $\delta$ in (\ref{Eq_pseudo-val})), but both IB-estimators should have a reduced finite sample bias where, in addition, the IB-ROB is also robust to data contamination.

\subsubsection{Simulation Study}
\label{Sec_sim-logistic}

We perform a simulation study to validate the properties of the IB-MLE and IB-ROB and compare their finite sample performance to other well established estimators. In particular, as a benchmark, we also compute the MLE, the bias reduced MLE (MLE-BR) using the \texttt{brglm} function (with default parameters) of the \texttt{brglm} package in R \citep{brglm2017}, as well as the robust estimator (\ref{Eq_rob-glm}) using the \texttt{glmrob} function in R without data transformation (ROB). We consider four situations that can occur with real data which result from the combinations of balanced outcome classes (Setting I) and unbalanced outcome classes (Setting II) with and without data contamination. We also consider a large model with $p=200$ and chose $n$ as to provide EPV of respectively $5$ and $3.75$, which are below the usually recommended value of 10. The parameter values for the simulations are provided in Table \ref{tab:sim-logistic}. 

 \begin{table}[!hb]
     \centering
     \caption{Simulation settings for the LRM.}
     \begin{tabular}{lrr}
 \toprule
 Parameters &  Setting I & Setting II \\
 \midrule
 $p=$ & $200$ & $200$ \\
 $n=$ & $2000$ & $3000$ \\
 $\sum_{i=1}^ny_i\approx$ & $1000$ & $750$ \\
 EPV $\approx$ & 5 & 3.75 \\
 $H=$ & $500$ & $500$ \\
 $\beta_1=\beta_2=$ & $5$ & $5$ \\
 $\beta_3=\beta_4=$ & $-7$ & $-7$ \\
 $\beta_5=\ldots=\beta_{200}=$ & $0$ & $0$ \\
 $\delta=$ & $0.01$ & $0.01$ \\
 Simulations & $500$ & $500$ \\
 \bottomrule
     \end{tabular}
     \label{tab:sim-logistic}
 \end{table}

\begin{figure}[!ht]
     \centering
     \includegraphics[width=15cm]{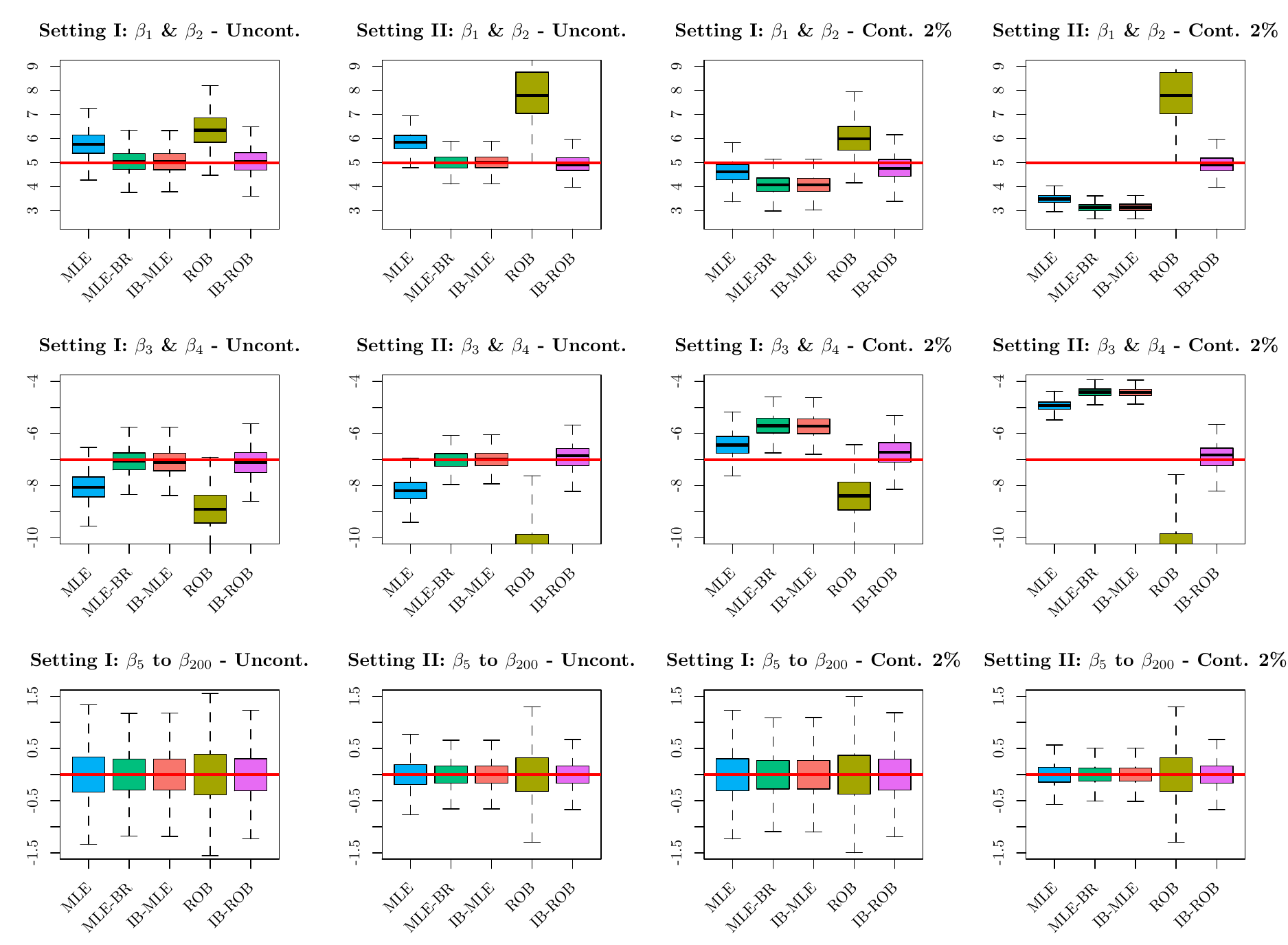}
     \caption{Finite sample distribution of estimators for the LRM using the simulation settings presented in Table \ref{tab:sim-logistic}. The estimators are the MLE (MLE), the Firth's bias reduced MLE (MLE-BR), the IB  with the MLE computed on the pseudo values \eqref{Eq_pseudo-val} as initial estimator (IB-MLE), the robust estimator in \eqref{Eq_rob-glm} (ROB) and the IB  with the robust estimator  computed on the pseudo values \eqref{Eq_pseudo-val} as initial estimator (IB-ROB). For each simulation setting, $1,000$ samples are generated.}
     \label{fig:sim-logistic-bxp}
 \end{figure}

\begin{figure}[!ht] 
     \centering
     \includegraphics[width=15cm]{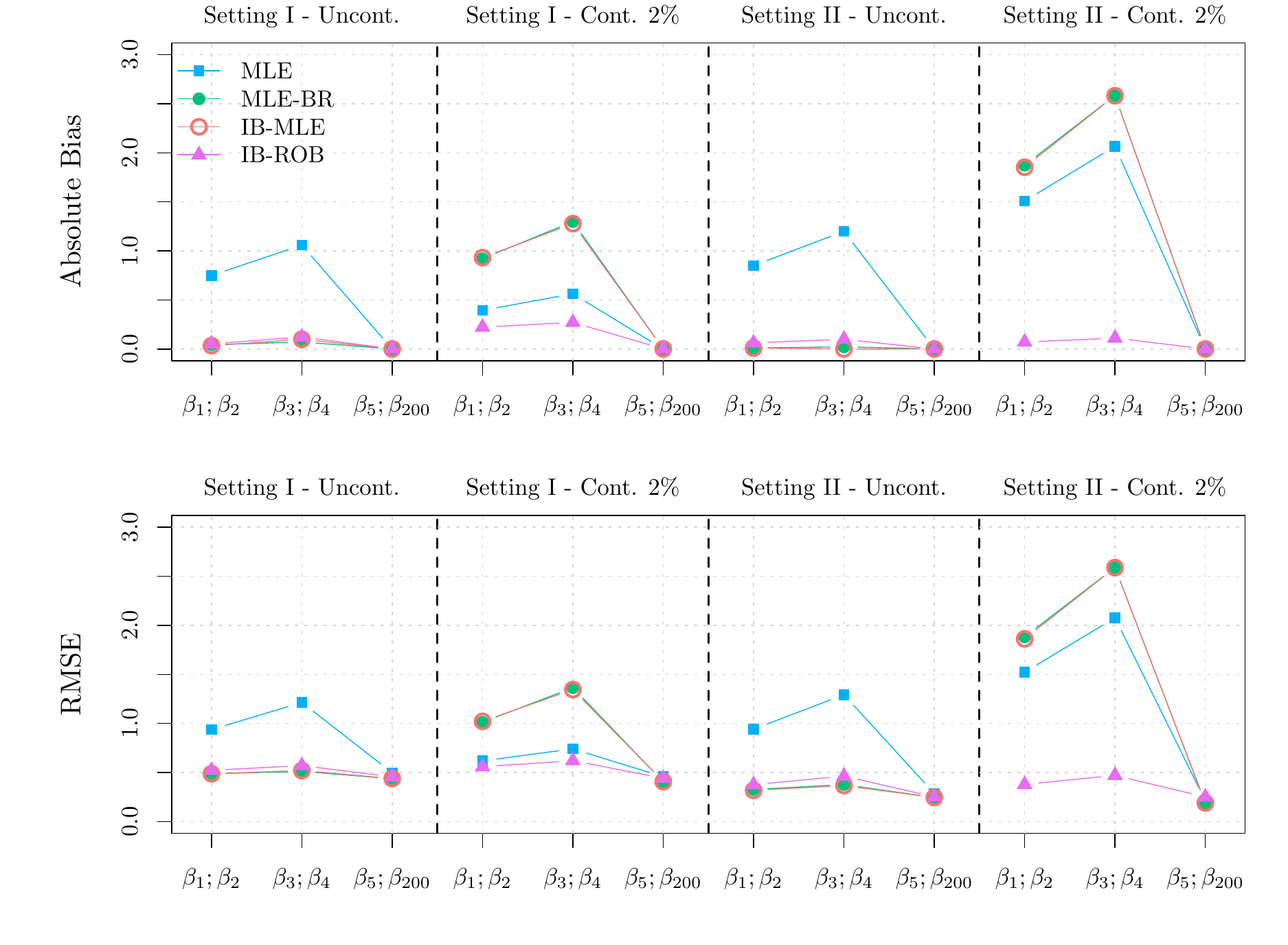}
     \caption{Finite sample bias and RMSE of estimators for the LRM using the simulation settings presented in Table \ref{tab:sim-logistic}. The estimators are the MLE (MLE), the Firth's bias reduced MLE (MLE-BR), the IB-estimator  with the MLE computed on the pseudo values (\ref{Eq_pseudo-val}) as initial estimator (IB-MLE) and the IB-estimator with the robust estimator in (\ref{Eq_rob-glm}) computed on the pseudo values (\ref{Eq_pseudo-val}) as initial estimator (IB-ROB). The bias and RMSE of the robust estimator in (\ref{Eq_rob-glm}) (ROB) was omitted in order to avoid an unsuitable scaling of the graphs. For each simulation setting, $1,000$ samples are generated.}
     \label{fig:sim-logistic-summary}
 \end{figure}

The covariates were simulated independently from distributions $\mathcal{N}(0,4/\sqrt{n})$ for Setting I and $\mathcal{N}(0.6,4/\sqrt{n})$ for Setting II, in order to ensure that the size of the log-odds ratio $\mathbf{x}_i\bbeta$ does not increase with $n$, so that $\mu_i(\bbeta)$ is not trivially equal to either $0$ or $1$. To contaminate the data, we chose a rather extreme misclassification error to show a noticeable effect on the different estimators, which consists in permuting 2\% of the responses with corresponding larger (smaller) fitted probabilities (expectations). The simulation results are presented in Figure \ref{fig:sim-logistic-bxp} as boxplots of the finite sample distribution, and in Figure \ref{fig:sim-logistic-summary} as the bias and Root Mean Squared Error (RMSE) of the different estimators.

 \clearpage

The finite sample distributions presented in Figure \ref{fig:sim-logistic-bxp}, as well as the summary statistics given by the bias and RMSE presented in Figure \ref{fig:sim-logistic-summary}, allow us to draw the following conclusions that support the theoretical results. In the uncontaminated case, as expected the MLE is biased (except when the slope parameters are zero) and the MLE-BR, IB-MLE and IB-ROB are all unbiased which is not the case for the robust estimator ROB. Moreover, the variability of all estimators is comparable, except for ROB which makes it rather inefficient in these settings. With 2\% of contaminated data (missclassification error), the only unbiased estimator is IB-ROB and its behaviour remains stable compared to the uncontaminated data setting. This is in line with a desirable property of robust estimators, that is stability with or without (slight) data contamination. The behaviour of all estimators remains the same in both settings, that is, whether or not the responses are balanced. Finally, as argued above, a better proposal for a robust, bias reduced and consistent estimator, as an alternative to IB-ROB, could in principle be proposed, but this is left for further research. 

\subsection{Bias Corrected Estimator for the Random Intercept Logistic Regression Model}
\label{sec:Mlogistic}

An interesting way of extending the LRM to account for the dependence structure between the observed responses is to use the GLMM family \citep[see e.g.][]{NeWe:72,BrCl:93,LeNe:01,McCuSe:01,JiangBook2007}. We consider here the case of the random intercept model, which is a suitable model in many applied settings. The binary response is denoted by $\Y_{ij}(\bm{\beta}_0, n)$, where $i=1,\ldots,m$ and $j=1,\ldots,n_i$. The expected value of the response is expressed as
\begin{equation}
    \boldsymbol{\mu}_{ij}(\beta_0,\bm{\beta}|U_i)\equiv\mathbb{E}\left[\Y_{ij}(\beta_0,\bm{\beta}, n)|U_i\right] = \frac{\exp{\left(\beta_0 + \mathbf{x}_{ij}^T\boldsymbol{\beta}+U_i\right)}}{1+\exp{\left(\beta_0 + \mathbf{x}_{ij}^T\boldsymbol{\beta}+U_i\right)}},
    \label{Eq_linpredGLMM}
\end{equation}
where $\beta_0$ is the intercept,
$\mathbf{x}_{ij}$ is a $q$-vector of covariates,
$\boldsymbol{\beta}$ is a $q$-vector of regression coefficients and the random effect  $
U_i,\; i = 1, ..., m$ is a normal random variable with zero mean and (unknown) variance $\sigma^2$.

Because the random effects are not observed, the MLE is derived on the marginal likelihood function
where the random effects are integrated out. These integrals have no known closed-form solutions, so approximations to the (marginal) likelihood function have been proposed, including the pseudo- and Penalized Quasi- Likelihood (PQL) \citep[see][]{Scha:91,WoOCo:93,BrCl:93}, LAplace approximations (LA) \citep{RaYaYo:00,huber2004estimation} and adaptive Gauss-Hermite Quadrature (GHQ) \citep{PiCh:06}. It is well known that PQL methods lead to biased estimators while LA and GHQ are more accurate \citep[for extensive accounts of methods across software and packages, see e.g.][]{BOLKER2009127,KiChEm:13}. The \texttt{lme4} R package \citep{BaMaBo:10} uses both the LA and GHQ to compute the likelihood function, while the \texttt{glmmPQL} R function of the \texttt{MASS} library \citep[that supports][]{VeRi:02} uses the PQL.

\begin{table}[!hb]
     \centering
     \caption{Simulation settings for the LRM with a random intercept.}
     \begin{tabular}{lrr}
 \toprule
 Parameters &  Setting I & Setting II \\
 \midrule
 $p=q + 2=$ & $31$ & $31$ \\
 $m=$ & $5$ & $50$ \\
 $\forall i,\;n_i=n_\circ =$ & $50$ & $5$\\
 $n=$ & $250$ & $250$ \\
 $\sum_{i=1}^m\sum_{j=1}^{n_i} y_{i,j}\approx$ & $125$ & $125$   \\
 EPV $\approx$ & $4$ & $4$ \\
 $H=$ & $200$ & $200$ \\
 $\beta_0=$ & $0$ & $0$ \\
 $\beta_1=\beta_2=$ & $5$ & $5$ \\
 $\beta_3=\beta_4=$ & $-7$ & $-7$ \\
 $\beta_5=\ldots=\beta_{30}=$ & $0$ & $0$ \\
 $\sigma^2=$ & $1.5$ & $1.5$ \\
 $\delta=$ & $0.01$ & $0.01$ \\
 Simulations & $1,000$ & $1,000$ \\
 \bottomrule 
     \end{tabular}
     \label{tab:sim-logistic2}
 \end{table}

In this section, we propose an IB-estimator that has a reduced finite sample bias compared to the (approximated) MLE using either the LA, GHQ or PQL approaches. More precisely, as an initial estimator, for computational efficiency, we consider an estimator defined through penalized iteratively reweighted least squares as described, for example, in \cite{BaMaBo:10} and implemented in the \texttt{glmer} function (with argument \texttt{nAGQ} set to 0) of the \texttt{lme4} R package. As done in Section \ref{Sec_IB-logistic} for the LRM, we transform the observed responses to get pseudo-values as in (\ref{Eq_pseudo-val}) with $\delta=0.01$. Hence, the initial estimator is a less accurate approximation to the MLE (compared e.g. to the LA or the GHQ approximations) although it is not affected by data separation which makes it asymptotically biased. However, we expect the IB to correct the bias induced through the transformed responses.

To study the behaviour of the IB-estimator and compare its performance in terms of bias and variance in finite samples to different approximations of the MLE, we perform a simulation study using the two settings described in Table \ref{tab:sim-logistic2}. Both settings can be considered as high dimensional in the sense that $p$ is relatively large with respect to $m$. Moreover, while Setting II ($m=50$, $n_i=n_\circ = 5, \forall i$) reflects a possibly more frequent situation, Setting I concerns the case where $m$ is small (and much smaller than $n_\circ$), a situation frequently encountered in cluster randomised trials  \citep[see e.g.][and the references therein]{CRT-Huang-2016,CRT-Leyrat-2017}. As for the simulation study in Section \ref{Sec_sim-logistic} studying the LRM, the covariates are simulated independently from the distribution $\mathcal{N}(0,4/\sqrt{n})$.   

\begin{figure}[!ht]
     \centering
     \includegraphics[width=15cm]{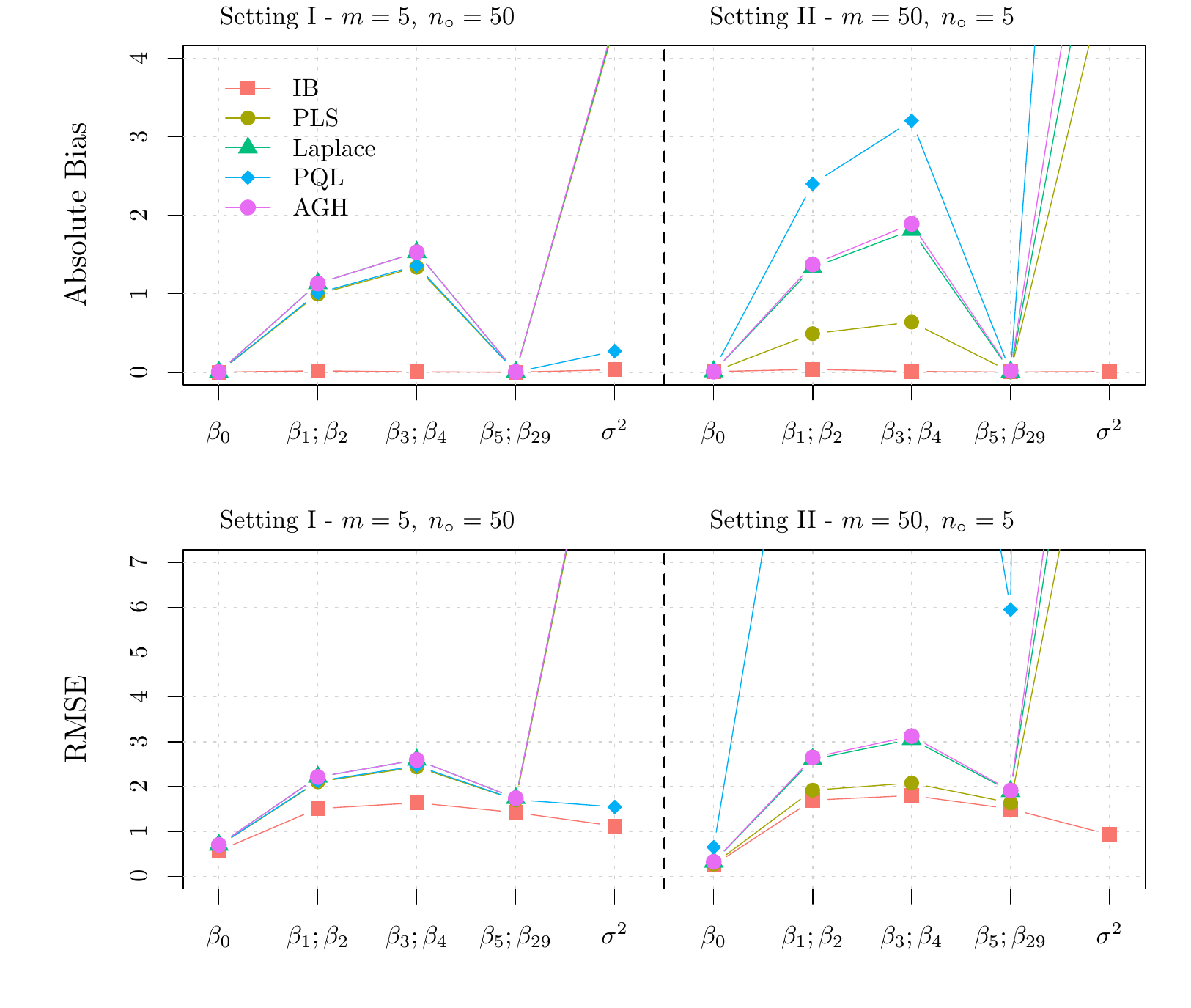}
      \caption{Finite sample bias and RMSE of estimators for the LRM with a random intercept, using the simulation settings presented in Table \ref{tab:sim-logistic2}. The estimators are the MLE with Laplace approximation (MLE-LA), the MLE with adaptive Gauss-Hermite quadratures (MLE-GHQ), the PQL (MLE-PQL) and the IB-estimator  with intital estimator explained in the main text (IB). For each simulation setting, $1,000$ samples are generated.}
     \label{fig:sim-glmm-summary}
 \end{figure}

Figure \ref{fig:sim-glmm-summary} presents the finite sample bias and RMSE of the approximated MLE estimators and the IB-estimator. It would appear evident that the proposed IB-estimator drastically reduces (removes) the finite sample estimation bias, especially for the random effect variance estimator. Moreover, the IB-estimator also achieves the lowest RMSE. These simulation results confirm the theoretical properties of the IB-estimator, namely finite sample (near) unbiasedness, which has important advantages when performing inference in practice, and/or when the parameter estimates, such as the random intercept variance estimate, are used, for example, to evaluate the sample size needed in subsequent randomized trials. The IB-estimator can eventually be based on another initial estimators in order to possibly improve efficiency for example, however this study is left for future research.

 \clearpage

\newpage
\bibliographystyle{agsm}
\bibliography{biblio}

\newpage
\phantom{a}
\newpage


\appendix
\titleformat{\section}[hang]{\large\center\scshape}{{\sc Appendix} \Alph{section}:}{1em}{}

\section{Notation and Organisation of the Appendices}\label{app:notations}

As mentioned in Remark \ref{Remark:notation}, in order to avoid potential confusing situations, we will use more precise notations in the sequel. Namely, we make our notations more precise and we use, as presented in Table  \ref{tab:notations}, a dictionary between the notations used in the main text and in the appendices. The justification for these precision in notation is given in  Appendix \ref{app:gen:setting}. 

\begin{table}[!hb]
    \centering
     \caption{Dictionary between notations used in the main text and the appendices.}
    \begin{tabular}{rcl}
    \toprule
Main &   & Appendices \\
\midrule
$\X(\bto,n)$ & $\equiv$ & $\X(\bto,n,\bw_0)$    \\
$\X^*_i(\bt,n)$ & $\equiv$ & $\X(\bt,n,\bw_i)$ \\
$\hbpi(\bto,n)$ & $\equiv$ & $\hbpi(\bto,n,\bw_0)$  \\
$\hbpi^*_h(\bt,n)$ & $\equiv$ & $\hbpi(\bto,n,\bw_h)$  \\
$\hbtheta^{(k)}$ & $\equiv$ & $\hbt^{(k)}$  \\
$\hbtheta$ & $\equiv$ & $\hbt$ \\
$\widehat{\bT}$ & $\equiv$ & $\widehat{\bT}_{(n,H)}$ \\
$\mathbf{v}(\btheta,n)$ & $\equiv$ & $\mathbf{v}(\btheta,n,\bw)$ \\
\bottomrule
    \end{tabular}
    \label{tab:notations}
\end{table}

Another notable difference between the main text and the appendices is that the results presented in the latter are proved independently in the following appendices. Indeed, in Appendix \ref{app:conv:IB} we prove the converenge of the IB sequence, in Appendix \ref{app:unbias:IB} the unbiasedness of the IB-estimator, in Appendix \ref{app:consist:IB} its consistency and finally in Appendix \ref{app:asymp:norm:IB} we state and prove its asymptotic normality. In each Appendix, the assumptions are tailored to the result we aim to prove resulting in a consequent number of assumptions. For this reason, an assumption framework is defined in Appendix \ref{app:gen:frame} to classify assumptions of the same type. In the last one, Appendix \ref{app:main:res}, these assumptions are combined to prove the main results of the article, namely Theorem \ref{thm:conv:consit:iter:boot} and Theorem \ref{THM:bias} that are stated in the main text, as well as Theorem \ref{app:thm:main:3} which formally states the asymptotic normality of the IB-estimator.

Although the results presented in Appendices \ref{app:conv:IB} to \ref{app:unbias:IB} are closely related to the ones of Theorem \ref{thm:conv:consit:iter:boot} and \ref{THM:bias}, we nevertheless provide a detailed discussion that is, in parts, redundant, but makes it complete within the material provided in these appendices.

\newpage
\setcounter{equation}{0}
\renewcommand{\theequation}{B.\arabic{equation}}
\section{Mathematical Setup}
\label{app:gen:setting}

In this appendix, we redefine the mathematical setup presented in the main text  in order to make it more precise, by using the notations presented in Table \ref{tab:notations}. Let us define $\mathbf{X}\left(\bm{\theta}, n, \bm{\omega}\right)\in\real^n$ as being a random sample generated under model $F_{\bm{\theta}}$ (possibly conditional on a set of fixed covariates), where $\bm{\theta}\in\bm{\Theta}\subset\real^p$ is the parameter vector of interest and $\bm{\omega}\in\real^m,\,m\geq n > 0$, represents a random variable explaining the source of randomness of the sample. More specifically, $\bm{\omega}$ can be considered as a random \emph{seed} that produces, conditionally on its value and given a value of $\bm{\theta}$, the sample $\mathbf{X}\left(\bm{\theta}, n, \bm{\omega}\right)$ of size $n$ in a deterministic manner. Indeed, $\bm{\omega}$ can be conceived as a random variable issued from a model $G$, thereby justifying the definition of the random sample as $\mathbf{X}\left(\bm{\theta}, n, \bm{\omega}\right)$. With this in mind, denoting ``$\overset{d}{=}$'' as ``equality in distribution'', there is no requirement for $\bm{\omega}$ to be unique since it is possible to have $\mathbf{X}\left(\bm{\theta}, n, \bm{\omega}\right)\overset{d}{=}\mathbf{X}\left(\bm{\theta}, n, \bm{\omega}^\ast\right)$ even though $\bm{\omega}\overset{d}{\neq}\bm{\omega}^\ast$. For simplicity, throughout this article we will consider $\bm{\omega}$ as only belonging to a fixed set $\bOmega \equiv \left\{\bm{\omega}_{h}\right\}_{h=1,\dots,H} \subset \real^m$. In this setting, defining $\bm{\theta}_0 \in \bm{\Theta}$, we let $\bm{\omega}_0$ and $\bm{\theta}_0$ denote respectively the unknown \textit{random} vector and \textit{fixed} parameter vector used to generate the random sample $\mathbf{X}\left(\bm{\theta}_0, n, \bm{\omega}_0\right)$ that will be used to estimate $\bm{\theta}_0$. Knowing that the sample $\mathbf{X}\left(\bm{\theta}_0, n, \bm{\omega}_0\right)$ is generated from a random seed $\bm{\omega}_0$, we can also consider other samples that we denote as $\mathbf{X}\left(\bm{\theta}, n, \bm{\omega}\right)$, where $(\bm{\theta}, n, \bm{\omega}) \in \bm{\Theta} \times \mathbb{N}^\ast \times \bm{\Omega}$, which can therefore be simulated based on different values of $\bm{\omega}$ (and $\bm{\theta})$. Notice that, based on this premise, we have that $\bm{\omega}_0 \not\in \bOmega$ which therefore implies that, in the context of this article, it is not possible to generate samples using $\bm{\omega}_0$. Hence, the use of $\bm{\omega}$ allows to explicitly display the randomness of a sample. In addition, with this notation it is not only possible to clearly distinguish the observed sample and the simulated ones but also to define the difference (or equivalence) between two simulated samples, say  $\mathbf{X}\left(\bm{\theta}, n, \bm{\omega}_j\right)$ and  $\mathbf{X}\left(\bm{\theta}, n, \bm{\omega}_l\right)$.  

Now, we focus on the estimation procedure aimed at targeting the value of $\bm{\theta}_0$ from $\mathbf{X}\left(\bm{\theta}_0, n, \bm{\omega}_0\right)$. For this purpose, let us define $\hat{\bm{\pi}}\left(\mathbf{X}\left(\bm{\theta}_0, n, \bm{\omega}_0\right)\right)$ as a biased estimator of $\bm{\theta}_0$ which, despite being possibly inconsistent, is either readily available or can easily be computed. For simplicity of notation, from this point onward we will refer to the estimator $\hat{\bm{\pi}}\left(\mathbf{X}\left(\bm{\theta}_0, n, \bm{\omega}_0\right)\right)$ as $\hat{\bm{\pi}}\left(\bm{\theta}_0, n, \bm{\omega}_0\right)$. Similarly, to the definitions of Section \ref{sec:setting} we define, for all $(n, H) \in  \mathbb{N}^\ast \times \mathbb{N}^\ast$, the following set
\begin{equation}
    \widehat{\bT}_{(n,H)} \equiv \argzero_{\bm{\theta} \in
      \bm{\Theta}} \;  \hat{\bm{\pi}}(\bm{\theta}_0, n, \bm{\omega}_0) - \frac{1}{H} \sum_{h = 1}^H  \hat{\bm{\pi}}(\bm{\theta}, n, \bm{\omega}_{h}),
    \label{app:eq:indirectInf:hTimesN}
\end{equation}
where $\hat{\bm{\pi}}(\bm{\theta}, n, \bm{\omega}_{h})$ is obtained on $\mathbf{X}\left(\bm{\theta}, n, \bm{\omega}_h \right)$. This definition assumes that the solution set $\widehat{\bT}_{(n,H)}$ is not empty for all $( n, H) \in \mathbb{N}^\ast \times \mathbb{N}^\ast$, which is reasonable when $\dim(\hat{\bm{\pi}}(\bm{\theta}, n, \bm{\omega}))= \dim(\bm{\theta})$. 
In similar fashion, the IB sequence $\left\{\hbt^{(k)}\right\}_{k \in \mathbb{N}}$ can defined using the notation of this appendix as
\begin{equation}
	\hbt^{(k)} \equiv \hbt^{(k-1)} + \left[ \hat{\bm{\pi}}(\bm{\theta}_0, n, \bm{\omega}_0) - \frac{1}{H} \sum_{h = 1}^H  \hat{\bm{\pi}}\left(\hbt^{(k-1)}, n, \bm{\omega}_{h}\right) \right],
	\label{app:eq:iterboot}
\end{equation}
where $\hbt^{(0)}\in\bm{\Theta}$. When this iterative procedure converges, we define the IB-estimator $\hbt\in\real^p$ as the limit in $k$ of $\hbt^{(k)}$.

\newpage
\setcounter{equation}{0}
\renewcommand{\theequation}{C.\arabic{equation}}
\section{Assumptions and Asymptotic Framework}
\label{app:gen:frame}

Various assumptions are considered and, in order to provide a general reference for the type of assumption, we employ the following conventions to name them:
\begin{itemize}
    \item Assumptions indexed by ``A'' refer to the topology of $\bm{\Theta}$;
    \item Assumptions indexed by ``B'' refer to the existence of the expected value of the initial estimator;
    \item Assumptions indexed by ``C'' refer to the random component of the initial estimator;
    \item Assumptions indexed by ``D'' refer to the bias of the initial estimator.
\end{itemize}

When needed, sub-indices or upper-indices will be used to distinguish assumptions of the same type (for example Assumptions \ref{assum:A:1} and \ref{assum:A'} are different but both related to the topology of $\bm{\Theta}$).
 
As it was mentioned in the main text, the asymptotic framework we use here is somewhat unusual as we always consider arbitrarily large but finite $n$ and $p$. Indeed, if $p$ is such that $p \to \infty$ as $n \to \infty$, our asymptotic results may not be valid when taking the limit in $n$ (and therefore in $p$) but are valid for all finite $n$ and $p$. The difference between the considered framework and others where limits are studied is rather subtle as, in practice, asymptotic results are typically used to derive approximations in finite samples for which the infinite (dimensional) case is not necessarily informative. Indeed, the topology of $\real^p$ for any finite $p$ and $\real^\infty$ are profoundly different. For example, the consistency of an estimator is generally dependent on the assumption that $\bm{\Theta}$ is compact. In the infinite dimensional setting closed and bounded sets are not necessarily compact. Therefore, this assumption becomes rather implausible for many statistical models and would imply, among other things, a detailed topological discussion of requirements imposed on $\bm{\Theta}$. Similarly, many of the mathematical arguments presented in this article may not apply in the infinite dimensional case. Naturally, when $p \to c <\infty$ as $n \to \infty$, limits in $n$ are allowed. Although not necessary, we abstain from using statements or definitions such as ``$\,\xrightarrow{p}$'' (convergence in probability) or ``$\,\xrightarrow{d}\,$'' (convergence in distribution) which may lead one to believe that the limit in $n$ exists. This choice is simply made to avoid confusion.

\newpage
\setcounter{equation}{0}
\renewcommand{\theequation}{D.\arabic{equation}}
\section{Converenge of the IB Sequence}\label{app:conv:IB}

In this appendix, we define new assumptions that are weaker than the ones presented in the main text. These assumptions are used to study the properties of the IB sequence $\hbtseq$ defined in (\ref{app:eq:iterboot}). Our first assumption concerns the topology of the parameter space $\bm{\Theta}$. 
	
 \setcounter{Assumption}{0}
    \renewcommand\theAssumption{\Alph{Assumption}$_1$}
    \begin{Assumption}
		\label{assum:A:1}
		Let $\bm{\Theta}$ be such that
		\begin{equation*}
			\argzero_{\bm{\theta} \in
		      \real^p \setminus \bm{\Theta}} \;  \hat{\bm{\pi}}(\bm{\theta}_0, n, \bm{\omega}_0) - \frac{1}{H} \sum_{h = 1}^H  \hat{\bm{\pi}}(\bm{\theta}, n, \bm{\omega}_{h})
		      = \emptyset\,.
		\end{equation*}
\end{Assumption}
\vspace{0.25cm}
Assumption \ref{assum:A:1} is clearly weaker than Assumption \ref{assum:A}. It is very mild and essentially ensures that no solution outside of $\bm{\Theta}$ exists for the IB-estimator defined in (\ref{app:eq:indirectInf:hTimesN}). 
In addition, Assumption \ref{assum:A:1} doesn't require $\bm{\Theta}$ to be a compact set since, for example, it can represent an open set. Next, we simply restate Assumption \ref{assum:B} with the notation used in the appendices.
    \setcounter{Assumption}{1}
    \renewcommand{\theHAssumption}{otherAssumption\theAssumption}
    \renewcommand\theAssumption{\Alph{Assumption}}
\begin{Assumption}
	\label{assum:B'}
	For all $\left(\bm{\theta}, \, n, \, \bm{\omega}\right) \in \bm{\Theta} \times \mathbb{N}^\ast \times \bOmega$, the expectation
		$
			\bm{\pi} \left( \bm{\theta}, n\right) \equiv \mathbb{E}\left[\hat{\bm{\pi}}(\bm{\theta},n, \bm{\omega})\right] 
		$
		exists and is finite, i.e. $\lvert\bm{\pi}_j \left(\bm{\theta}, n  \right)\rvert < \infty$ for all $j = 1, \ldots, p$.
		Moreover, for all $j = 1, \ldots, p$ $\bm{\pi}_j \left( \bm{\theta}\right) \equiv \displaystyle{\lim_{n \to \infty}} \bm{\pi}_j \left( \bm{\theta}, n\right)$ exists.
\end{Assumption}
\vspace{0.25cm}
Using Assumption \ref{assum:B'}, we can write:
	\begin{equation}
	   \hat{\bm{\pi}}(\bm{\theta},n,\bm{\omega})  = \bm{\pi} \left( \bm{\theta}, n\right)+ \mathbf{v} \left(
	    \bm{\theta},n, \bm{\omega}\right),
	  \label{app:bias:estim}
	\end{equation}
where $\mathbf{v} \left(\bm{\theta},n, \bm{\omega}\right) \equiv \hat{\bm{\pi}}(\bm{\theta},n,\bm{\omega}) - \bm{\pi} \left( \bm{\theta}, n\right)$ is a zero-mean random vector. 

Similarly, we restate Assumption \ref{assum:C} using the notation of the appendices. 
    \setcounter{Assumption}{2}
    \renewcommand{\theHAssumption}{otherAssumption\theAssumption}
    \renewcommand\theAssumption{\Alph{Assumption}}
	\begin{Assumption}
		\label{assum:C'}
		The random vector $\mathbf{v}(\bm{\theta},n, \bm{\omega})$ is such that its second moment exists and such that there exists a real $\alpha > 0$ such that for all $\left(\bm{\theta}, \, \bm{\omega}\right) \in \bm{\Theta} \times \bOmega$ and all $j = 1,\,\ldots,\,p$, we have
		\begin{equation*}
			\mathbf{v}_j(\bm{\theta},n, \bm{\omega}) = \mathcal{O}_{\rm p}(n^{-\alpha})\;\;\;\;
			 \text{and}\;\;\;\;
			 \lim_{n \to \infty} \; \frac{p^{\nicefrac{1}{2}}}{n^\alpha} = 0.
		\end{equation*}
	\end{Assumption}
	\vspace{0.25cm}
	
Next, we adapt the bias terms to the notation used in the appendices. Indeed, we consider the bias of $\hat{\bm{\pi}}(\bm{\theta},n,\bm{\omega})$ for $\bm{\theta}$ and we let $\mathbf{d}\left(\bm{\theta}, n\right) \equiv \mathbb{E} \left[\hat{\bm{\pi}}(\bm{\theta},n,\bm{\omega})\right] - \bm{\theta}$, allowing to write
	\begin{equation*}
	   \hat{\bm{\pi}}(\bm{\theta},n,\bm{\omega})  = \bm{\theta} + \mathbf{d} \left( \bm{\theta}, n\right)+ \mathbf{v} \left(
	    \bm{\theta},n, \bm{\omega}\right).
	\end{equation*}
	The functions $\a{\bt}, \mathbf{b}(\bt,n), \cn, \Ln$ and $\r{\bt}$ are defined in the same way as in \eqref{def:fct:bias:d} and \eqref{def:fct:bias:b}.
	Assumption \ref{assum:D:1} (below) imposes restrictions on the bias function that are weaker than the ones proposed in Assumption \ref{assum:D}. 
	
    \setcounter{Assumption}{3}
    \renewcommand\theAssumption{\Alph{Assumption}$_1$}
	\begin{Assumption}
		\label{assum:D:1}
		The bias function $\mathbf{d}\left(\bm{\theta}, n\right)$ is such that:
		\begin{enumerate}
		    \item The function $\mathbf{a}(\bm{\theta})$ is a contraction map in that for all
		    $\bm{\theta}_1,\bm{\theta}_2\in\bm{\Theta}$ such that $\bm{\theta}_1 \neq \bm{\theta}_2$ we have 
		    \begin{equation*}
		        \Big\lVert \mathbf{a}(\bm{\theta}_2)-\mathbf{a}(\bm{\theta}_1)  \Big\rVert_2 <  \big\rVert \bm{\theta}_2-\bm{\theta}_1 \big\lVert_2\,.
		    \end{equation*}
		    
		\item There exist real $\beta, \gamma > 0$ such that for all $\bm{\theta} \in \bm{\Theta}$ and all $j,l = 1,\,\ldots, \, p$, we have
		\begin{equation*}
		\begin{aligned}
			  \mathbf{L}_{j,l}(n) = \mathcal{O}(n^{-\beta}), \;\;\;  \mathbf{r}_j\left(\bm{\theta}, n\right) = \mathcal{O}(n^{-\gamma}),\;\;\; 
			  \lim_{n \to \infty} \; \frac{p}{n^\beta} = 0,  \;\;\;\ \text{and} \;\;\;
			  \lim_{n \to \infty} \; \frac{p^{\nicefrac{1}{2}}}{n^\gamma} = 0.
		 \end{aligned}
		\end{equation*}
		\end{enumerate}
	\end{Assumption}
    \vspace{0.25cm}
    
    There are two differences between Assumptions \ref{assum:D} and \ref{assum:D:1}. The first one is that Assumption \ref{assum:D} requires 
    \begin{equation*}
        \lim_{n \to \infty}\; p^{3/2}/n^{\beta} = 0,
    \end{equation*}
    whereas in Assumption \ref{assum:D:1} this requirement is simply
    \begin{equation*}
        \lim_{n \to \infty}\; p/n^{\beta} = 0.
    \end{equation*}
    This modification of Assumption \ref{assum:D} implies that under the form of the bias postulated in (\ref{Eq_bais-standard}), where we have that $\beta, \gamma \geq 1$ and for $\sqrt{n}$-consistent initial estimator, where we have $\alpha = \nicefrac{1}{2}$, the requirements of Assumptions \ref{assum:C'} and \ref{assum:D:1}	are satisfied if
	\begin{equation*}
	    \lim_{n \to \infty} \; \frac{p}{n} = 0.
	    \label{app:eq:prat:sett:D1}
	\end{equation*}
    This is a weaker condition than (\ref{equ:p/n}), obtained under Assumption \ref{assum:D}. The second difference between Assumptions \ref{assum:D} and \ref{assum:D:1} is that the requirement on the term $\mathbf{c}(n)$ are not needed.
    
We now study the convergence properties of the IB sequence defined in (\ref{app:eq:iterboot}) when used to obtained the IB-estimator presented in (\ref{app:eq:indirectInf:hTimesN}). In Lemma \ref{lemma:unique:iter:boot}, we show that when $n$ is sufficiently large, the solution space $\widehat{\bT}_{(n,H)}$ contains only one element. In other words, this result ensures that the function $\frac{1}{H} \sum_{h = 1}^H  \hat{\bm{\pi}}(\bm{\theta}, n, \bm{\omega}_{h})$  in (\ref{app:eq:indirectInf:hTimesN}) is injective for fixed (but possibly large) sample size $n$ and fixed $H$. Lemma \ref{lemma:unique:iter:boot} formally states this result.

	\begin{customlemma}{2}
		\label{lemma:unique:iter:boot}
		Under Assumptions \ref{assum:A:1}, \ref{assum:B'}, \ref{assum:C'} and \ref{assum:D:1}, for all $H \in \mathbb{N}^\ast$ and $n \in \mathbb{N}^\ast$ sufficiently large, the solution space $\widehat{\bT}_{(n,H)}$ is a singleton.
	\end{customlemma}
	\vspace{0.25cm}
	
	\begin{proof}
	We define the function 
		\begin{equation}
			\begin{aligned}
				T_{(n,H)} \;\;:\;\; & \bm{\Theta}  & \longrightarrow  &\;\;\;\; \real^p\\
				& \bm{\theta}  & \longmapsto & \;\;\;\; T_{(n,H)}(\bm{\theta}),\\
			\end{aligned}
			\label{eq:define:T}
		\end{equation}
		where 
		\begin{equation*}
			T_{(n,H)}(\bm{\theta}) \equiv \bm{\theta} + \hat{\bm{\pi}}(\bm{\theta}_0, n, \bm{\omega}_0) - \frac{1}{H} \sum_{h = 1}^H  \hat{\bm{\pi}}(\bm{\theta}, n, \bm{\omega}_{h}).
		\end{equation*}
		We define the set of fixed points of $T_{(n,H)}$ as follows
		\begin{equation*}
			\widetilde{\bm{\Theta}}_{(n,H)} = \left\{\bm\theta \in \bm{\Theta} \;\; \Big| \;\;T_{(n,H)}(\bm{\theta}) = \bm{\theta} \right\}.
		\end{equation*}
		Clearly, we have that 
		\begin{equation}
			\widetilde{\bm{\Theta}}_{(n,H)} = \widehat{\bT}_{(n,H)},
			\label{eq:setfixed_point}
		\end{equation}
		as defined in (\ref{app:eq:indirectInf:hTimesN}). In the remainder of the proof we will show that an extension of the function $T_{(n,H)}$ admits a unique fixed point in $\real^p$ by applying Kirszbraun theorem (see \citealp{federer2014geometric}) and Banach fixed-point theorem. Then, using (\ref{eq:setfixed_point}) and Assumption \ref{assum:A:1} we will be able to conclude that the set $\widehat{\bT}_{(n,H)}$ only contains this fixed point. 
		Let us start by using (\ref{app:bias:estim}) and Assumptions \ref{assum:B'}, \ref{assum:C'} and \ref{assum:D:1}, we have that
		\begin{equation*}
			\begin{aligned}
	\hat{\bm{\pi}}(\bm{\theta}_0, n, \bm{\omega}_0) &= \bm{\pi} \left( \bm{\theta}_0, n\right)+ \mathbf{v} \left(\bm{\theta}_0,n, \bm{\omega}_0\right)
				=  \bm{\theta}_0 + \mathbf{a}(\btheta_0) + \mathbf{c}(n) + \mathbf{L}(n) \bm{\theta}_0 + \bm{\delta}^{(1)},
			\end{aligned}
		\end{equation*}
		where $\bm{\delta}_j^{(1)} = \mathcal{O}\left(n^{-\gamma}\right) +  \mathcal{O}_{\rm p}\left(n^{-\alpha}\right)$ for $j = 1,\,\ldots,\, p$ and
		\begin{equation*}
			\begin{aligned}
				\frac{1}{H} \sum_{h = 1}^H  \hat{\bm{\pi}}(\bm{\theta}, n, \bm{\omega}_{h}) &= \bm{\pi} \left( \bm{\theta}, n\right) + \frac{1}{H} \sum_{h = 1}^H  \mathbf{v}(\bm{\theta}, n, \bm{\omega}_{h})\\
				&= \bm{\theta} + \mathbf{a}(\btheta) + \mathbf{c}(n) + \mathbf{L}(n) \bm{\theta} + \bm{\delta}^{(2)},
			\end{aligned}
		\end{equation*}
		where $\bm{\delta}^{(2)}_j =	\mathcal{O}\left(n^{-\gamma}\right) + \mathcal{O}_{\rm p}\left(H^{-\nicefrac{1}{2}} n^{-\alpha}\right)$ for $j = 1,\,\ldots,\, p$. Therefore, we have
		\begin{equation*}
			T_{(n,H)}(\bm{\theta}) = \bm{\theta}_0 + \mathbf{a}(\btheta_0) - \mathbf{a}(\btheta) + \mathbf{L}(n) \left(\bm{\theta}_0 - \bm{\theta}\right) +  \bm{\delta}^{(1)} - \bm{\delta}^{(2)}.
		\end{equation*}
		Next, we consider the following quantity for $\bm{\theta}_1, \bm{\theta}_2 \in \bm{\Theta}$, 
		\begin{equation*}
			\begin{aligned}
				\left\lVert T_{(n,H)}(\bm{\theta}_1) - T_{(n,H)}(\bm{\theta}_2)\right\rVert_2^2 &= \left\lVert\mathbf{a}(\btheta_2) - \mathbf{a}(\btheta_1) + \mathbf{L}(n) \left(\bm{\theta}_2 - \bm{\theta}_1\right) +  \bm{\delta}^{(3)} \right\rVert_2^2\\
				 &\leq \left\lVert\mathbf{a}(\btheta_2) - \mathbf{a}(\btheta_1) \right\rVert_2^2 + \left\lVert\mathbf{L}(n) \left(\bm{\theta}_2 - \bm{\theta}_1\right) \right\rVert_2^2 + \left\lVert\bm{\delta}^{(3)} \right\rVert_2^2, 
			\end{aligned}
		\end{equation*}
		where $\bm{\delta}_j^{(3)} = \mathcal{O}\left(n^{-\gamma}\right) +  \mathcal{O}_{\rm p}\left(n^{-\alpha}\right)$ for $j = 1,\,\ldots,\, p$. Then, using the fact that 
		\begin{equation*}
		    \left\lVert\mathbf{L}(n) \left(\bm{\theta}_2 - \bm{\theta}_1\right) \right\rVert_2^2 \leq \left\lVert\mathbf{L}(n)^T \mathbf{L}(n)\right\rVert_F \,  \left\lVert\bm{\theta}_2 - \bm{\theta}_1 \right\rVert_2^2
		\end{equation*}
		where $\lVert\cdot\rVert_F$ is the Frobenius norm, we obtain 
		\begin{equation*}
			\begin{aligned}
				\left\lVert T_{(n,H)}(\bm{\theta}_1) - T_{(n,H)}(\bm{\theta}_2) \right\rVert_2^2 \leq \left\lVert\mathbf{a}(\btheta_2) - \mathbf{a}(\btheta_1) \right\rVert_2^2 + \left\lVert\mathbf{L}(n)^T \mathbf{L}(n)\right\rVert_F \,  \left\lVert\bm{\theta}_2 - \bm{\theta}_1 \right\rVert_2^2 + \left\lVert\bm{\delta}^{(3)} \right\rVert_2^2 .
			\end{aligned}
		\end{equation*}
		We now consider each term of the above equation separately. First, since $\bm{a}(\btheta)$ is a contraction map, there exists $\varepsilon\in (0,1)$ such that %
	    \begin{equation*}
	    \left\lVert\mathbf{a}(\btheta_2) - \mathbf{a}(\btheta_1)\right\rVert_2 \leq \varepsilon\left\lVert\btheta_2 - \btheta_1\right\rVert_2.   
	    \end{equation*}
		Secondly, we have
		\begin{equation*}
			\lVert\mathbf{L}(n)^T \mathbf{L}(n)\rVert_F \,  \left\lVert\bm{\theta}_2 - \bm{\theta}_1 \right\rVert_2^2 = \Delta \; \left\lVert\bm{\theta}_2 - \bm{\theta}_1 \right\rVert_2^2,
		\end{equation*}
		where $\Delta = \mathcal{O}\left(p^2 n^{-2\beta} \right)$. Indeed, by writing $\mathbf{B} = \mathbf{L}(n)^T \mathbf{L}(n)$ we obtain
		\begin{equation*}
			\begin{aligned}
				\Delta &= \lVert \mathbf{B} \rVert_F = \sqrt{\sum_{j = 1}^p \sum_{l = 1}^p B_{j,l}^2} \leq p \max_{j,l = 1,\, \ldots, \, p} |B_{j,l}| = p \max_{j,l = 1,\, \ldots, \, p} |\sum_{m = 1}^p \mathbf{L}_{m,j}(n) \mathbf{L}_{m,l}(n)|\\
	 & \leq p^2	\max_{j,l,m = 1,\, \ldots, \, p} |\mathbf{L}_{m,j}(n) \mathbf{L}_{m,l}(n)| = \mathcal{O}\left(p^2 n^{-2\beta}\right).
			\end{aligned}
	\end{equation*}
	Finally, we have
	\begin{equation*}
		\lVert\bm{\delta}^{(3)} \rVert_2^2 = \sum_{j = 1}^p \left(\bm{\delta}_j^{(3)}\right)^2 \leq p \max_{j = 1,\,\ldots,\,p} \left(\bm{\delta}_j^{(3)}\right)^2 = \mathcal{O}\left(p n^{-2 \gamma}\right) + \mathcal{O}_{\rm p}\left(p n^{-2 \alpha}\right).
	\end{equation*}
	%
	%
    By combining these results, we have
	\begin{equation*}
			\left\Vert T_{(n,H)}(\bm{\theta}_1) - T_{(n,H)}(\bm{\theta}_2) \right\Vert_2^2 \leq (\varepsilon^2+\Delta) \left\Vert\bm{\theta}_2 - \bm{\theta}_1 \right\Vert_2^2 + \mathcal{O}\left(p n^{-2 \gamma}\right) + \mathcal{O}_{\rm p}\left(p n^{-2 \alpha}\right).
	\end{equation*}
	Since $\alpha, \beta, \gamma> 0$ by Assumptions \ref{assum:C'} and \ref{assum:D:1}
	, for sufficiently large $n$ we have that for all $\bm{\theta}_1, \, \bm{\theta}_2 \in \bm{\Theta}$ 
	\begin{equation}
	\label{T:contract}
		\left\Vert T_{(n,H)}(\bm{\theta}_1) - T_{(n,H)}(\bm{\theta}_2) \right\Vert_2 <  \left\Vert\bm{\theta}_2 - \bm{\theta}_1\right\Vert_2.
	\end{equation}
   Using Kirszbraun theorem, we can extend $T_{(n,H)}$ to a contraction map from $\real^p$ to itself. Therefore, applying Banach fixed-point theorem, there exists a unique fixed point $\bt^\ast_{(n,H)} \in \real^p$ such that 
	%
	\begin{equation*}
		T_{(n,H)} \left(\bm{\theta}^\ast_{(j,n,H)}\right) = \bm{\theta}^\ast_{(j,n,H)}\,.
	\end{equation*}
	By Assumption \ref{assum:A:1} we have that $\bm{\theta}^\ast_{(j,n,H)} \in \bm{\Theta}$, which implies by (\ref{eq:setfixed_point}) that
	\begin{equation*}
		\widehat{\bT}_{(n,H)} = \left\{ \bm{\theta}^\ast_{(j,n,H)} \right\},
	\end{equation*}
	which concludes the proof. 
	
	\end{proof}

    Building on the identification result given in the above lemma, the following proposition states that the IB-estimator is the limit of the IB sequence.
	\begin{Proposition}
		\label{thm:iter:boot}
		Under Assumptions \ref{assum:A:1}, \ref{assum:B'}, \ref{assum:C'} and \ref{assum:D:1}, for all $ H \in \mathbb{N}^\ast$ and $n \in \mathbb{N}^\ast$ sufficiently large, the IB sequence $\hbtseq$ has a limit which is the IB-estimator, i.e.
		\begin{equation*}
    	 \lim_{k \to \infty} \; \hbt^{(k)} = \hbt \in \widehat{\bT}_{(n,H)}.
		\end{equation*}
		Moreover, there exists a real $ \epsilon \in (0, \, 1)$ such that for any $k \in \mathbb{N}^\ast$ 
		\begin{equation*}
			\left\lVert\hbt^{(k)} - \hbt\right\rVert_2  =\mathcal{O}_{\rm p}({p}^{\nicefrac{1}{2}}\,\epsilon^k).
		\end{equation*}
	\end{Proposition}
	\vspace{0.25cm}
	\begin{proof}
	We start by recalling that  $\left\{ \hbt^{(k)} \right\}_{k \in \mathbb{N}}$ is defined as
	\begin{equation*}
		\hbt^{(k)} = T_{(n,H)} \left(\hbt^{(k-1)}\right)
	\end{equation*}
	where $\hbt^{(0)} \in \bm{\Theta}$ and the function $T_{(n,H)}$ is defined in (\ref{eq:define:T}).
	Using the same arguments used in the proof of Lemma \ref{lemma:unique:iter:boot} we have that $T_{(n,H)}$ allows to apply Banach fixed-point theorem implying that 
	\begin{equation*}
		\lim_{k \to \infty}  \; \hbt^{(k)} = \hbt,
	\end{equation*}
	which concludes the first part of the proof.
	
	For the second part, consider a $k\in \mathbb{N}^\ast$. By the inequality (\ref{T:contract}) of the proof of Lemma \ref{lemma:unique:iter:boot}, there exists an $\epsilon \in (0, \, 1)$ such that 
	\begin{equation*}
\left\Vert T_{(n,H)}\left(\hbt^{(k)}\right) - T_{(n,H)}\left(\hbt\right) \right\Vert_2 \leq \epsilon \left\Vert\hbt^{(k)} - \hbt\right\Vert_2.
	\end{equation*}
	By construction of the sequence $\left\{ \hbt^{(k)} \right\}_{k \in \mathbb{N}}$ and since $\hbt$ is a fixed point of $T_{(n,H)}$, we can derive that
	\begin{equation*}
	    \left\Vert T_{(n,H)}\left(\hbt^{(k)}\right) - T_{(n,H)}\left(\hbt\right) \right\Vert_2 \leq \frac{\epsilon^k}{1-\epsilon} \left\Vert\hbt^{(0)} - \hbt^{(1)}\right\Vert_2  =\mathcal{O}_{\rm p}({p}^{\nicefrac{1}{2}}\epsilon^k),
	\end{equation*}
	which concludes the second part of the proof.
	\end{proof}
	
\newpage
\setcounter{equation}{0}
\renewcommand{\theequation}{E.\arabic{equation}}
\section{Consistency of the IB-Estimator}\label{app:consist:IB}

	In this appendix, we consider the consistency property of the IB-estimator and set weaker versions of Assumptions \ref{assum:A} and \ref{assum:D}, used in Theorem \ref{thm:conv:consit:iter:boot}.

	\setcounter{Assumption}{0}
	\renewcommand{\theHAssumption}{otherAssumption\theAssumption}
	\renewcommand\theAssumption{\Alph{Assumption}$_2$}
	\begin{Assumption}
	\label{assum:A:2}
	Preserving the same requirement given in Assumption \ref{assum:A:1}, we add the condition that $\bm{\Theta}$ is compact.
	\end{Assumption}
	\vspace{0.25cm}
	
	Assumption \ref{assum:A:2} ensures that $\bm{\Theta}$ is compact but remains weaker than Assumption \ref{assum:A}. We now introduce Assumption \ref{assum:D:2}.
	
	\setcounter{Assumption}{3}
	\renewcommand{\theHAssumption}{otherAssumption\theAssumption}
	\renewcommand\theAssumption{\Alph{Assumption}$_2$}
	\begin{Assumption}
	\label{assum:D:2}
	Preserving the same definitions given in Assumption \ref{assum:D:1} and defining $c_n \equiv \displaystyle{\max_{j=1,\dots,p}} \mathbf{c}_j(n)$, we require that:
	\begin{enumerate}
	    \item $\mathbf{a}(\btheta)$ is continuous and such that the function $\btheta + \mathbf{a}(\btheta)$ is injective,
	    
	    \item the constant $\beta$ and the sequence $\left\{c_n\right\}_{n\in\mathbb{N}^*}$ are such that
	\begin{equation*}
	   \lim_{n \to \infty} \; \frac{p^{\nicefrac{3}{2}}}{n^\beta} = 0, \;\;\;\;\text{and} \;\;\;\; \lim_{n \to \infty} \; p^{\nicefrac{1}{2}}c_n = 0.
	\end{equation*}
	\end{enumerate}
    The other requirements of Assumption \ref{assum:D:1} remain unchanged.
	\end{Assumption}
	\vspace{0.25cm}
	
	As stated in Assumption \ref{assum:D:2}, some requirements are the same as those in Assumption \ref{assum:D:1}, namely that there exist real $\beta, \gamma > 0$ such that for all $\bm{\theta} \in \bm{\Theta}$, we have, for all $j,l = 1,\,\ldots, \, p$,
		\begin{equation*}
		\begin{aligned}
			  \mathbf{L}_{j,l}(n) = \mathcal{O}(n^{-\beta}), \;\;\;  \mathbf{r}_j\left(\bm{\theta}, n\right) = \mathcal{O}(n^{-\gamma}) \;\;\;\ \text{and} \;\;\;
			  \lim_{n \to \infty} \; \frac{p^{\nicefrac{1}{2}}}{n^\gamma} = 0.
		 \end{aligned}
		\end{equation*}
		Compared to Assumptions \ref{assum:D} and \ref{assum:D:1}, Assumption \ref{assum:D:2} relaxes the condition on $\mathbf{a}(\btheta)$. Apart from this difference, Assumptions \ref{assum:D:2} and \ref{assum:D} impose the same requirements. Compared to Assumption \ref{assum:D:1}, Assumption \ref{assum:D:2} imposes stronger requirements on $\beta$ and finally imposes a condition on the vector $\mathbf{c}(n)$. Clearly, Assumption \ref{assum:A:2} implies Assumption \ref{assum:A:1} while, on the contrary, Assumptions \ref{assum:D:1} and \ref{assum:D:2} don't mutually imply each other. In the situation where $\alpha = \nicefrac{1}{2}$, $\beta = 1$, $\gamma = 2$ and $\mathbf{c}(n) = \mathbf{0}$, Assumption \ref{assum:D:2} (on which Proposition \ref{THM:consistency} is based) is satisfied if
	\begin{equation}
	    \lim_{n \to \infty} \; \frac{p^{\nicefrac{3}{2}}}{n} = 0.
	    \label{eq:n:p:app}
	\end{equation}
	\begin{Proposition}
	\label{THM:consistency}
	Under Assumptions \ref{assum:A:2}, \ref{assum:B'}, \ref{assum:C'} and \ref{assum:D:2}, $\hbt$ is a consistent estimator of $\bm{\theta}_0$ for all $H \in \mathbb{N}^\ast$ in that for all $\varepsilon > 0$ and all $\delta > 0$, there exists a sample size $n^\ast \in\mathbb{N}^\ast$ such for all $n \in \mathbb{N}^\ast$ satisfying $n \geq n^\ast$ we have:
	\begin{equation*}
	    \Pr \left(|| \hbt - \bm{\theta}_0 ||_2 \geq \varepsilon \right) \leq   \delta.
	\end{equation*}
	\end{Proposition}
	\vspace{0.25cm}
	
	\begin{proof}
    	This proof is directly obtained by verifying the conditions of Theorem 2.1 of \cite{newey1994large} on the functions $\hat{Q}(\bm{\theta}, n)$ and $Q(\bm{\theta})$ defined as follow:

	\begin{equation*}
	    \hat{Q}(\bm{\theta}, n) \equiv \Big\Vert \hat{\bm{\pi}}(\bm{\theta}_0, n, \bm{\omega}_0) - \frac{1}{H} \sum_{h = 1}^H  \hat{\bm{\pi}}(\bm{\theta}, n, \bm{\omega}_{h}) \Big\Vert_2,
	\end{equation*}
	and 
	\begin{equation*}
	     {Q}(\bm{\theta}) = \big\Vert {\bm{\pi}}(\bm{\theta}_0) -  {\bm{\pi}}(\bm{\theta}) \big\Vert_2.
	\end{equation*}
	Reformulating the requirements of this theorem to our setting, we have to show that  (i) $\bm{\Theta}$ is compact, (ii) ${Q}(\bm{\theta})$ is continuous, (iii) ${Q}(\bm{\theta})$ is uniquely minimized at $\btheta_0$, (iv) $\hat{Q}(\bm{\theta}, n)$ converges uniformly in probability to $Q(\bm{\theta})$.

	On one hand, Assumptions \ref{assum:A:2} ensures that $\bm{\Theta}$ is compact. On the other hand, Assumption \ref{assum:D:2} guarantees that ${Q}(\bm{\theta})$ is continuous and uniquely minimized at $\btheta_0$ since $\pi(\btheta)=\btheta+\mathbf{a}(\btheta)$ is required to be continuous and injective. What remains to be shown is that $\hat{Q}(\bm{\theta}, n)$ converges uniformly in probability to $Q(\bm{\theta})$, which is equivalent to show that: $\forall \varepsilon > 0$ and $\forall \delta > 0$, there exists a sample size $n^\ast \in\mathbb{N}^\ast$ such that for all $n \geq n^\ast$
	\begin{equation*}
	    \Pr \left( \sup_{\bm{\theta} \in \bm{\Theta}} \; \Big| \hat{Q}(\bm{\theta}, n) - Q(\bm{\theta}) \Big| \geq \varepsilon \right) \leq   \delta.
	\end{equation*}

	Fix $\varepsilon > 0$ and $\delta > 0$. Using the above definitions, we have that
	\begin{equation}
	    \sup_{\bm{\theta} \in \bm{\Theta}} \; \Big| \hat{Q}(\bm{\theta}, n) - Q(\bm{\theta}) \Big| \leq \sup_{\bm{\theta} \in \bm{\Theta}} \; \left[ \Big| \hat{Q}(\bm{\theta}, n) - Q(\bm{\theta}, n) \Big| + \Big| {Q}(\bm{\theta}, n) - Q(\bm{\theta}) \Big| \right],
	    \label{eq:convergen_proof_eq}
	\end{equation}
	where
	\begin{equation}
	        {Q}(\bm{\theta}, n) \equiv \big\Vert {\bm{\pi}}(\bm{\theta}_0, n) -  {\bm{\pi}}(\bm{\theta}, n) \big\Vert_2.
	\end{equation}
	We now consider each term of (\ref{eq:convergen_proof_eq}) separately. For simplicity, we define
		\begin{equation*}
	    \bar{\bm{\pi}}\left(\bm{\theta}, n, \bm{\omega}^{(H)}\right) \equiv \frac{1}{H} \sum_{h = 1}^H  \hat{\bm{\pi}}(\bm{\theta}, n, \bm{\omega}_{h}),
	\end{equation*}
	and, considering the first term of \eqref{eq:convergen_proof_eq}, we have
	\begin{equation*}
	   \begin{aligned}
	       \Big| \hat{Q}(\bm{\theta}, n) - Q(\bm{\theta}, n) \Big| 
	       &\leq  \big\Vert \hat{\bm{\pi}}(\bm{\theta}_0, n, \bm{\omega}_0) - {\bm{\pi}}(\bm{\theta}_0, n) \big\Vert_2 + \big\Vert {\bm{\pi}}(\bm{\theta}, n) - \bar{\bm{\pi}}\left(\bm{\theta}, n, \bm{\omega}^{(H)}\right) \big\Vert_2 \\
	       &= \big\Vert \mathbf{v} \left(\bm{\theta}_0, n, \bm{\omega}_0\right) \big\Vert_2 + \Big\Vert \frac{1}{H}\sum_{h=1}^H\mathbf{v} \left(\bm{\theta}, n, \bm{\omega}_{h}\right) \Big\Vert_2 \\
	       &= \mathcal{O}_{\rm p} \left(\sqrt{p} n^{-\alpha} \right) + \mathcal{O}_{\rm p} \left(\sqrt{p} n^{-\alpha} H^{-\nicefrac{1}{2}}\right) = \mathcal{O}_{\rm p} \left(\sqrt{p} n^{-\alpha} \right),
	   \end{aligned} 
	\end{equation*}
	%
	%
	by Assumption \ref{assum:C'}. Similarly, we have
	\begin{equation*}
	   \begin{aligned}
	    \Big| {Q}(\bm{\theta}, n) - Q(\bm{\theta}) \Big| &= \Big| \big\Vert {\bm{\pi}}(\bm{\theta}_0, n) - {\bm{\pi}}(\bm{\theta}_0) + {\bm{\pi}}(\bm{\theta}_0) - {\bm{\pi}}(\bm{\theta})  + {\bm{\pi}}(\bm{\theta}) - {\bm{\pi}}(\bm{\theta}, n) \big\Vert_2 \\
	       &- \big\Vert {\bm{\pi}}(\bm{\theta}_0) -  {\bm{\pi}}(\bm{\theta}) \big\Vert_2 \Big| \\
	       &\leq  \big\Vert {\bm{\pi}}(\bm{\theta}_0, n) - {\bm{\pi}}(\bm{\theta}_0) \big\Vert_2 + \big\Vert {\bm{\pi}}(\bm{\theta}) - {\bm{\pi}}(\bm{\theta}, n) \big\Vert_2 \\
	       &=   \big\Vert \mathbf{c}(n) + \mathbf{L}(n)\bto + \mathbf{r}(\bto,n) \big\Vert_2 + \big\Vert \mathbf{c}(n) + \mathbf{L}(n)\bm{\theta} + \mathbf{r}(\bm{\theta},n) \big\Vert_2 \\
	       &= \mathcal{O} \left(\sqrt{p} \max \left( c_n, p n^{-\beta}, n^{-\gamma} \right)\right),
	   \end{aligned} 
	\end{equation*}
	by Assumption \ref{assum:D:2}. Therefore, we obtain
    \begin{equation*}
        \sup_{\bm{\theta} \in \bm{\Theta}} \; \Big| \hat{Q}(\bm{\theta}, n) - Q(\bm{\theta}) \Big| = \mathcal{O}_{\rm p} \left(\sqrt{p} n^{-\alpha} \right) + \mathcal{O} \left(\sqrt{p} \max \left(c_n, p n^{-\beta}, n^{-\gamma} \right)\right).
    \end{equation*}
	By Assumptions \ref{assum:C'} and \ref{assum:D:2}, there exists a sample size $n^\ast \in\mathbb{N}^\ast$ such that for all $n \in \mathbb{N}^*$ satisfying $n \geq n^\ast$ we have
	\begin{equation*}
	    \Pr \left( \sup_{\bm{\theta} \in \bm{\Theta}} \; \Big| \hat{Q}(\bm{\theta}, n) - Q(\bm{\theta}) \Big| \geq \varepsilon \right) \leq   \delta.
	\end{equation*}
	Therefore, the four condition of Theorem 2.1 of \cite{newey1994large} are verified implying the result.
	\end{proof}
	
	Proposition \ref{THM:consistency} implies that when $p \to c <\infty$ as $n \to \infty$ we can simply write:
	\begin{equation*}
	   \lim_{n \to \infty} \Pr \left(|| \hbt - \bm{\theta}_0 ||_2 > \varepsilon \right) = 0.
	\end{equation*}
	However, in the case where $p \to \infty$ as $n \to \infty$, the statement is weaker in the sense that we cannot conclude that the limit exist but only that $ \Pr (|| \hbt - \bm{\theta}_0 ||_2 > \varepsilon )$ is arbitrarily small.
	\newpage
\setcounter{equation}{0}
\renewcommand{\theequation}{F.\arabic{equation}}
\section{Phase Transition Unbiasedness of the IB-Estimator}\label{app:unbias:IB}

In this appendix, we study the finite sample bias of the IB-estimator $\hbt$ and show that it is a PT-unbiased estimator. For this purpose, a slightly modified assumption framework for Theorem~\ref{THM:bias} is considered, where a modified version of Assumption \ref{assum:D'} is introduced, while Assumption \ref{assum:A} and \ref{assum:B'} remains unchanged. Assumption \ref{assum:C'} is not needed here since we are considering the bias of an estimator. 
For completness, we first restate Assumption \ref{assum:A} using our notation.
    \setcounter{Assumption}{0}
    \renewcommand{\theHAssumption}{otherAssumption\theAssumption}
    \renewcommand\theAssumption{\Alph{Assumption}}
    \begin{Assumption}
    \label{assum:A'}
     Let $\bm{\Theta}$ be a convex and compact subset of $\real^p$ such that $\btheta_0 \in \operatorname{Int}(\bm{\Theta})$ and
    \begin{equation*}
			\argzero_{\bm{\theta} \in
		      \real^p \setminus \bm{\Theta}} \;  \hat{\bm{\pi}}(\bm{\theta}_0, n, \bm{\omega}_0) - \frac{1}{H} \sum_{h = 1}^H  \hat{\bm{\pi}}(\bm{\theta}, n, \bm{\omega}_{h})
		      = \emptyset\,.
	\end{equation*}
    \end{Assumption}
    \vspace{0.25cm}
    Next, we introduce Assumption \ref{assum:D:3}.
    \setcounter{Assumption}{3}
    \renewcommand{\theHAssumption}{otherAssumption\theAssumption}
    \renewcommand\theAssumption{\Alph{Assumption}$_3$}
    \begin{Assumption}
    \label{assum:D:3}
    Consider the functions $\mathbf{d}(\bm{\theta}, n)$, $\mathbf{a}(\btheta)$, $\mathbf{b}(\bm{\theta}, n)$ and $\mathbf{r}(\bm{\theta}, n)$ defined in \eqref{def:fct:bias:d} and \eqref{def:fct:bias:b} and  let $\mathbf{R} (\bm{\theta}, n) \equiv \frac{\partial}{\partial \, \btheta^T}  \mathbf{r}(\bm{\theta}, n) \in \real^{p \times p}$ be the Jacobian matrix of $\mathbf{r} (\bm{\theta}, n)$ in $\btheta$.  These functions are such that for all $\btheta \in \bm{\Theta}$:
    \begin{enumerate}
    \item There exists a matrix $\mathbf{M} \in \real^{p \times p}$ and a vector $\mathbf{s} \in \real^{p}$ such that $\mathbf{a}(\btheta) = \mathbf{M}\btheta + \mathbf{s}$.
    \item There exists a sample size $n^* \in \mathbb{N}^*$ such that for all $n \in \mathbb{N}^*$ satisfying $n \geq n^*$ the matrix $(\mathbf{M} + \mathbf{L}(n)+\mathbf{I})^{-1}$ exists.
    \item There exists a real $\gamma > 0$ such that for all $j = 1, \dots, p$ we have
    \begin{equation*}
	  \begin{aligned}
			 \mathbf{r}_j\left(\bm{\theta}, n\right) = \mathcal{O}(n^{-\gamma}) \;\;\; \text{and} \;\;\;
			 \displaystyle\lim_{n\to\infty}\frac{p^2}{n^{\gamma}} = 0.%
		\end{aligned}
	\end{equation*}
   \item $\mathbf{R}(\bm{\theta}, n)$ exists and is continuous in $\btheta\in\bm{\Theta}$ 
   for all $n \in \Ns$ with $n \geq n^\ast$.
    \end{enumerate}
    \end{Assumption}
    \vspace{0.25cm}
    
     Compared to Assumption \ref{assum:D'}, Assumption \ref{assum:D:3} is weaker. First, the matrix $\mathbf{M}$ may be such that $\lVert\mathbf{M}\rVert_F \geq 1$, thereby relaxing the contraction mapping hypothesis on $\mathbf{a}(\btheta)$ given in Assumptions \ref{assum:D},
      \ref{assum:D'} and \ref{assum:D:1}. Second, conditions 3 and 5 of Assumption \ref{assum:D'} are removed.

     Clearly, neither Assumption \ref{assum:D:2} 
     nor Assumption \ref{assum:D:3} imply each other. The first part of Assumption \ref{assum:D:3} imposes a restrictive form for the bias function $\mathbf{a}(\btheta)$. If $\hat{\bm{\pi}}(\bm{\theta}_0, n, \bm{\omega})$ is a consistent estimator of $\bm{\theta}_0$, then this restriction is automatically satisfied. 
    Under these new conditions, the PT-unbiasedness of $\hbt$ is established in Proposition \ref{app:prop:bias} below. 
    
    \begin{Proposition}  
		\label{app:prop:bias}
		Under Assumptions \ref{assum:A'}, \ref{assum:B'} and \ref{assum:D:3}, $\hbt$ is a PT-unbiased estimator, i.e. there exists a sample size $n^* \in \mathbb{N}^*$ such that for all $n \in \mathbb{N}^*$ satisfying $n \geq n^*$ and for all $H \in\mathbb{N}^\ast$, we have $\Big\lVert	\mathbb{E} \left[\hbt \right] - \bto  \Big\rVert_2 = 0$. 
	\end{Proposition}
	\vspace{0.25cm}
	
	  Before giving the proof we introduce the following useful notation. Let $\mathbf{f}\,  : \bm{\Theta} \to \real^p$ and $\mathbf{F}(\btheta)~\equiv~ \frac{\partial}{\partial \, \btheta^T} \mathbf{f}\left(\bm{\theta}\right) \in \real^{p \times p}$ be its Jacobian matrix. Then, we define $\mathbf{F}\big(\btheta^{(\mathbf{f})}\big)$ such that when using the multivariate mean value theorem between $\hbt$ and ${\btheta}_{0}$ we obtain
\begin{equation}
    \mathbf{f} \left(\hbt\right) =  \mathbf{f} \left({\btheta}_{0}\right) + \mathbf{F}\left(\btheta^{(\mathbf{f})}\right) \left( \hbt - {\btheta}_{0} \right).
    \label{def:MVT:multi_vari}
\end{equation}

Therefore, $\btheta^{(\mathbf{f})}$ corresponds to a set of $p$ vectors lying in the segment $(1 - \lambda) \hbt + \lambda {\btheta}_{0}$ for $\lambda \in [0,1]$ (with respect to the function $\mathbf{f}$). Keeping the latter notation in mind (for a generic function $\mathbf{f}$), we provide the proof of Proposition \ref{app:prop:bias}.
	
	\begin{proof}
        From (\ref{app:eq:indirectInf:hTimesN}) we have that $\hbt$ is such that for $(n, H) \in \mathbb{N}^\ast \times \mathbb{N}^\ast$ 
	\begin{equation*}
	    \hat{\bm{\pi}}(\bm{\theta}_0, n, \bm{\omega}_0) = \frac{1}{H} \sum_{h = 1}^H  \hat{\bm{\pi}}(\hbt, n, \bm{\omega}_{h}).
	\end{equation*}
	Using Assumptions \ref{assum:B'} and \ref{assum:D:3}, we expand each side of the above equation as follows
	\begin{equation}
	    \begin{aligned}
	     \hat{\bm{\pi}}(\bm{\theta}_0, n, \bm{\omega}_0) &= \bm{\theta}_0 + \mathbf{a}(\bm{\theta}_0) + \mathbf{c}(n) + \mathbf{L}(n) \bm{\theta}_0 + \mathbf{r} (\bm{\theta}_0, n) + \mathbf{v} \left(\bm{\theta}_0, n, \bm{\omega}_0\right)\\
	     \frac{1}{H} \sum_{h = 1}^H  \hat{\bm{\pi}}(\hbt, n, \bm{\omega}_{h}) &= \hbt + \mathbf{a}(\hbt) +\mathbf{c}(n) + \mathbf{L}(n) \hbt + \mathbf{r} (\hbt, n) \\
	     &\quad +  \frac{1}{H} \sum_{h = 1}^H  \mathbf{v}(\hbt, n, \bm{\omega}_{h}).
	    \end{aligned}
	    \label{eq:proof:bias:inter:1}
	\end{equation}
	Therefore, we obtain
	\begin{equation*}
	    \begin{aligned}
	    &\mathbb{E}\left[\frac{1}{H} \sum_{h = 1}^H  \hat{\bm{\pi}}(\hbt, n, \bm{\omega}_{h}) - \hat{\bm{\pi}}(\bm{\theta}_0, n, \bm{\omega}_0) \right]\\
	    &= \left(\mathbf{I} + \mathbf{L}(n)\right)\mathbb{E}\left[\hbt - \bm{\theta}_0\right] + \mathbb{E}\left[\mathbf{a}(\hbt) - \mathbf{a}(\bm{\theta}_0)\right] + \mathbb{E} \left[\mathbf{r} (\hbt, n) - \mathbf{r} (\bm{\theta}_0, n)\right]\\
	    &= \left(\mathbf{I} + \mathbf{L}(n) + \mathbf{M}\right) \mathbb{E}\left[\hbt - \bm{\theta}_0\right] +\mathbb{E} \left[\mathbf{r} (\hbt, n) - \mathbf{r} (\bm{\theta}_0, n)\right]\\
	    &= \mathbf{0}.
	    \end{aligned} 
	\end{equation*}
	By hypothesis, $n$ is such that the inverse of $\mathbf{B} \equiv \mathbf{I} + \mathbf{L}(n) + \mathbf{M}$ exists and we obtain
	\begin{equation}
	    \mathbb{E}\left[\hbt - \bm{\theta}_0\right] = -\mathbf{B}^{-1} \mathbb{E} \left[\mathbf{r} (\hbt, n) - \mathbf{r} (\bm{\theta}_0, n)\right]. 
	    \label{proof:bias:eq:inter2}
	\end{equation}
    By Assumptions \ref{assum:A'} and \ref{assum:D:3}, $\mathbf{r} (\hbt, n)$ is a bounded random variable on a compact set. Moreover, since $\mathbf{r}(\bm{\theta},n) = \mathcal{O}(n^{-\gamma})$ elementwise by Assumption \ref{assum:D:3}, we have 
    \begin{equation}
        \mathbb{E} \left[\mathbf{r} (\hbt, n) - \mathbf{r} (\bto,n) \right] = \mathcal{O}(n^{-\gamma})
    \end{equation} elementwise. Consequently, we deduce from \eqref{proof:bias:eq:inter2} that 
    \begin{equation}
        \mathbb{E} \left[ \hbt - \bto \right] = \mathcal{O}(pn^{-\gamma})
        \label{proof:bias:eq:inter3}
    \end{equation}
    elementwise. 
    
    The idea now is to re-evaluate $\mathbb{E} \left[\mathbf{r} (\hbt, n) - \mathbf{r} (\bm{\theta}_0, n)\right]$ using the mean value theorem. This will allow us to make an induction that will show that for all $\delta\in\mathbb{N}$ 
	\begin{equation*}
	     \Big\lVert \mathbb{E}\left[ \hbt - \bto \right] \Big\rVert_2 = \mathcal{O}\left((p^2n^{-\gamma})^\delta\right).
	\end{equation*}
	We will then use Lemma \ref{lemma:newO} to conclude the proof.
	
    Applying the mean value theorem for vector-valued functions to $\mathbf{r} (\hbt, n) - \mathbf{r} (\bm{\theta}_0, n)$ we have
	\begin{equation*}
	    \mathbf{r} (\hbt, n)- \mathbf{r} (\bm{\theta}_0, n) = \mathbf{R}\left(\bt^{(\mathbf{r})},n\right) \left(\hbt - \bm{\theta}_0 \right).
	\end{equation*}
	where the matrix $\mathbf{R}\left(\bt^{(\mathbf{r})},n\right)$ is defined in \eqref{def:MVT:multi_vari}. 
	By Assumptions \ref{assum:A'} and \ref{assum:D:3}, $\mathbf{R}\left(\bt^{(\mathbf{r})},n\right)$ is also a bounded random variable. Moreover, we have $\mathbb{E}\left[\mathbf{R}\left(\bt^{(\mathbf{r})},n\right)\right] = \mathcal{O}(n^{-\gamma})$ elementwise since by Assumption \ref{assum:D:3} $\mathbf{r}(\bm{\theta},n) = \mathcal{O}(n^{-\gamma})$ elementwise, and 
	\begin{equation*}
	    \mathbf{r} (\bm{\theta}, n) = \mathbf{r} (\bm{\theta}_0, n) + \mathbf{R}\left(\bt^{(\mathbf{r})},n\right) \left(\bm{\theta} - \bm{\theta}_0 \right).
	\end{equation*}
	for all $\bm{\theta} \in \bm{\Theta}$.
	For simplicity, we denote 
	\begin{equation*}
	    \mathbf{R} \equiv \mathbf{R}\left(\bt^{(\mathbf{r})},n\right),\ \ \
	    \mathbf{\Delta}^{(\mathbf{r})} \equiv  \mathbf{r} (\hbt, n) - \mathbf{r} (\bm{\theta}_0, n)
	    \ \ \ \text{and} \ \ \ \mathbf{\Delta} \equiv \hbt - \bm{\theta}_0,
	\end{equation*} 
	which implies that $\mathbf{\Delta}^{(\mathbf{r})}=\mathbf{R}\mathbf{\Delta}$. Moreover, a consequence of \eqref{proof:bias:eq:inter3} is that 
		\begin{equation*}
	    \big|\mathbb{E}\left[\mathbf{\Delta}_j\right]\big|=\mathcal{O}(pn^{-\gamma}),
	\end{equation*}
	for any $j=1,\dots,p$ and hence 
	\begin{equation*}
	    \Vert \mathbb{E}\left[\mathbf{\Delta}\right] \Vert_2 = \mathcal{O}\left(p^2n^{-\gamma}\right).
	\end{equation*}
	Now, for all $l=1,\dots,p$ we have 
	\begin{equation}
	    \mathbf{\Delta}^{(\mathbf{r})}_l=\sum^{p}_{m=1}\mathbf{R}_{l,m}\mathbf{\Delta}_m \leq p \max_{m}\mathbf{R}_{l,m}\mathbf{\Delta}_m.
	    \label{proof:bias:eq:inter4}
	\end{equation}

  Using Cauchy-Schwarz inequality, we have 
   \begin{equation*}
       \mathbb{E} \left[ | \mathbf{R}_{l,m}\mathbf{\Delta}_m |\right] \leq  \mathbb{E} \left[ \mathbf{R}_{l,m}^2 \right]^{\nicefrac{1}{2}} \mathbb{E} \left[ \mathbf{\Delta}_m^2 \right]^{\nicefrac{1}{2}}.
    \end{equation*}
	Since $\mathbf{R}_{l,m}$ and $\mathbf{\Delta}_m$ are bounded random variables on a compact set, $\mathbb{E}\left[\mathbf{R}_{l,m}\right] = \mathcal{O}\left(n^{-\gamma}\right)$ and $\mathbb{E}\left[\mathbf{\Delta}_m\right] = \mathcal{O}\left(pn^{-\gamma}\right)$ we have 
	\begin{equation}
	\label{eq:unif:O:bound}
   \mathbb{E} \left[ \mathbf{R}_{l,m}^2 \right]^{\nicefrac{1}{2}} \mathbb{E} \left[ \mathbf{\Delta}_m^2 \right]^{\nicefrac{1}{2}} = \mathcal{O}\left(n^{-\gamma}\right)\mathcal{O}\left(p n^{-\gamma}\right)=\mathcal{O}\left(p n^{-2\gamma}\right).
	\end{equation}
	Therefore, we obtain 
	\begin{equation*}
    \mathbb{E} \left[ | \mathbf{R}_{l,m}\mathbf{\Delta}_m |\right] \leq  \mathbb{E} \left[ \mathbf{R}_{l,m}^2 \right]^{\nicefrac{1}{2}} \mathbb{E} \left[ \mathbf{\Delta}_m^2 \right]^{\nicefrac{1}{2}} = \mathcal{O}\left(p n^{-2\gamma}\right),
    \end{equation*}
	and hence, using \eqref{proof:bias:eq:inter4} we deduce that
	\begin{equation*}
        \Big| \mathbb{E} \left[  \mathbf{\Delta}^{(\mathbf{r})}_l \right] \Big| \leq p \mathbb{E} \left[ \max_{m} | \mathbf{R}_{l,m}\mathbf{\Delta}_m | \right] = \mathcal{O}\left((pn^{-\gamma})^2\right). 
    \end{equation*}
	Considering this, and since $\mathbb{E}\left[\mathbf{\Delta}\right] = - \mathbf{B}^{-1}\mathbb{E}\left[\mathbf{\Delta}^{(\mathbf{r})}\right]$, we have, for all $j=1,\dots,p$
	\begin{equation*}
	   \big| \mathbb{E} \left[  \mathbf{\Delta}_j\right] \big| =  \mathcal{O}\left(p^3n^{-2\gamma}\right),
	\end{equation*}
	and consequently,
	\begin{equation*}
	   \lVert \mathbb{E} \left[  \mathbf{\Delta}\right] \rVert_2 =  \mathcal{O}\left((p^2n^{-\gamma})^2\right).
	\end{equation*}
	Since $\mathbb{E}\left[\mathbf{\Delta}\right] = -\mathbf{B}^{-1}\mathbb{E}\left[\mathbf{\Delta}^{(\mathbf{r})}\right]$ and $\mathbb{E}\left[\mathbf{\Delta}^{(\mathbf{r})}\right] = \mathbb{E}\left[\mathbf{R}\mathbf{\Delta}\right]$, one can repeat the same computations and deduce by induction that for all $\delta \in \mathbb{N}$
	\begin{equation*}
	    \Vert \mathbb{E}\left[\mathbf{\Delta}\right] \Vert_2 = \mathcal{O}\left((p^2n^{-\gamma})^\delta\right).
	\end{equation*} 
	Using equation \eqref{eq:unif:O:bound} we see that     
	\begin{equation*}
	    \Vert \mathbb{E}\left[\mathbf{\Delta}\right] \Vert_2 = \mathcal{O}_{\delta\in\mathbb{N}}\left((p^2n^{-\gamma})^\delta\right),
	\end{equation*}
	which ends the proof by Lemma \ref{lemma:newO}.	
	\end{proof}
	
	Proposition \ref{app:prop:bias} provides a strategy for proving the consistency of the IB-estimator $\hbt$. However, this requires slight modifications of the assumptions considered in Proposition \ref{app:prop:bias}. More in detail, we impose stronger requirements on $\alpha$ and $\gamma$ in  Assumptions \ref{assum:C'} and \ref{assum:D:3}. 
	\setcounter{Assumption}{2}
	\renewcommand{\theHAssumption}{otherAssumption\theAssumption}
	\renewcommand\theAssumption{\Alph{Assumption}$_1$}
	\begin{Assumption}
	\label{assum:C:1}
	Preserving the same definitions and requirements of Assumption \ref{assum:C'}, we additionally require that the constant $\alpha$ is such that
	\begin{equation*}
	   \lim_{n \to \infty} \frac{p^{\nicefrac{5}{2}}}{n^\alpha}=0.
	\end{equation*}
	\end{Assumption}
	\vspace{0.25cm}
	
	\setcounter{Assumption}{3}
	\renewcommand{\theHAssumption}{otherAssumption\theAssumption}
	\renewcommand\theAssumption{\Alph{Assumption}$_4$}
	\begin{Assumption}
	\label{assum:D:4}
	
	Preserving the same definitions and requirements of Assumption \ref{assum:D:3}, we additionally require that the constant $\gamma$ is such that
	\begin{equation*}
	   \lim_{n \to \infty} \frac{p^{\nicefrac{5}{2}}}{n^\gamma}=0.
	\end{equation*}
	\end{Assumption}
	\vspace{0.25cm}
	
	Clearly, Assumption \ref{assum:C'} is weaker than Assumption \ref{assum:C:1}. Neither Assumption \ref{assum:D:4} nor Assumption \ref{assum:D} imply each other. Indeed, in the situations where $\alpha = \nicefrac{1}{2}$ and $\gamma = 2$, the estimator $\hbt$ is consistent for $\bm{\theta}_0$ under the conditions of Proposition \ref{app:prop:bias} and if
   	\begin{equation*}
	    \lim_{n \to \infty} \; \frac{p^5}{n} = 0,
	\end{equation*}
	which is clearly a stronger requirement than (\ref{eq:n:p:app}) based on Proposition \ref{THM:consistency} (using Assumption \ref{assum:D:2} which is weaker than Assumption \ref{assum:D}). However, some requirements of Assumption \ref{assum:D} are stronger than the ones in Assumption \ref{assum:D:4}. An important difference is that the function $\mathbf{a}(\bt)$ doesn't need to be a contraction map in Assumption \ref{assum:D:4}. In addition, most of the requirements on the matrix $\mathbf{L}(n)$ and the vector $\mathbf{c}(n)$ are removed. Using these new assumptions, we present the following corollary to Proposition \ref{app:prop:bias}. 
	
	\begin{Corollary}
	\label{coro:consist}
	Under Assumptions \ref{assum:A'}, \ref{assum:B'}, \ref{assum:C:1}, and \ref{assum:D:4}, $\hbt$ is a consistent estimator of $\bm{\theta}_0$ for all $H \in \mathbb{N}^\ast$ in the sense of Proposition \ref{THM:consistency}.
	\end{Corollary}
	\vspace{0.25cm}
    \begin{proof}
     Fix $\varepsilon > 0$ and $\delta > 0$. We need to show that there exists a sample size $n^\ast \in\mathbb{N}^\ast$ such that for all $n \in \mathbb{N}^*$ satisfying $n \geq n^\ast$
	\begin{equation*}
	    \Pr \left(|| \hbt - \bm{\theta}_0 ||_2 \geq \varepsilon \right) \leq   \delta.
	\end{equation*}
	From Chebyshev's inequality we have that
	\begin{equation*}
	    \Pr \left(|| \hbt - \bm{\theta}_0 ||_2  \geq \varepsilon \right) \leq \frac{\mathbb{E}\left[|| \hbt - \bm{\theta}_0 ||_2^2\right]}{\varepsilon^2}.
	\end{equation*}
	Therefore, we only need to show that there exists a sample size $n^\ast \in\mathbb{N}^\ast$ such that for all $n \geq n^\ast$
	\begin{equation*}
	    \mathbb{E}\left[|| \hbt - \bm{\theta}_0 ||_2^2\right] \leq \varepsilon^2\delta .
	\end{equation*}

	Using the same notation as in the proof of Proposition \ref{app:prop:bias}, we have from (\ref{eq:proof:bias:inter:1}) and for $n$ sufficiently large that
	\begin{equation}
	\begin{aligned}
	   \mathbf{\Delta} & \equiv \hbt - \bm{\theta}_0 \\
	    & = \mathbf{B}^{-1} \left(\mathbf{r} (\bm{\theta}_0, n) - \mathbf{r} (\hbt, n) + \mathbf{v} \left(\bm{\theta}_0, n, \bm{\omega}_0\right) - \frac{1}{H} \sum_{h = 1}^H  \mathbf{v}(\hbt, n, \bm{\omega}_{h})\right)\\
	    & =  \mathbf{B}^{-1} \left( \mathbf{R}\mathbf{\Delta} + \mathbf{v} \right), 
	    \end{aligned}
	    \label{eq:comment:consist}
	\end{equation}
	where $\mathbf{v}$ is zero mean random vector of order $\mathcal{O}_{\rm p}\left(n^{-\alpha}\right)$ elementewise by Assumption \ref{assum:C'}. Thus, by Assumption \ref{assum:D:4} we have that for all $j = 1,\dots,p$ 
	\begin{equation*}
	   \begin{aligned}
	        \mathbf{\Delta}_j &= \sum^p_{l=1} \left( \mathbf{B}^{-1} \right)_{j,l}\left(\mathbf{R}\mathbf{\Delta} + \mathbf{v}\right)_l = 
	        \sum^p_{l=1} \sum^p_{m=1} \left( \mathbf{B}^{-1} \right)_{j,l}\left(\mathbf{R}_{l,m}\mathbf{\Delta}_m + \mathbf{v}_m\right) \\
	         &= \sum^p_{l=1} \sum^p_{m=1} \mathcal{O_{\rm p}}(n^{-\gamma}) + \mathcal{O_{\rm p}}(n^{-\alpha}) =
	         \sum^p_{l=1} \sum^p_{m=1} \mathcal{O_{\rm p}}\left(\max(n^{-\alpha} , n^{-\gamma})\right) 
	         \\
	         &= \mathcal{O_{\rm p}}\left(p^2\max(n^{-\alpha} , n^{-\gamma})\right).
	   \end{aligned}
	\end{equation*}
	Therefore we have,
	%
	\begin{equation*}
	  \left\lVert \mathbf{\Delta} \right\rVert_2^2 = \sum_{k=1}^p \mathbf{\Delta}_k^2 = p\mathcal{O}_{\rm  p}\left(p^4\max(n^{-2\alpha}  , n^{-2\gamma})\right) = \mathcal{O}_{\rm  p}\left(p^5 n^{-2\min(\alpha, \gamma)}\right),
	\end{equation*}
	and thus since $\bm{\Theta}$ is compact by Assumption \ref{assum:A'}, we have
	\begin{equation*}
	\begin{aligned}
	    \mathbb{E}\left[\lVert \mathbf{\Delta} \rVert_2^2\right] = \mathcal{O}\left(p^5 n^{-2\min(\alpha, \gamma)}\right).
	    \end{aligned}
	\end{equation*}
	Using Assumptions \ref{assum:C:1} and \ref{assum:D:4} the last equality implies that there exists a sample size $n^\ast \in\mathbb{N}^\ast$ such that for all $n \in \mathbb{N}^*$ satisfying $n \geq n^\ast$, we obtain 
	\begin{equation*}
	\begin{aligned}
	    \mathbb{E}\left[\lVert \mathbf{\Delta} \rVert_2^2\right] \leq \varepsilon^2\delta.
	    \end{aligned}
	\end{equation*}

    \end{proof}
	
\newpage
\setcounter{equation}{0}
\renewcommand{\theequation}{G.\arabic{equation}}
\section{Asymptotic Distribution of the IB-Estimator}\label{app:asymp:norm:IB}

In order to study the asymptotic normality of the IB-estimator $\hbt$ we  modify some of the assumptions in order to incorporate additional requirements on the auxiliary estimator. 

Before introducing the next assumption we define various quantities. First, we let
 \begin{equation*}
    \bm{\Sigma}(\btheta,n) \equiv \var \left( \sqrt{n} \mathbf{v}(\bm{\theta}, n, \bm{\omega}) \right),
\end{equation*}
where $\bm{\Sigma}(\btheta,n)$ is nonsingular.
Then, using this definition, we introduce the following quantity
\begin{equation*}
    Y(\btheta, n, \bm{\omega}, \mathbf{u}) \equiv \sqrt{n} \mathbf{u}^T \bm{\Sigma}(\btheta,n)^{-\nicefrac{1}{2}} \mathbf{v}(\bm{\theta}, n, \bm{\omega}),
\end{equation*}
where $\mathbf{u} \in \real^p$ is such that $||\mathbf{u} ||_2 = 1$. Clearly, from the definition of $\mathbf{v}(\bm{\theta}, n, \bm{\omega})$ we have that $\mathbb{E}[Y(\btheta, n, \bm{\omega}, \mathbf{u})] = 0$. Therefore, without loss of generality we can always decompose $Y(\btheta, n, \bm{\omega}, \mathbf{u})$ as follows:
\begin{equation}
    Y(\btheta, n, \bm{\omega}, \mathbf{u}) =  Z(\btheta, n, \bm{\omega}, \mathbf{u}) + \delta_n W(\btheta, n, \bm{\omega}, \mathbf{u}),
    \label{eq:gauss:approx}
\end{equation}
where $W(\btheta, n, \bm{\omega}, \mathbf{u})$ and $Z(\btheta, n, \bm{\omega}, \mathbf{u})$ are zero mean random variables and $\big(\delta_n\big) \in \mathfrak{D}$ where
\begin{equation*}
	    \mathfrak{D} \equiv \left\{ \left( \delta_n \right)_{n \in \mathbb{N}^\ast} \, \big| \; \delta_n \in \real \setminus (-\infty, 0) \; \forall n, \; \delta_l \geq \delta_m \;\; \text{if} \;\; l < m \;\; \text{and} \;\; \lim_{n \to \infty} \delta_n = 0\right\}.
\end{equation*}
In Assumption \ref{assum:C:ast} below, we restrict the behaviour of $\mathbf{v}(\bm{\theta}, n, \bm{\omega})$ and in particular we require that it satisfies a specific Gaussian approximation based on the decomposition considered in (\ref{eq:gauss:approx}).

\setcounter{Assumption}{2}
\renewcommand{\theHAssumption}{otherAssumption\theAssumption}
\renewcommand\theAssumption{\Alph{Assumption}$^\ast$}

\begin{Assumption}
\label{assum:C:ast}
Preserving the requirements of Assumption \ref{assum:C'} we add the following conditions on $\mathbf{v}(\bm{\theta}, n, \bm{\omega})$. There exists a sample size $n^\ast \in \mathbb{N}^*$ such that:
\begin{enumerate}
\item For all $\bm{\theta}\in \bm{\Theta}$ and all $n \in \mathbb{N}^*$ with $n \geq n^\ast$, the matrix $\bm{\Sigma}(\btheta,n)$ exists, is nonsingular and is such that $\bm{\Sigma}_{j,l}(\btheta,n) = \mathcal{O}(1)$ for all $j,l=1,\dots,p$.
\item For all $\bm{\omega}\in\bm{\Omega}$ and all $n \in \mathbb{N}^*$ with $n \geq n^\ast$, the Jacobian matrix $\mathbf{V}(\bm{\theta}, n , \bm{\omega}) \equiv \frac{\partial}{\partial \, \btheta^T} \mathbf{v}\left(\bm{\theta}, n, \bm{\omega}\right)$ exists and is continuous in $\btheta \in \bm{\Theta}$. 
\item  Considering the decomposition in \eqref{eq:gauss:approx}, then for all $\mathbf{u} \in \real^p$ such that $||\mathbf{u} ||_2 = 1$, there exist sequences $\big(\delta_n\big) \in \mathfrak{D}$, as well as a random variable $W(\btheta, n, \bm{\omega}, \mathbf{u}) = \mathcal{O}_{\rm p}(1)$, such that
$Z(\btheta, n, \bm{\omega}, \mathbf{u})$ is a standard normal random variable.
\end{enumerate}
\end{Assumption}

The first requirement of Assumption \ref{assum:C:ast} is quite mild and commonly assumed as it simply requires that $\bm{\Sigma}(\btheta, n)$ is a suitable covariance matrix. In addition, it implies that $\mathbf{v}(\bm{\theta}, n, \bm{\omega}) = \mathcal{O}_{\rm p} \left(n^{-\nicefrac{1}{2}}\right)$, elementwise. Similarly to the discussion following Assumption \ref{assum:D'}, our second requirement  ensures that $\mathbf{V}_{j,l}(\hbt, n , \bm{\omega})$ is a bounded random variable for all $j,l=1,\dots,p$ and all $(n, \bm{\omega}) \in \mathbb{N}^* \times \bm{\Omega}$. 
The third condition of Assumption \ref{assum:C:ast} describes how ``close'' the distribution of $\mathbf{v}(\bm{\theta}, n, \bm{\omega})$ is to a multivariate normal distribution. Such an assumption is quite strong and may not always be satisfied. In the case where $p \to c < \infty$ as $n \to \infty$, our requirement on the distribution of $\mathbf{v}(\bm{\theta}, n, \bm{\omega})$ can simply be expressed as
\begin{equation*}
    \sqrt{n} \mathbf{v}(\bm{\theta}, n, \bm{\omega}) \xrightarrow{d} \mathcal{N}(\mathbf{0}, \bm{\Sigma}(\btheta)),
\end{equation*}
where $\bm{\Sigma}(\btheta) \equiv \displaystyle \lim_{n \to \infty} \bm{\Sigma}(\btheta, n)$. Finally, Assumption \ref{assum:C:ast} clearly implies Assumption \ref{assum:C'} but not necessarily Assumption~\ref{assum:C:1}.

In Assumption \ref{assum:D:5} below, we impose an additional requirement on the bias of the auxiliary estimator.

\setcounter{Assumption}{3}
\renewcommand{\theHAssumption}{otherAssumption\theAssumption}
\renewcommand\theAssumption{\Alph{Assumption}$_5$}
\begin{Assumption}
\label{assum:D:5}
There exists a sample size $n^\ast \in \mathbb{N}^\ast$ such that for all $n \in \mathbb{N}^*$ satisfying $n \geq n^\ast$  the Jacobian matrices $\mathbf{A}(\bm{\theta}) \equiv \frac{\partial}{\partial \, \btheta^T} \mathbf{a}\left(\bm{\theta}\right)$ and $\mathbf{R}(\bm{\theta}, n) \equiv \frac{\partial}{\partial \, \btheta^T} \mathbf{r}\left(\bm{\theta}, n\right)$ exist and are continuous in $\btheta\in\bm{\Theta}$. Moreover, we require that for all $\mathbf{u}\in\real^p$ with $||\mathbf{u}||_2 = 1$, there exists a sequence $(\delta_n) \in\mathfrak{D}$ such that 
\begin{equation}
    \scalebox{0.98}{$
    \sqrt{n}\left(1+\frac{1}{H}\right)^{-\nicefrac{1}{2}}\mathbf{u}^T\bm{\Sigma}(\bto,n)^{-\nicefrac{1}{2}}\Big[\mathbf{B}(\bto,n)-\mathbf{B}\left(\btheta^{(\mathbf{a}, \mathbf{r})},n\right)\Big]\left(\hbt-\bto\right) = \mathcal{O}_{\rm p}(\delta_n)
    $},
    \label{eq:assum:d:5}
\end{equation}
where 
\begin{equation*}
    \mathbf{B}(\bt,n) \equiv \mathbf{I} + \mathbf{A}(\bt) + \mathbf{L}(n) + \mathbf{R}(\bt,n),
\end{equation*}
\begin{equation*}
    \mathbf{B}\big(\btheta^{(\mathbf{a}, \mathbf{r})},n\big) \equiv \mathbf{I} + \mathbf{A}(\btheta^{(\mathbf{a})}) + \mathbf{L}(n) + \mathbf{R}\left(\btheta^{(\mathbf{r})},n\right).
\end{equation*}
while $\mathbf{A}(\btheta^{(\mathbf{a})})$ and $\mathbf{R}(\btheta^{(\mathbf{r})})$ are defined using notation in \eqref{def:MVT:multi_vari}.
\end{Assumption}
\vspace{0.25cm}

Assumption \ref{assum:D:5} allows us to quantify how ``far'' the matrices $\mathbf{B}(\bto,n)$ and $\mathbf{B}\big(\hbt, n\big)$ are from each other. In the case where $p \to c < \infty$ as $n \to \infty$, the consistency of $\hbt$ and the continuity of $\mathbf{B}(\bt) \equiv \displaystyle\lim_{n \to \infty} \mathbf{B}(\bt,n)$ in $\bm{\theta}_0$ would be sufficient so that Assumption \ref{assum:D:5} would not be needed to establish the asymptotic distribution of $\hbt$. However, when $p \to \infty$ as $n \to \infty$ the requirement in (\ref{eq:assum:d:5}) can be strong and difficult to verify for a specific model and auxiliary estimator. Having stated this, the following proposition defines the distribution of the IB-estimator.

\begin{Proposition}
\label{Thm:Gauss:approx}
Under Assumptions \ref{assum:A'}, \ref{assum:B'}, \ref{assum:C:ast} and \ref{assum:D:5}, for all $\mathbf{u} \in \real^p$ such that $||\mathbf{u}||_2 = 1$, there exist a sample size $n^\ast \in \mathbb{N}^\ast$ and a sequence $\left( \delta_n\right)\in\mathfrak{D}$ such that for all $n \in \mathbb{N}^\ast$ satisfying $n \geq n^\ast$ we have
\begin{equation*}
    \scalebox{0.99}{$
    \sqrt{n}\left(1+\frac{1}{H}\right)^{-\nicefrac{1}{2}}\mathbf{u}^T\bm{\Sigma}(\bto,n)^{-\nicefrac{1}{2}}\mathbf{B}(\bto, n)\left(\hbt-\bto\right) \overset{d}{=} Z + \delta_n\mathcal{O}_{\rm p}\left( \max\left(1,\frac{p^2}{\sqrt{H}}\right)\right)$},
\end{equation*}
where $Z \sim \mathcal{N}(0,1)$.
\end{Proposition}
\vspace{0.25cm}

\begin{proof}
Without loss of generality we can assume that the sample sizes $n^* \in \mathbb{N}^*$ of Assumption \ref{assum:C:ast} and \ref{assum:D:5} are the same. Let $n \in \mathbb{N}^*$ be such that $n \geq n^*$.
By the definition of the IB-estimator 
we have 
\begin{equation}\label{def:JINI}
     \hp{\bto}{\bwo} = \hpH{\hbt}.
\end{equation}
Since the auxiliary estimator may be expressed, for any $(\btheta, n, \bm{\omega}) \in \bm{\Theta} \times \mathbb{N}^* \times \bm{\Omega}$, as
\begin{equation}\label{eq:hat:pi}
    \hp{\bt}{\bw} = \bt + \a{\bt} + \cn + \mathbf{L}_n\bt + \r{\bt} + \v{\bt}{\bw}.
\end{equation}
we deduce from \eqref{def:JINI}, by using the mean value theorem on $\mathbf{a}(\btheta)$ and $\mathbf{r}(\btheta,n)$, that 
\begin{equation}
    \mathbf{B}\left(\btheta^{(\mathbf{a},\mathbf{r})},n\right)\bD= \v{\bto}{\bwo} - \vH{\hbt},
     \label{proof:prop4:B_THETA}
\end{equation}
where $\mathbf{B}\left(\btheta^{(\mathbf{a},\mathbf{r})},n\right) \equiv \mathbf{I} + \mathbf{A}(\btheta^{(\mathbf{a})}) + \mathbf{L}(n) + \mathbf{R}\left(\btheta^{(\mathbf{r})},n\right)$ and  $\bD \equiv \hbt - \bto$. Applying the mean value theorem on $\v{\hbt}{\bw_{h}}$ and using the notation defined in \eqref{def:MVT:multi_vari}, we obtain for all $h=1,\dots,H$
\begin{equation*}
    \v{\hbt}{\bw_{h}} = \v{\bto}{\bw_{h}} + \mathbf{V}\left(\btheta^{(\mathbf{v})},n,\bw_{h}\right)\bD.
\end{equation*}
Hence, let $\mathbf{u} \in \real^p$ be such that $|| \mathbf{u} ||_2 = 1$, then by Assumption \ref{assum:C:ast} we have, 
\begin{equation*}
    \sqrt{n} \mathbf{u}^T \bm{\Sigma}(\bto,n)^{-\nicefrac{1}{2}}\v{\bto}{\bwo}   = Z(\bto, n, \bwo, \mathbf{u}) + \delta_n^{(0)} W(\bto, n, \bwo, \mathbf{u}) 
\end{equation*}
and for all $h=1,\dots,H$
\begin{equation*}
    \sqrt{n} \mathbf{u}^T \bm{\Sigma}(\bto,n)^{-\nicefrac{1}{2}}\v{\bto}{\bw_{h}}   = Z(\bto, n, \bw_{h}, \mathbf{u}) + \delta_n^{(h)} W(\bto, n, \bw_{h}, \mathbf{u}) 
\end{equation*}
where $(\delta_n^{(0)}),(\delta_n^{(1)}),\dots, (\delta_n^{(H)}) \in \mathfrak{D}$ and 
\begin{equation*}
\begin{aligned}
        Z(\bto, n, \bwo, \mathbf{u})\sim \mathcal{N}(0,1), \;\;\; \;\;\;  Z(\bto, n, \bw_{h}, \mathbf{u}) \sim \mathcal{N}(0,1), \\
        W(\bto, n, \bwo, \mathbf{u}) = \mathcal{O}_{\rm p}(1), \;\;\; \;\;\;  W(\bto, n, \bw_{h}, \mathbf{u}) = \mathcal{O}_{\rm p}(1).
\end{aligned}
\end{equation*}
Without loss of generality, we can suppose that $(\delta_n) \equiv (\delta_n^{(0)})=(\delta_n^{(1)})=\dots=(\delta_n^{(H)})$ as one may simply modify $W(\bto, n, \bwo, \mathbf{u})$ and $W(\bto, n, \bw_{h}, \mathbf{u})$ for all $h=1,\dots,H$. For simplicity, we write $\mathfrak{X}_{\bto,n,\mathbf{u}} \equiv \sqrt{n} \mathbf{u}^T \bm{\Sigma}(\bto,n)^{-\nicefrac{1}{2}}$, $Z_0 \equiv Z(\bto, n, \bwo, \mathbf{u})$ and $Z_h \equiv Z(\bto, n, \bw_{h}, \mathbf{u})$.
Multiplying the right hand side of \eqref{proof:prop4:B_THETA} by $\mathfrak{X}_{\bto,n,\mathbf{u}}$, we have 
\begin{equation}
    \begin{aligned}
        &\mathfrak{X}_{\bto,n,\mathbf{u}}\left(\v{\bt}{\bwo} - \vH{\hbt}\right)\\ &=   
                Z_0 + \mathcal{O}_{\rm p}(\delta_n) -\dfrac{1}{H}\sum_{h=1}^{H} \big(Z_h + \mathcal{O}_{\rm p}(\delta_n)\big)  
         - \dfrac{1}{H}\sum_{h=1}^H \mathfrak{X}_{\bto,n,\mathbf{u}} \mathbf{V}\left(\btheta^{(\mathbf{v})},n,\bw_{h}\right)\bD
    \end{aligned}
    \label{proof:prop4:XV_THETA}
\end{equation}
The first three terms of \eqref{proof:prop4:XV_THETA} are considered separately in the following computation 
\begin{equation}
    \begin{aligned}
        Z_0 + \mathcal{O}_{\rm p}(\delta_n) -\dfrac{1}{H}\sum_{h=1}^{H} \big(Z_h + \mathcal{O}_{\rm p}(\delta_n)\big) &=  Z_0 - \dfrac{1}{H}\sum_{h=1}^H Z_h + \mathcal{O}_{\rm p}(\delta_n) - \dfrac{1}{H}\sum_{h=1}^H \mathcal{O}_{\rm p}(\delta_n) 
        \\
        & \overset{d}{=} \left(1+\frac{1}{H}\right)^{\nicefrac{1}{2}} Z + \mathcal{O}_{\rm p}(\delta_n) - \mathcal{O}_{\rm p}\left(\frac{\delta_n}{\sqrt{H}}\right) \\
        &= \left(1+\frac{1}{H}\right)^{\nicefrac{1}{2}} Z + \mathcal{O}_{\rm p}(\delta_n),
    \end{aligned}
    \label{proof:prop4:XVI_THETA}
\end{equation}
where $Z \sim \mathcal{N}\left(0, 1\right)$. We now consider the last term of \eqref{proof:prop4:XV_THETA}. Since $\bm{\Sigma}(\btheta,n)^{-\nicefrac{1}{2}} = \mathcal{O}(1)$ and $\mathbf{V}\left(\btheta^{(\mathbf{v})},n,\bw_{h}\right) = \mathcal{O}_p\left(n^{-\nicefrac{1}{2}}\right)$ elementwise by Assumption \ref{assum:C:ast} and $|| \mathbf{u} ||_2 = 1$ we have for all $h=1,\dots,H$
\begin{equation*}
          \mathfrak{X}_{\bto,n,\mathbf{u}} \mathbf{V}\left(\btheta^{(\mathbf{v})},n,\bw_{h}\right)\bD = \mathcal{O}_{\rm p}(p)\mathbf{1}_p^T\bD,
\end{equation*}

where $\mathbf{1}_p \in \real^p$ is a vector of ones. Since Assumptions \ref{assum:A'}, \ref{assum:B'}, \ref{assum:C:ast} and \ref{assum:D:5} imply that $\hbt$ is consistent estimator of $\bto$ (see Figure \ref{Fig_assumptions-mess}), there exists a sequence $(\delta'_n) \in \mathfrak{D}$ such that $\bD=\mathcal{O}_{\rm p}(\delta'_n)$ elementwise for all $n \in \mathbb{N}^*$. Without loss of generality, we may assume that $(\delta'_n)=(\delta_n)$. Therefore, we have 
\begin{equation*}
\begin{aligned}
\mathfrak{X}_{\bto,n,\mathbf{u}} \mathbf{V}\left(\btheta^{(\mathbf{v})},n,\bw_{h}\right)\bD 
& = 
\mathcal{O}_{\rm p}(p)\mathbf{1}_p^T\bD =  
\mathcal{O}_{\rm p}(p)\mathbf{1}_p^T\mathcal{O}_{\rm p}(\delta_n)\mathbf{1}_p \\ 
& =
\mathcal{O}_{\rm p}(p\delta_n)\mathbf{1}_p^T\mathbf{1}_p
=
\mathcal{O}_{\rm p}(p^2\delta_n),
\end{aligned}
\end{equation*}
which leads to the computation for the last term of \eqref{proof:prop4:XV_THETA},
\begin{equation}
    \dfrac{1}{H}\sum_{h=1}^H \mathfrak{X}_{\bto,n,\mathbf{u}} \mathbf{V}\left(\btheta^{(\mathbf{v})},n,\bw_{h}\right)\bD =
    \dfrac{1}{H}\sum_{h=1}^H \mathcal{O}_{\rm p}(p^2\delta_n) = \mathcal{O}_{\rm p}\left(\delta_n\frac{p^2}{\sqrt{H}}\right).
    \label{proof:prop4:XV_THETA_FINAL}
\end{equation}

Combining \eqref{proof:prop4:B_THETA}, \eqref{proof:prop4:XV_THETA} and \eqref{proof:prop4:XV_THETA_FINAL}, we obtain
\begin{equation*}
\begin{aligned}
   &\left(1+\frac{1}{H}\right)^{-\nicefrac{1}{2}}  \mathfrak{X}_{\bto,n,\mathbf{u}}\mathbf{B}\left(\btheta^{(\mathbf{a},\mathbf{r})},n\right)\bD   
   \\=
   &\left(1+\frac{1}{H}\right)^{-\nicefrac{1}{2}}  \mathfrak{X}_{\bto,n,\mathbf{u}}\left(\v{\bt}{\bwo} - \vH{\hbt}\right) 
    \\
    &\overset{d}{=}
    \left(1+\frac{1}{H}\right)^{-\nicefrac{1}{2}} \left(\left(1+\frac{1}{H}\right)^{\nicefrac{1}{2}} Z + \mathcal{O}_{\rm p}(\delta_n) + \mathcal{O}_{\rm p}\left(\delta_n\frac{p^2}{\sqrt{H}}\right)\right) 
    \\
    &\overset{d}{=}
    Z + \mathcal{O}_{\rm p}\left(\delta_n\sqrt{\frac{H}{H+1}}\right) + \mathcal{O}_{\rm p}\left(\delta_n\frac{p^2}{\sqrt{H+1}}\right)
    \\
    &=
    Z + \delta_n\mathcal{O}_{\rm p}\left(1\right) + \delta_n\mathcal{O}_{\rm p}\left(\frac{p^2}{\sqrt{H}}\right)
    =
    Z + \delta_n\mathcal{O}_{\rm p}\left(\max \left(1,\frac{p^2}{\sqrt{H}}\right)\right).
\end{aligned}
\end{equation*}
The proof is completed using \eqref{eq:assum:d:5} of Assumption \ref{assum:D:5}.
\end{proof}

\begin{Remark}
From the proof of Proposition \ref{Thm:Gauss:approx}, we have that
\begin{equation}
    \scalebox{0.99}{$
    \sqrt{n}\left(1+\frac{1}{H}\right)^{-\nicefrac{1}{2}}\mathbf{u}^T\bm{\Sigma}(\bto,n)^{-\nicefrac{1}{2}}\mathbf{B}(\bt^{(\mathbf{a}, \mathbf{r})}, n)\left(\hbt-\bto\right) \overset{d}{=} Z + \delta_n\mathcal{O}_{\rm p}\left( \max\left(1,\frac{p^2}{\sqrt{H}}\right)\right)$},
    \label{eq:gauss:approx:sm}
\end{equation}
without using \eqref{eq:assum:d:5} in Assumption \ref{assum:D:5}. As a result, Assumption \ref{assum:D:5} is not strictly necessary for the Gaussian approximation of the distribution of $\hbt$. However, while this assumption is strong, it is quite convenient to deliver a reasonable approximation of $\mathbf{B}(\bt^{(\mathbf{a}, \mathbf{r})}, n)$ in \eqref{eq:gauss:approx:sm} when $p$ is allowed to diverge, since its applicability is quite limited in practice.  

\end{Remark}
\vspace{0.25cm}
\newpage
\setcounter{equation}{0}
\renewcommand{\theequation}{H.\arabic{equation}}
\section{Proofs of the Main Results}
\label{app:main:res}

The results and assumptions presented so far can now be combined to deliver the three theorems of this article. While Theorem \ref{thm:conv:consit:iter:boot} and \ref{THM:bias} were stated in the main text, Theorem \ref{app:thm:main:3}, which is stated below, presents the asymptotic distribution of the IB-estimator.
Theorem \ref{app:thm:main:3} is based on Assumption \ref{assum:D''}, which is the union of Assumption \ref{assum:D'} and Assumption \ref{assum:D:5}. Assumption \ref{assum:D''} is provided below for clarity of exposition. 

\setcounter{Assumption}{3}
\renewcommand\theAssumption{\Alph{Assumption}$^{\ast\ast}$}
\begin{Assumption}
\label{assum:D''}
    The bias function $\mathbf{d}\left(\bm{\theta}, n\right)$ is such that:
		\begin{enumerate}
		    \item The asymptotic bias function $\mathbf{a}(\bt)$ can be written as \begin{equation*}
		        \mathbf{a}(\btheta) = \mathbf{M}\btheta + \mathbf{s}, \end{equation*} where $\mathbf{M} \in \real^{p \times p}$ with $|| \mathbf{M} ||_F < 1$ and $\mathbf{s} \in \real^{p}$.
		    
		 \item  There exists a sample size $n^* \in \mathbb{N}^*$ such that for all $n \in \mathbb{N}^*$ satisfying $n \geq n^*$ the matrix $(\mathbf{M} + \mathbf{L}(n)+\mathbf{I})^{-1}$ exists.
         
		    \item There exist real $\beta, \gamma > 0$ such that for all $\bm{\theta} \in \bm{\Theta}$ and any $j,l=1,\dots,p$, we have
		    \begin{equation*}
		    \begin{aligned}
			  \mathbf{L}_{j,l}(n) = \mathcal{O}(n^{-\beta}), \;\;\; \mathbf{r}_j\left(\bm{\theta}, n\right) = \mathcal{O}(n^{-\gamma}),\;\;\;   \lim_{n \to \infty} \; \frac{p^{\nicefrac{3}{2}}}{n^\beta} = 0 \;\;\;
		      \text{and} \;\;\; \lim_{n \to \infty} \; \frac{p^{2}}{n^\gamma} = 0.
		    \end{aligned}
		    \end{equation*}
		    
	   \item Defining $c_n \equiv \displaystyle{\max_{j=1,\dots,p}} \mathbf{c}_j(n)$ for all $n \in \mathbb{N}^*$, we require that the sequence $\left\{c_n\right\}_{n\in\mathbb{N}^*}$ is such that
            \begin{equation*}
	        \lim_{n \to \infty} \; p^{\nicefrac{1}{2}}c_n = 0.
	        \end{equation*}    
		
       \item There exists a sample size $n^\ast \in \mathbb{N}^\ast$ such that for all $n \in \mathbb{N}^*$ satisfying $n \geq n^*$ the Jacobian matrices $\mathbf{A}(\bm{\theta}) \equiv \frac{\partial}{\partial \, \btheta^T} \mathbf{a}\left(\bm{\theta}\right)$ and $\mathbf{R}(\bm{\theta}, n) \equiv \frac{\partial}{\partial \, \btheta^T} \mathbf{r}\left(\bm{\theta}, n\right)$ exist and are continuous in $\btheta\in\bm{\Theta}$. Moreover, we require that for all $\mathbf{u}\in\real^p$ with $||\mathbf{u}||_2 = 1$, there exists a sequence $(\delta_n) \in\mathfrak{D}$ such that
        \begin{equation*}
        \scalebox{0.99}{$
         \sqrt{n}\left(1+\frac{1}{H}\right)^{-\nicefrac{1}{2}}\mathbf{u}^T\bm{\Sigma}(\bto,n)^{-\nicefrac{1}{2}}\Big[\mathbf{B}(\bto,n)-\mathbf{B}\left(\btheta^{(\mathbf{a}, \mathbf{r})},n\right)\Big]\left(\hbt-\bto\right) = \mathcal{O}_{\rm p}(\delta_n)$},
        \end{equation*}
        where 
\begin{equation*}
    \mathbf{B}(\bt,n) \equiv \mathbf{I} + \mathbf{A}(\bt) + \mathbf{L}(n) + \mathbf{R}(\bt,n)
\end{equation*}
and 
\begin{equation*}
    \mathbf{B}\big(\btheta^{(\mathbf{a}, \mathbf{r})},n\big) \equiv \mathbf{I} + \mathbf{A}(\btheta^{(\mathbf{a})}) + \mathbf{L}(n) + \mathbf{R}\left(\btheta^{(\mathbf{r})},n\right).
\end{equation*}
    \end{enumerate}
\end{Assumption}
\vspace{0.25cm}

Using this new assumption, we now restate Theorem \ref{thm:conv:consit:iter:boot} and \ref{THM:bias}, and state Theorem \ref{app:thm:main:3}, followed by their proof.

\setcounter{Theorem}{0}
\begin{Theorem}
\label{app:thm:conv:consit:iter:boot}
		Under Assumptions \ref{assum:A'}, \ref{assum:B'}, \ref{assum:C'} and \ref{assum:D}, for all $H \in  \mathbb{N}^\ast$,
		
		\begin{enumerate}
		    
		    \item There exist a $n^* \in \mathbb{N}^\ast$ such that for all $n \in \mathbb{N}$ with $n \geq n^*$, $\left\{\hbt\right\} = \widehat{\bT}_{(n,H)}$, i.e. the set $\widehat{\bT}_{(n,H)}$ is a singleton.
		    
		    \item There exist a $n^* \in \mathbb{N}^\ast$ such that for all $n \in \mathbb{N}$ with $n \geq n^*$, the sequence $\left\{\hbt^{(k)}\right\}_{k \in \mathbb{N}}$ has the following limit
		    \begin{equation*}
            \lim_{k \to \infty} \; \hbt^{(k)} = \hbt.
         	\end{equation*}
         	
         	 Moreover, there exists a real $ \epsilon \in (0, \, 1)$ such that for any $k \in \mathbb{N}^\ast$ 
         	       \begin{equation*}
         	           \left\lVert\hbt^{(k)} - \hbt\right\rVert_2  =\mathcal{O}_{\rm p}({p}^{\nicefrac{1}{2}}\,\epsilon^k).
         	       \end{equation*}
         	
         	\item $\hbt$ is a consistent estimator of $\bto$, in that for all $\varepsilon > 0$ and all $\delta > 0$, there exists a sample size $n^\ast \in\mathbb{N}^\ast$ such that for all $n \in \mathbb{N}^\ast$ satisfying $n \geq n^\ast$ we have:
	\begin{equation*}
	    \Pr \left(|| \hbt - \bm{\theta}_0 ||_2 \geq \varepsilon \right) \leq   \delta.
	\end{equation*}
	\end{enumerate} 
	\end{Theorem}
\vspace{0.25cm}

\begin{Theorem}
\label{app:THM:bias}
	Under Assumptions \ref{assum:A'}, \ref{assum:B'}, \ref{assum:C'} and \ref{assum:D'}, and for all $H \in \Ns$, the IB-estimator $\hbt$ is consistent and PT-unibiased. {Therefore, the IB-estimator satisfies:}
	\begin{enumerate}
	    \item $\hbt$ is a consistent estimator of $\bto$, in that for all $\varepsilon > 0$ and all $\delta > 0$, there exists a sample size $n^\ast \in\mathbb{N}^\ast$ such for all $n \in \mathbb{N}^\ast$ satisfying $n \geq n^\ast$ we have:
	\begin{equation*}
	    \Pr \left(|| \hbt - \bm{\theta}_0 ||_2 \geq \varepsilon \right) \leq   \delta.
	\end{equation*}
	    \item $\hbt$ is PT-unbiased, that is there exists a $n^* \in \mathbb{N}^*$ such that for all $n\in\mathbb{N}^*$ satisfying $n\geq n^*$, we have 
		$\big\lVert	\mathbb{E} [\hbtheta] - \bto  \big\rVert_2 = 0$. 
	\end{enumerate}
\end{Theorem}
\vspace{0.25cm}

\begin{Theorem}
\label{app:thm:main:3}
Under Assumptions \ref{assum:A'}, \ref{assum:B'}, \ref{assum:C:ast} and \ref{assum:D''}, and for all $H \in \Ns$, the IB-estimator $\hbt$ satisfies: 
	\begin{enumerate}
	    \item $\hbt$ is a consistent estimator of $\bto$, in that for all $\varepsilon > 0$ and all $\delta > 0$, there exists a sample size $n^\ast \in\mathbb{N}^\ast$ such for all $n \in \mathbb{N}^\ast$ satisfying $n \geq n^\ast$ we have:
	\begin{equation*}
	    \Pr \left(|| \hbt - \bm{\theta}_0 ||_2 \geq \varepsilon \right) \leq   \delta.
	\end{equation*}
	    \item $\hbt$ is PT-unbiased, that is there exists a $n^* \in \mathbb{N}^*$ such that for all $n\in\mathbb{N}^*$ satisfying $n\geq n^*$, we have $\big\lVert	\mathbb{E} [\hbtheta] - \bto  \big\rVert_2 = 0$. 
	    \item For all $\mathbf{u} \in \real^p$ such that $||\mathbf{u}||_2 = 1$, there exist a sample size $n^\ast \in \mathbb{N}^\ast$ and a sequence $\left( \delta_n\right)\in\mathfrak{D}$ such that for all for all $n \in \mathbb{N}^*$ satisfying $n \geq n^*$ we have
\begin{equation*}
    \scalebox{0.99}{$
    \sqrt{n}\left(1+\frac{1}{H}\right)^{-\nicefrac{1}{2}}\mathbf{u}^T\bm{\Sigma}(\bto,n)^{-\nicefrac{1}{2}}\mathbf{B}(\bto, n)\left(\hbt-\bto\right) \overset{d}{=} Z + \delta_n\mathcal{O}_{\rm p}\left( \max\left(1,\frac{p^2}{\sqrt{H}}\right)\right)$}, 
\end{equation*}
where $Z \sim\mathcal{N}(0,1)$.
	\end{enumerate}
\end{Theorem}

\newpage
\begin{proof}
The figure below delivers a proof of Theorems \ref{app:thm:conv:consit:iter:boot}, \ref{app:THM:bias} and \ref{app:thm:main:3}.

\begin{figure}[!ht]
\centering
\begin{adjustbox}{width=.85\textwidth}
\begin{tikzpicture}
\node (rectD14) at (4.625,-2.75) [fill = blue!10, draw = blue!45, very thick,  minimum width=0.8cm, minimum height=0.8cm, rounded corners=0.2cm] {};
\node (rectAplus) at (4.625,2.75) [fill = blue!10, draw = blue!45, very thick,  minimum width=0.8cm, minimum height=0.8cm, rounded corners=0.2cm] {};
\node (rectD15) at (6.875,-2.75) [fill = orange!10, draw = orange!70, very thick, minimum width=0.8cm, minimum height=0.8cm, rounded corners=0.2cm] {};
\node (rectD12) at (2.375,-2.75) [fill = green!20, draw = green!80!black!80, very thick, minimum width=0.8cm, minimum height=0.8cm, rounded corners=0.2cm] {};

\node (rect1) at (-1,0) [myblue, draw, thick, fill = mygrey, minimum width=1cm, minimum height=4cm] {};
\node (rect2) at (1.25,0) [myblue, draw, thick, fill = mygrey, minimum width=1cm, minimum height=4cm] {};
\node (rect3) at (3.5,0) [myblue, draw, fill = mygrey, thick, minimum width=1cm, minimum height=4cm] {};
\node (rect4) at (5.75,0) [myblue, draw, fill = mygrey, thick, minimum width=1cm, minimum height=4cm] {};
\node (rect5) at (8,0) [myblue, draw, fill = mygrey, thick, minimum width=1cm, minimum height=4cm] {};
\node (rectA2) at (1.25,1.5) [fill = green!20, draw = green!80!black!80, very thick, minimum width=0.8cm, minimum height=0.8cm, rounded corners=0.2cm] {};

\node [myblue, text width=1.8cm,align=center] (thm2) at (1.25,4.1) {\scriptsize Proposition \ref{app:prop:bias}\\PT-Unbiased};
\node [myblue, text width=1.8cm,align=center] (coro1) at (-1,4.1) {\scriptsize Corollary \ref{coro:consist}\\Consistency};
\node [myblue, text width=2.2cm,align=center] (thm1) at (3.5,4.1) {\scriptsize Proposition \ref{thm:iter:boot}\\IB convergence};
\node [myblue, text width=1.8cm,align=center] (thm3) at (5.75,4.1) {\scriptsize Proposition \ref{THM:consistency}\\Consistency};
\node [myblue, text width=2.25cm,align=center] (thm4) at (8,4.1) {\scriptsize Proposition \ref{Thm:Gauss:approx}\\Asym.~normality};

\draw [myblue, -] (coro1) -- (rect1);
\draw [myblue, -] (thm2) -- (rect2);
\draw [myblue, -] (thm1) -- (rect3);
\draw [myblue, -] (thm3) -- (rect4);
\draw [myblue, -] (thm4) -- (rect5);

\node [aqua] (D) at (2.375,-2.75) {\ref{assum:D'}};
\node [mypurple] (Dast) at (4.625,-2.75) {\ref{assum:D}};
\node [mypurple] (Dstar) at (6.875,-2.75) {\ref{assum:D''}};
\node [mypurple] (Astar12) at (4.625,2.75) {\ref{assum:A'}};
\node (rectC3) at (8,-0.5) [fill = orange!10, draw = orange!70, very thick,  minimum width=0.7cm, minimum height=0.7cm, rounded corners=0.2cm] {};
\node (rectB1) at (8,0.5) [fill = orange!10, draw = orange!70, very thick,  minimum width=0.7cm, minimum height=0.7cm, rounded corners=0.2cm] {};
\node (rectC12) at (5.75,-0.5) [fill = blue!10, draw = blue!45, very thick,  minimum width=0.7cm, minimum height=0.7cm, rounded corners=0.2cm] {};
\node (rectB2) at (5.75,0.5) [fill = blue!10, draw = blue!45, very thick, minimum width=0.7cm, minimum height=0.7cm, rounded corners=0.2cm] {};
\node (rectA4) at (8,1.5) [fill = orange!10, draw = orange!70, very thick,  minimum width=0.7cm, minimum height=0.7cm, rounded corners=0.2cm] {};
\node (rectC1) at (3.5,-0.5) [fill = green!20, draw = green!80!black!80, very thick, minimum width=0.7cm, minimum height=0.7cm, rounded corners=0.2cm] {};
\node (rectB2) at (3.5,0.5) [fill = green!20, draw = green!80!black!80, very thick, minimum width=0.7cm, minimum height=0.7cm, rounded corners=0.2cm] {};
\node [myorange] (A2) at (1.25,1.5) {\ref{assum:A'}};
\node [myorange] (B1) at (1.25,0.5) {\ref{assum:B'}};
\node [myorange] (D2) at (1.25,-1.5) {\ref{assum:D:3}};
\node [myorange] (A22) at (-1,1.5) {\ref{assum:A'}};
\node [myorange] (B2) at (-1,0.5) {\ref{assum:B'}};
\node [myorange] (C2) at (-1,-0.5) {\ref{assum:C:1}};
\node [myorange] (D3) at (-1,-1.5) {\ref{assum:D:4}};
\node [myorange] (A1) at (3.5,1.5) {\ref{assum:A:1}};
\node [myorange] (B3) at (3.5,0.5) {\ref{assum:B'}};
\node [myorange] (C1) at (3.5,-0.5) {\ref{assum:C'}};
\node [myorange] (D1) at (3.5,-1.5) {\ref{assum:D:1}};
\node [myorange] (A3) at (5.75,1.5) {\ref{assum:A:2}};
\node [myorange] (B4) at (5.75,0.5) {\ref{assum:B'}};
\node [myorange] (C12) at (5.75,-0.5) {\ref{assum:C'}};
\node [myorange] (D4) at (5.75,-1.5) {\ref{assum:D:2}};
\node [myorange] (A4) at (8,1.5) {\ref{assum:A'}};
\node [myorange] (B5) at (8,0.5) {\ref{assum:B'}};
\node [myorange] (C3) at (8,-0.5) {\ref{assum:C:ast}};
\node [myorange] (D5) at (8,-1.5) {\ref{assum:D:5}};

\draw [myorange, ->] (D3) -- (D2);
\draw [myorange, ->] (C2) -- (rectC1);
\draw [myorange, ->] (A3) -- (A1);
\draw [myorange, ->] (rectC3) -- (rectC12);
\draw [myorange, ->] (D5) -- (D4);
\draw [myorange, ->] (rectAplus.south east) -- (A3);
\draw [myorange, ->] (rectAplus.south west) -- (A1);
\draw [mypurple, ->] (rectD12.north west) -- (D2);
\draw [mypurple, ->] (rectD14.north west) -- (D1);
\draw [mypurple, ->] (rectD14.north east) -- (D4);
\draw [mypurple, <-] (rectD14.west) -- (rectD12.east);
\draw [myorange, <-] (rectD12) edge[bend left = -25] (rectD15);
\draw [mypurple, ->] (rectD15.north east) -- (D5);

\draw [myblue, double, <-] (6.3,0.15) -- (7.45,0.15);
\draw [myblue, double, <-] (0.7,0.15) -- (-0.45,0.15);

\node (rectred) at (-1,-4) [fill = blue!10, draw = blue!45, thick, minimum width=0.4cm, minimum height=0.4cm, rounded corners=0.1cm] {};
\node (rectblue) at (-1,-4.5) [fill = green!20, draw = green!80!black!80, thick,  minimum width=0.4cm, minimum height=0.4cm, rounded corners=0.1cm] {};
\node (rectgreen) at (-1,-5) [fill = orange!10, draw = orange!70, thick, minimum width=0.4cm, minimum height=0.4cm, rounded corners=0.1cm] {};
\node [myblue, text width=8.8cm] (thm2) at (3.8,-4) {\scriptsize Assumptions used in Theorem \ref{app:thm:conv:consit:iter:boot}};
\node [myblue, text width=8.8cm] (thm2) at (3.8,-4.5) {\scriptsize Assumptions used in Theorem \ref{app:THM:bias}};
\node [myblue, text width=8.8cm] (thm2) at (3.8,-5) {\scriptsize Assumptions used in Theorem \ref{app:thm:main:3}};
\end{tikzpicture}
\end{adjustbox}
\caption{Illustration of the implication links ($\longrightarrow$) of the different assumptions used in Appendices \ref{app:conv:IB} to \ref{app:asymp:norm:IB} and the sufficient conditions 
\ref{assum:A'}, \ref{assum:B'}, \ref{assum:C'}, \ref{assum:C:ast}, \ref{assum:D}, \ref{assum:D'} and \ref{assum:D''},  needed for Theorems \ref{app:thm:conv:consit:iter:boot} to \ref{app:thm:main:3}. The double arrows ($\Longrightarrow$) are used to denote implication of boxes (i.e. a set of assumptions) while simple arrows are used to represent implications between assumptions. The color boxes are used to highlight the assumptions used in Theorems \ref{app:thm:conv:consit:iter:boot} to \ref{app:thm:main:3}.}
\label{Fig_assumptions-mess}
\end{figure}
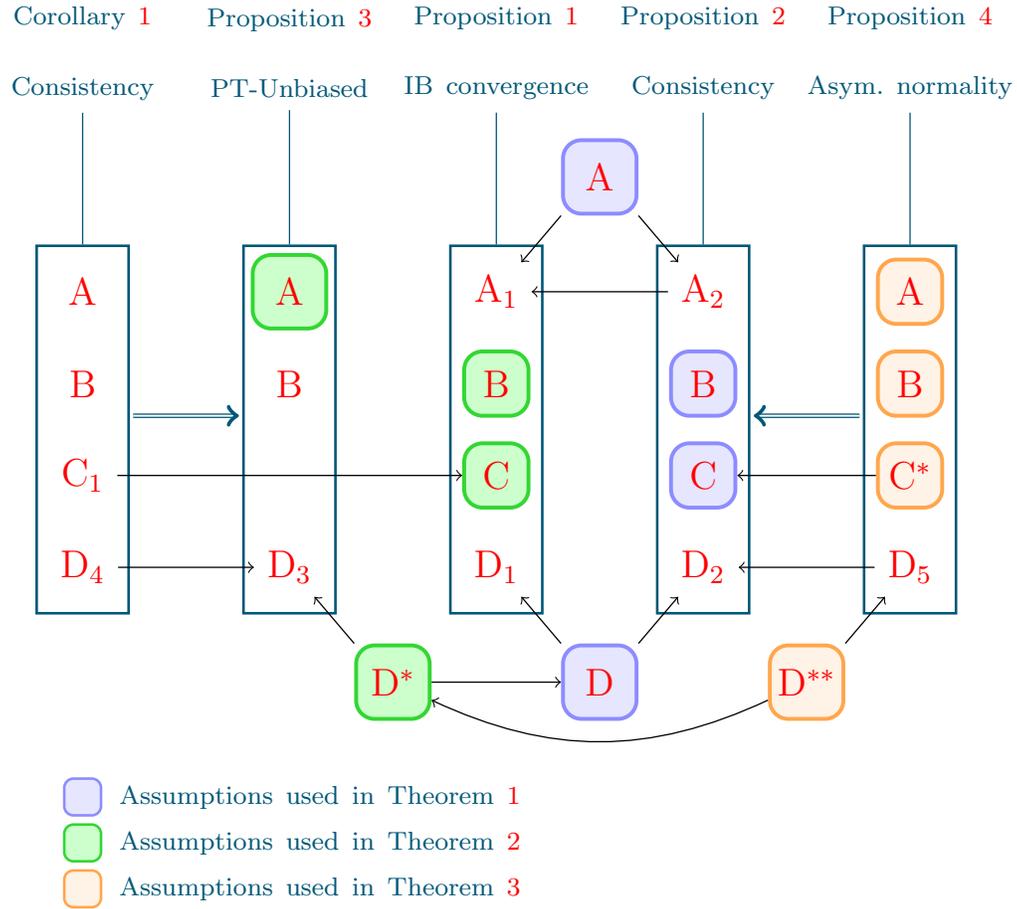
\end{proof}

\begin{Remark}
\label{rem:weaker:assum:for:thms}
Other results can be derived from Figure \ref{Fig_assumptions-mess}. For example, one can state a theorem ensuring that the IB-estimator $\hbt$ is unbiased and asymptotically normal (and therefore consistent) without the guarantee that the IB sequence $\left\{\hbt^{(k)}\right\}_{k\in\N}$ is convergent. This can be of interest if another method is used to compute the IB-estimator. Another case is when one is only interested in the convergence of the IB sequence $\left\{\hbt^{(k)}\right\}_{k\in\N}$ and the consistency of the IB-estimator $\hbt$. Assumption \ref{assum:A'} is stronger than necessary for such a result. The reason we use this assumption to prove Theorem \ref{app:thm:conv:consit:iter:boot} is guided by the rationale that we want a restricted number of general assumptions for the main Theorems \ref{app:thm:conv:consit:iter:boot}, \ref{app:THM:bias} and \ref{app:thm:main:3}.
\end{Remark}


\end{document}